\documentclass[12pt]{article}

\usepackage{graphicx} 

\usepackage[english]{babel}
\usepackage{amsfonts,amsmath, amssymb}
\usepackage{dutchcal}

\usepackage{xcolor}

\usepackage{enumerate}
\usepackage{epsf,epsfig,amsfonts,a4wide}
\usepackage{amsfonts, mathdots}
\usepackage{amsmath,amssymb,amsthm}
\usepackage{url}
\usepackage{amsmath,amssymb,amsthm}

\usepackage{amssymb}
\usepackage{amsthm}
\usepackage{epsf,epsfig,amsfonts,graphicx}

\DeclareUnicodeCharacter{1ECD}{\d{o}}

\newtheorem{Theorem}{Theorem}[section]

\newtheorem{Lemma}{Lemma}[section]
\newtheorem{Proposition}{Proposition}[section]

\newtheorem{Problem}{Problem}

\newtheorem{Definition}{Definition}[section]
\newtheorem{Ex}{Example}[section]
\newtheorem{Remark}{Remark}[section]

\theoremstyle{remark}

\newcommand{\be}{\begin{equation}}
\newcommand{\ee}{\end{equation}}

\newcommand{\R}{\mathbb{R}}\newcommand{\Id}{\textrm{\rm Id}}

\newcommand {\Lfirst} {{L_{\mathsf{comp1}}}}
\newcommand {\Lsecond} {{L_{\mathsf{comp2}}}}

\newcommand{\gl}{\mathrm{gl}}

\newcommand{\ddd}{\mathrm{d}}

\newcommand{\mM}{{\mathsf M}}

\newcommand{\dd}{{\mathrm d}\,}

\newcommand{\ii}[0]{\mathrm{i \,}}

\newcommand{\tr}{\operatorname{tr}}

\newcommand{\trace}{\operatorname{tr}}
\newcommand{\const}{\mathrm{const}}

\newcommand{\weg}[1]{}

\usepackage{csquotes}

\usepackage[
  backend=biber,
  style=numeric,
  sorting=nyt,
  isbn=false, doi=true, url=false  
]{biblatex}

\title{Research problems on relations of  Integrable Systems, Liouville Foliations, and Nijenhuis Geometry}
\author{Alexey V. Bolsinov  and 
Andrey Yu. Konyaev and 
Vladimir S. Matveev}
\date{}

\begin{document}

\author{Alexey V.\ Bolsinov\footnote{School of Mathematics,
 Loughborough University,
 LE11 3TU, UK  and Institute of Mathematics and Mathematical Modeling, Almaty, Kazakhstan\ \ \quad {\tt A.Bolsinov@lboro.ac.uk}},    
Andrey Yu.\  Konyaev\footnote{Faculty of Mechanics and Mathematics and Center for Fundamental and Applied Mathematics, Moscow State University, 119992, Moscow, Russia  
 \ \ \quad{\tt  maodzund@yandex.ru}}  \, and 
    Vladimir S.\ Matveev\footnote{
Institut f\"ur Mathematik, Friedrich Schiller Universit\"at Jena,
07737 Jena,  Germany  \ \ \quad {\tt  vladimir.matveev@uni-jena.de}}}
\maketitle

\abstract{ The paper surveys open problems and questions
related to the interplay between the theory of integrable systems with infinitely and 
finitely many degrees of freedom, Liouville foliations of finite-dimensional systems,  Nijenhuis geometry and Haantjes geometry, and integrable and superintegrable systems.
This text has grown from the  preparation for  the workshops Differential Geometry,  Integrable Systems, Liouville Foliations, and Nijenhuis Geometry, CIMAT, Guanajuato, August 2025 and Differential  Geometry and Integrable Systems, La Trobe University, Melbourne, February 2026. 
It  includes both relatively simple questions to get familiar with the topic and challenging problems that are of great importance to the area. }

\tableofcontents 

\section{Introduction}

The aim of this paper is to present a collection of open problems at the intersection of differential geometry, the theory of integrable systems, and related areas of mathematics. These problems are motivated by recent developments in Nijenhuis and Haantjes geometry, the theory of finite- and infinite-dimensional integrable systems, and the geometry of Killing tensors and superintegrable systems. Although the problems arise in different contexts and can often be studied independently, they share a common theme: the use of geometric and algebraic structures to understand integrability, its underlying mechanisms, and the behaviour of integrable systems near singularities.

The problems vary considerably in difficulty and scope. Some are concrete questions that can be approached by direct calculations or computer experiments and may serve as an introduction to the subject. Others concern fundamental geometric structures or suggest new connections between different branches of mathematics. We hope that this collection will be useful both to researchers already working in these areas and to mathematicians interested in entering the field.

One of the main lines of the suggested problems is due to a recently found relation between a  certain class of infinite-dimensional integrable systems, the so-called   BKM systems  constructed in  \cite{nijapp4}, and certain explicit finite-dimensional integrable systems of Benenti   type. 

The relation will be called {\it finite-dimensional reduction}. The parameter of this finite-dimensional reduction is $N$, which corresponds, in the KdV case,  to the number of gaps in the finite-gap solutions.

We suggest studying   \begin{itemize}
\item These finite-dimensional systems. 
  \begin{itemize}\item In particular, we suggest  studying the  Liouville foliations of these systems.  Note that the system ``lives'' on a noncompact manifold, so the leaves of Liouville foliations are not
  necessarily  compact. 
  \item In particular, we would like  to study singular points of the Liouville foliation. Note that  the most famous solutions of the KdV system (which is a special simple case of BKM systems), solitons, correspond to singular Liouville leaves.  Note also that general BKM systems  allow more freedom than the KdV system.  An interesting point could be a study of the  existence of a focus-type singularities and the corresponding monodromy effects effects~\cite{BolsinovMatveev1996, Duistermaat1980,  Matveev1996,VuNgoc2003, Zung1997}. 
\end{itemize}
  \item We suggest  understanding  how the (symplectic) topology of the 
  Liouville foliation affects the dynamics of the corresponding system.

  \item An interesting point of the study is how noncompactness of Liouville foliation affects the dynamics. Of course, it allows the point to escape to infinity, which is not very interesting. But in certain cases the point remains in a bounded region. We conjecture
  that these cases correspond to the so-called peakon solutions of the  Camassa-Holm equation, and  suggest to carefully study them, and also their analogy for other BKM systems. 
\begin{itemize}
\item In particular, we suggest  studying    the transformation which sends solutions of the  corresponding Benenti system to the solutions of the BKM systems. This transformation is explicit, but  is  given by 
a quite complicated formula which should be studied in more detail. 

\item The previous point is important  to understand whether the set of solutions coming from finite-dimensional reduction is dense in the set of all solutions.  This may lead, in perspective, to a quantization of BKM systems.
\end{itemize} 
\item  The Benenti system is closely related to Nijenhuis geometry and 
geodesically equivalent metrics. A related, but still independent direction of study is the  investigation of  singular points of gl-regular Nijenhuis operators, in  dimension  3. 
\end{itemize}

The relation between BKM systems and finite-dimensional integrable systems is only one of the themes of this paper. We also discuss a number of other problems, which are largely independent of BKM systems but are related to the same general circle of ideas. These include:

\begin{itemize} \item Haantjes geometry:
the relation between Haantjes and Nijenhuis operators, in particular in the presence of compatible symplectic structures and conservation laws.
\item Superintegrable systems: algebraic structures associated with maximal superintegrability, magnetic geodesic flows, superintegrable metrics with closed geodesics, and the relation between quadratic superintegrability and separation of variables.
\item Killing tensors and Frobenius structures: their connections with integrable systems of hydrodynamic type, Hamiltonian operators, separation of variables, and the geometry of symmetric spaces and Lie groups.
\item Singularities of Nijenhuis operators: local classification problems for gl-regular Nijenhuis operators, including the behaviour of their algebraic types near singular points.
\item Algebraic aspects of Nijenhuis geometry: the classification of left-symmetric algebras and their relation to linear Nijenhuis operators.
\item Characteristic maps:
geometric properties and singularities of characteristic maps associated with Nijenhuis operators.
\end{itemize}
The problems collected here are not intended to form a complete list of open questions in these areas. Rather, they reflect directions in which we believe that further progress is both possible and interesting. We hope that the connections between the different subjects will stimulate new approaches and collaborations.

\section{Basics of Nijenhuis Geometry and BKM systems} \label{basics}

BKM systems  were constructed within the Nijenhuis geometry  program initiated in \cite{nij1,BMMT}. Benenti systems and corresponding separation of variables are also closely related to Nijenhuis geometry. Let us give a short overview of related questions and results  of Nijenhuis geometry.  

Basic facts  on Nijenhuis geometry  are systematically presented in the series of four papers
\cite{nij1, nij2, nij3, nij4}.   The paper \cite{nij1} is introductory and, we hope, reader-friendly.  In addition to well-known classical facts,  it introduces several new concepts as well as some new terminology (e.g. related to singularities of Nijenhuis operators).  It also contains several new (sometimes quite surprising) results.  

The main objects  in Nijenhuis Geometry are  $(1,1)$-tensor fields. We also call   them {\it operator fields} or simply {\it operators}.    

\begin{Definition}\label{def:nijop}
{\rm
A \emph{Nijenhuis operator} is a $(1,1)$-tensor $L$ with vanishing Nijenhuis torsion $\mathcal{N}_L$, that is, 
$$
 \mathcal{N}_L(v, w) \overset{\mathrm{def}}{=} L^2 [v, w] + [Lv, Lw] - L [Lv, w] - L [v, Lw] =0 
$$
for arbitrary vector fields $v$ and $w$. 

A manifold $\mM$ endowed with such an operator is called a {\it Nijenhuis manifold}.   
}\end{Definition} 

Notice that vanishing of the Nijenhuis torsion is the simplest geometric condition on a $(1,1)$-tensor, and this is the reason why Nijenhuis operators naturally appear in many areas of geometry, mathematical physics and even algebra.

The ultimate goal of Nijenhuis Geometry  is to answer the following natural  questions:  

\begin{enumerate} 

\item Local description: to what form can one bring a Nijenhuis operator near almost every point by a local coordinate change?  \label{q:1} 

\item Singular points: what does it mean for a point to be generic or singular in the context of Nijenhuis  geometry?
What singularities are non-degenerate? What singularities are stable?  How do Nijenhuis operators behave near non-degenerate and stable singular points?  
 
 \item Global properties: what restrictions on 
a Nijenhuis operator are imposed by compactness of the manifold?  And conversely,  what are topological obstructions to a Nijenhuis manifold carrying a Nijenhuis operator with specific properties (e.g. with no singular points)?    \label{q:3}

\item Apply results on Nijenhuis geometry to  branches of mathematics and mathematical physics, where Nijenhuis operators naturally appear. 
\end{enumerate}

 The relation of Nijenhuis geometry to geometric  theory of finite- and infinite-dimensional  integrable systems was observed independently in different setups, and also was touched upon in the series of ``Application of Nijenhuis Geometry'' papers \cite{nijapp2, nijapp1, nijapp3, nijapp4, nijapp5} and also in the papers \cite{BKMreduction,separation}.   A reader-friendly survey of related results of this  relatively long  list of papers is in the introduction section of \cite{BKMmatrix}.

Let us give some definitions which we will use later.

\begin{Definition}\label{def:alggen}
{\rm
A point $p\in \mM$ is called {\it algebraically generic}, if the algebraic type (equivalently, Segre characteristic) of $L$ does not change in some neighbourhood $U(p)\subset \mM$.  In such a situation, we will also say that $L$ is algebraically generic at $p\in \mM$.
}\end{Definition}
  
\begin{Definition} \label{def:singpoint}
{\rm
 A point  $p\in \mM$ is called {\it singular}, if it is not algebraically generic.   
}\end{Definition}

Note that the manifold and the Nijenhuis operator are always  assumed to be as smooth as necessary, say  $C^\infty$-smooth or even analytic.  The word ``singular'' is used as an antonym to ``generic''.  We also always assume that the manifold is connected.

It follows immediately from the continuity of $L$ that algebraically generic points form an open everywhere dense subset $\mM_\circ\subseteq \mM$.

Let us emphasise that in the context of our paper,  ``singular'' is a geometric property and refers to a bifurcation of the algebraic type of $L$  at a given point $p\in \mM$.   In algebra,  the property of ``being singular or regular''  is usually understood in the context of the representation theory.

\begin{Definition}
\label{def:algregular}
{\rm
An operator  $L(p)$   (and the corresponding point $p\in \mM$) is called  {\it $\mathrm{gl}$-regular}, if its $GL(n)$-orbit $$\mathcal O(L(p))=\{ XL(p)X^{-1}~|~ X\in GL(n)\}$$ has maximal dimension, namely, $$\dim \mathcal O(L(p)) = n^2- n .$$  Equivalently,  this means that the geometric multiplicity of each eigenvalue  $\lambda_i$ of $L(p)$ (i.e., the number of Jordan $\lambda_i$-blocks in the canonical Jordan decomposition over $\mathbb C$) equals one. }\end{Definition}

One more notion is related to analytic properties of the coefficients of the characteristic polynomial $$\chi_{L(x)}(t) = \det \bigl(t\cdot\Id - L(x)\bigr)= \sum_{k=0}^n t^k \sigma_{n-k}(x) \textrm{  of } L.$$

\begin{Definition}
\label{def:nondeg}
{\rm
A point $p\in \mM$ is called  {\it differentially non-degenerate}, if the differentials $$\ddd\sigma_1(x),\dots,\ddd\sigma_n(x)$$ of the coefficients of the characteristic polynomial $\chi_{L(x)}(t)$  are linearly independent at this point.  An operator is {\it differentially non-degenerate}, if every point is differentially non-degenerate. 
}\end{Definition}

\begin{Theorem}[Theorem 4.5 of \cite{nij1}]
 If $L$ is differentially non-degenerate then there exists a local coordinate system $w_1,\dots,w_n$ such that
\begin{equation}\label{dnd1}
L = \begin{pmatrix}
-w_1  &  1  &    &    & \\ 
-w_2 & 0 & 1&  & \\  
\vdots   &  \vdots &\ddots &\ddots &  \\
-w_{n-1}  & 0 &\dots & 0 &1 \\
-w_n & 0 & \dots & 0 & 0
\end{pmatrix}.
\end{equation}
The characteristic polynomial of such $L$ is 
\begin{equation} \label{eq:M1}
t^n + w_1 t^{n-1}+ w_2 t^{n-2}+\cdots + w_n,
\end{equation}
so the coordinates in which the operator has the above form are   the coefficients of the characteristic polynomial. 
\end{Theorem}

In other coordinates the operator has a different form. A useful coordinate system for the differentially 
nondegenerate operator  is the so-called diagonal coordinate  system. In this coordinate system, the operator has the form 
\begin{equation}\label{dnd}
   L_{\mathsf{diag}} = \left( \begin{array}{ccccc}
     x_1 & 0  & 0   & \dots & 0  \\
     0 &  x_2 &  0 & \dots  & 0  \\
     \vdots & & & \vdots & \\
     0 & \dots  & 0 & x_{n-1} & 0 \\
     0 &  \dots  & \dots  & 0 & x_n \end{array}\right).  \end{equation}
The coordinates $x_1,\dots, x_n$ have a clear geometric meaning: they are eigenvalues of $L$. Note that 
this  local 
coordinate system  is  not defined at all points, as at the points where the eigenvalues bifurcate their differentials  fail to be linearly independent. It was shown in \cite[\S 4.2]{nij1} that such a diagonal coordinate system exists, for a differentially nondegenerate Nijenhuis operator,  near almost every point, provided we allow some of the coordinates to be  complex-valued.

\subsection{Short description  of BKM systems} \label{s:descriptionBKM}

BKM systems were introduced in \cite{nijapp4}. Besides this paper,  the description is available in two other  places: one is \cite{BKMreduction}, and another is the research problem paper \cite{BKMmatrix}  prepared for  another Nijenhuis geometry workshop. 
BKM systems form  a  family of  $n$-component  integrable  systems of PDEs  
constructed by the authors in \cite{nijapp4} within the Nijenhuis geometry project \cite{nij1}.  This family includes, as particular cases with appropriately chosen parameters, many well-known equations such as KdV \cite{Deift2019}, coupled KdV \cite{fordy, BM2008, BS2023}, Harry Dym~\cite{Kruskal1975}, coupled Harry Dym \cite{fordy, MB2010}, Camassa--Holm \cite{FuchssteinerFokas1981, CH93}, multicomponent Camassa--Holm \cite{ChenLiuZhang2006, HI10}, Dullin--Gottwald--Holm \cite{gdh}, and Kaup--Boussinesq equations \cite{Kaup1975,MY79}.

 BKM systems admit (sufficiently many) conservation laws and symmetries of increasing order \cite{nijapp4}. Moreover, infinitely many essentially different “finite-gap  type” 
solutions of these equations can be constructed by reduction to finite-dimensional systems that are integrable in the sense of Liouville, whose solutions can be found by quadratures \cite{BKMreduction}, and, in certain cases, exactly.     
Moreover, in \cite{KM25},  a Lax pair for the BKM systems was  introduced. The Lax pairs are given by $2\times 2$ matrices whose entries are differential operators  acting on vector-valued functions.

The family 
is constructed by  the following data:

\begin{itemize} 
\item Number of components, $n\in \mathbb N$. 
  \item Differentially-nondegenerate Nijenhuis operator
  $L$   in dimension $n$. 
  \item  Polynomial 
  $m(\mu)$ of degree    $\le n$. 
\end{itemize}

Systems in the family are parametrised by a real (optionally, complex) parameter $\lambda$, including $\lambda = \infty$.  Depending on this parameter we distinguish 4 types of systems which, for brevity,  we refer to as BKM I -- IV. All these 4 types of systems can be studied using the methods sketched above, in particular 
the finite-dimensional reduction method applies to
all four types of BKM systems. We first recall
the construction of BKM IV systems, since this
class contains several familiar integrable PDEs,
including KdV and Kaup--Boussinesq, and therefore
provides a convenient starting point for comparison
with the classical theory.

Consider $\mathbb{R}^{n}(u^1, \dots, u^n)$.
Choose a   polynomial
$$m(\mu)= m_{n-1}\mu^{n-1} + m_{n-2}\mu^{n-2} \cdots + m_{0}$$ of degree $\le n-1$    and  an explicitly given     differentially nondegenerate Nijenhuis operator $L^i_j= L(u)$.    Set 
\begin{equation} \label{eq:condabove}
 \quad  \sigma(\mu)=\det (\mu \Id - L),  \quad \mathcal L_{\zeta_0} \sigma=1, \quad   \zeta = m(L) \zeta_0.
\end{equation}
As $L$ is differentially nondegenerate, the  differentials  of the coefficients of $\sigma$ are  linearly  independent. Then, the second equation  from \eqref{eq:condabove}  is a system of $n$ linear  equations on $n$ components of $\zeta_0$ with a nondegenerate coefficient matrix. Hence, it 
uniquely determines  the vector field $\zeta_0$. Therefore,   \eqref{eq:condabove} uniquely determines  the vector field $\zeta$. The explicit form of $\zeta$ in coordinates will be given in examples below.

Next, consider the following system of $n$ PDEs: 
\begin{equation}
\label{eq:explform}
u_{t} =  \tfrac{1}{2} \left(\trace L \right)_{xxx}   \zeta + 
\left(L + \tfrac{1}{2} \trace L \cdot \Id\right) u_x.   
\end{equation}
Since  $L$ is explicitly  given in terms of $u$ and   $\trace L$  and the components of $\zeta$ are explicit functions of $u$,    \eqref{eq:explform} can be  rewritten in the ``dynamical system'' form 
\begin{equation}
\frac{\partial u_i(x,t)}{\partial t}  = V_i[u], \label{eq:dynform} \end{equation}
where the components of $V_i[u]$  are  explicit polynomials  in $$ \tfrac{\partial u_i}{\partial x},\dots, \tfrac{\partial^3 u_i}{\partial x^3}$$ with  coefficients  depending  on $u$.

    We may view each solution $u(x,t)$ 
 as a family of curves $x\mapsto u(x,t)$ with $t$ being a parameter of the family. From this viewpoint,   the equation \eqref{eq:explform}  describes the evolution of a curve with time $t$. We do not discuss the well-posedness of the Cauchy problem here.

\begin{Ex}[KdV as BKM IV]{\rm
Take  $n = 1$ and differentially non-degenerate 1-dimensional Nijenhuis operator   $L = (u)$; the polynomial  $m$ has   degree $\le 0$ so that we set $m(\mu) = m_0\in \mathbb{R}$. Then,
$$\sigma(\mu)= \mu - u \, ,  \  \zeta_0= -\frac{\partial}{\partial u}  \textrm{ and   } \zeta = - m_0 \frac{\partial}{\partial u}.$$ 
Indeed,  $$\mathcal L_{\zeta_0} \sigma(\mu) = -\tfrac{\partial }{\partial u}(\mu-u) =  1$$ 
and  $$\zeta=m(L)\zeta_0 = m_0 \Id \left(-\frac{\partial}{\partial u}\right) = - m_0 \frac{\partial}{\partial u}.$$ 

The equation \eqref{eq:explform} reads then 
\begin{eqnarray}
 u_{t} 
&    =&   -\tfrac{m_0}{2} u_{xxx}+ \tfrac{3}{2} uu_x ,    \label{eq:KdV}
\end{eqnarray}
which is one of the 
equivalent versions of the famous KdV equation.

}\end{Ex} 

 We see that
 even  for $n=1$, the construction gives something interesting.  Let us recall another physically interesting example with   $n=2$. 

\begin{Ex}[Kaup-Boussinesq as BKM IV]{\rm \label{ex:KB}

Take  $n = 2$ and differentially non-degenerate Nijenhuis operator   
$$
L =   \begin{pmatrix} u_1 &1 \\ u_2 & 0\end{pmatrix}.
$$

As the polynomial  $m$ of degree $\le 1$ we again take  the constant polynomial $$m(\mu)=m_0\in \mathbb{R}$$ of degree $0$. We have $$\sigma(\mu)= \mu^2 - u_1\mu - u_2$$ so that 
$$
\zeta_0=\begin{pmatrix}0\\-1\end{pmatrix}   \ \textrm{ and }   \ \ \ \zeta = \begin{pmatrix}0\\
-m_0\end{pmatrix}.
$$

Equation \eqref{eq:explform} is then the system of two PDEs 
\begin{equation} \label{eq:KB}
\begin{aligned}
(u_1)_t & = (u_2)_x + \tfrac{3}{2} u_1 (u_1)_x, \\
(u_2)_t & = -\tfrac{m_0}{2} (u_1)_{xxx} + u_2 (u_1)_x + \tfrac{1}{2} u_1 (u_2)_x.
\end{aligned}
\end{equation}

 This is the famous  Kaup-Boussinesq  system \cite{Kaup1975},  see e.g. \cite[eqn. (4)]{Pavlov2014} or \cite[eqn. (2.32)]{MY79}. This system is also known as {\it  the  dispersive water wave} system, see e.g. \cite{fordy}. 
}\end{Ex}

Definition and construction of BKM I, II and III systems can be found in \cite[\S 2.1]{nijapp4} and will be recalled below.   The  BKM II system  also has the ``dynamical system'' form \eqref{eq:dynform}. BKM I and BKM III systems are  the so-called {\it evolutionary systems with constraints}, see \cite[eqn. (2)]{nijapp4}.

\begin{Remark} \label{rem:1}{\rm
    BKM systems constructed in \cite{nijapp4} are more general than those considered in the present paper. 
    Informally speaking, the systems in \cite{nijapp4}  are constructed from block-diagonal Nijenhuis operators $$L=L_0\oplus L_1\oplus \dots\oplus L_N$$ where each block $L_i$, $i=1,\dots, N$ contributes to the construction with a certain natural weight $\ell_i$. In the present paper, we restrict ourselves to  the  case $L=L_0$, i.e., $N=0$.
}\end{Remark}

\subsection{BKM systems of types I, II, III and IV} \label{sec:1.3}
Here we repeat once again the construction of BKM systems, now for all four types  (see more details in \cite{nijapp4, BKMreduction}). They form a family of pairwise commuting evolutionary 
flows parametrised by $\lambda\in \R \cup\{\infty\}$.  The systems of type I or II correspond to finite $\lambda\in\R$.  Systems of type III and IV  correspond to $\lambda = \infty$. 

The main ingredients of the construction are a differentially non-degenerate Nijenhuis operator $L$ on $\R^n(u_1,\dots,u_n)$, an arbitrary polynomial  
$$
m(\mu)= m_0 + m_1 \mu + m_2 \mu^2 + \dots + m_n \mu^n
$$ 
of degree $\leq n$ with constant coefficients, and  the vector field $\zeta$ uniquely defined by the relation 
$$
    \mathcal L_\zeta \sigma(\mu, u) = m(\mu) - m_n \sigma(\mu, u),
$$
where 
$ \mathcal L_\zeta$ denotes the Lie derivative and $$\sigma(\mu, u) = \det (\mu \operatorname{Id}-L)$$
is the characteristic polynomial of  $L$.

For each  $\lambda\in\R$ \footnote{We may also consider $\lambda\in\mathbb C$}, we define 
\begin{equation}\label{eq:BKM}
\begin{aligned}
    u_{t_\lambda} & = {q}_{xxx} \Bigl(L - \lambda \operatorname{Id}\Bigr)^{-1}\zeta + {q} \Bigl(L - \lambda \operatorname{Id}\Bigr)^{-1}u_x, \\
    1 & = m(\lambda) \Bigl( {q}_{xx} {q} - \tfrac{1}{2}\, {q}_x^2\Bigr) + \sigma(\lambda, u) {q}^2.
\end{aligned}    
\end{equation}
It is an evolutionary PDE with a single differential constraint (system of type I), which can be equally understood as a system of 
two PDEs on $n+1$  unknown functions $$u_1(x,t_\lambda),..., u_n(x,t_\lambda), q(x, t_\lambda)$$ of two variables.   If $\lambda$ is a root of $m(\cdot)$, then the constraint gives 
$$q = \sigma(\lambda, u)^{-1/2}$$ and hence can be ignored. The first equation of \eqref{eq:BKM} gives then an evolutionary system of PDEs on $n$  unknown functions $u_1(x,t_\lambda),..., u_n(x,t_\lambda)$ of two variables with no constraint (system of type II):
$$
 u_{t_\lambda}  = \left(\tfrac{1}{\sqrt{ \sigma(\lambda, u)}} \right)_{xxx}\Bigl(L - \lambda \operatorname{Id}\Bigr)^{-1}\zeta + \tfrac{1}{\sqrt{ \sigma(\lambda, u)}} \Bigl(L - \lambda \operatorname{Id}\Bigr)^{-1}u_x.
$$  
In order to keep uniformity of the construction, we will not specially distinguish this case.

The case  $\lambda = \infty$ still makes  sense, and the corresponding equations will be understood as systems of type III and IV.   To obtain them from  \eqref{eq:BKM}, we replace $\lambda$ with $\lambda^{-1}$ and in the new system, after appropriate rescaling,  take the linear term in its Taylor expansion in $\lambda$ at the point $\lambda = 0$.  These operations lead to the following system (of type III):
\begin{equation}
\label{eq:BKM3}
\begin{aligned}
u_{t_\infty} &= q_{xxx} \zeta + (L + q\,\Id) u_x,\\
0 &= 2q +  m_n q_{xx} - \tr L.
\end{aligned}
\end{equation}
Here all the ingredients have the same meaning as above in \eqref{eq:BKM} and  $m_n\in \R$ is the highest coefficient of the polynomial $m(\mu)= m_n t^n + \dots$    

If $m_n = 0$ (we may interpret this as a root at infinity), then the second relation in \eqref{eq:BKM3} implies  $q = \frac{1}{2}\tr L$ and \eqref{eq:BKM3}  becomes a usual multicomponent evolutionary equation (system of type IV): 
$$
u_{t_\infty} = \tfrac{1}{2} (\tr L)_{xxx} \, \zeta + \left(L + \tfrac{1}{2} \tr L \cdot \Id\right)\, u_x.
$$
Although in terms of applications, the latter equation is very important, below we will consider it as just a special case of \eqref{eq:BKM3}.

\subsection{BKM reduction and $\mathcal R$-map.} \label{sec:Rmap}

The reduced BKM system is an integrable Hamiltonian system on $T^*\R^N$, where $N$ can be arbitrarily large. By 
$(w_1,\dots,w_N,\ p_1,\dots p_N)$ we 
denote canonical coordinates on $T^*\R^N$.
The Hamiltonian of this system is of the form
$$
H = K + V,
$$
where $$K=\frac{1}{2} \sum_{i,j=1}^N g^{ij} p_ip_j$$ is the kinetic energy defined by the contravariant metric $g$ as in \eqref{eq:gM} (see below) and 
the potential $V$ depends on the choice of an arbitrary monic polynomial $$c(\mu)=\mu^{2N+n} + c_1 \mu^{2N+n-1}+\dots$$  of degree $2N+n$ (number $N$ and polynomial $c(\mu)$ serve as parameters of the reduced BKM system).  The potential is defined from the matrix relation 
$$
f (M) = U_0 M^{N-1} + U_1 M^{N-2} + \dots + U_{N-2} M + U_{N-1} \Id, 
$$
where $M$ is the Nijenhuis operator of the form \eqref{eq:matrixM} and $f(\cdot)$ is the rational function  $f(t)=\frac{c(t)}{m(t)}$.  From this relation
we can uniquely reconstruct the coefficients $$U_0, U_1, \dots, U_{N-1}$$ as functions of $w=(w_1,\dots,w_N)$.  Then by definition we set $V= U_0$.

This system is Liouville integrable and its first integrals can be written in the form 
 \begin{equation}
\label{eq:Flambda}
F_\mu (w,\dot w) = \tfrac{1}{2} g_0 (M_\mu \dot w, \dot w) - V_\mu(w),   \quad \mu \in\R,
\end{equation}
where $M_\mu$ is defined by
\begin{equation}
\label{eq:Mmu}
M_\mu = \begin{cases}
\det(\mu\Id - M) \cdot (M-\mu\Id)^{-1}, & \mbox{if $\mu \in \R$}, \\
 M -\tr M \cdot\Id , & \mbox{if $\mu = \infty$.}
\end{cases}
\end{equation}
and the potentials $V_\mu:\R^N \to \R$ are defined as 
\begin{equation}
\label{eq:Vmu}
V_\mu =  \begin{cases} U_0\mu^{N-1} + U_1\mu^{N-2} + \dots + U_{N-1}, & \mbox{if $\mu \in \R$}, \\   
-U_1, & \mbox{if $\mu = \infty$.}
\end{cases}
\end{equation}
In the formulas related to $\mu= \infty$, we assumed $N\ge 2$. 

Notice that $M_\mu$ and $V_\mu$ are both polynomials of degree $N-1$ in $\mu$ so that alternatively we can replace $M_\mu$ and $V_\mu$ by their coefficients to get $N$ independent first integrals. Also, one can easily check that
$\mu^{N-1} F_{\frac{1}{\mu}} =  -H  - \mu F_\infty + \dots$  where {\it dots} denote higher order terms in $\mu$ and $$H=\frac{1}{2} g_0(\dot w, \dot w) + U_0(w),$$ so that $F_\infty$ indeed is a first integral of our system.

To summarize, on $T^*\R^N$ we have a Liouville integrable system with Hamiltonian
$$
H(w,p)=\frac{1}{2} g_0^{-1}(p, p) + U_0(w)
$$
and commuting integrals
$$
F_\mu (w,p) = \tfrac{1}{2} g_0^{-1} (M_\mu^*\, p, p) - V_\mu(w), \quad \mu\in\R\cup \{\infty\}.
$$

To get a solution of the BKM system with parameter $\lambda \in\R\cup \{\infty\}$, we need to find a common solution $w(t,x)$ of the reduced system for two commuting Hamiltonian flows related to $H$ and $F_\lambda$.   More precisely, if $\Phi_H^x$ and $\Phi_{F_\lambda}^t$ denote the corresponding Hamiltonian flows, then we are interested in the common orbit
$$
\bigl( w(t,x),p(t,x)\bigr) = \Phi_H^x \circ \Phi_{F_\lambda}^t (w_0, p_0),
$$
where $(w_0,p_0)$ is a certain initial condition.

After this we will send/transform  $w(t,x)$ into a solution $u(x,t)$ of the BKM system with parameter $\lambda$ by using a map
$ \mathcal R:\mathcal U\subseteq \R^N \to \R^n$ defined by a certain algebraic condition which we formally describe now, where \(\mathcal U\) is an open subset on which the construction is defined.

Recall that we have two polynomials $c(\mu)$ (of degree $2N+n$) and  $m(\mu)$ of degree $\le n$.  We also introduce the polynomial  
$$w(\mu)= \mu^N + w_1 \mu^{N-1} + \dots + w_N$$  (characteristic polynomial of $M$).  Then there exists  a unique polynomial function $\rho(\mu, w)$  (monic polynomial in $\mu$ of degree $n$ with coefficients depending on $w$) such that
the polynomial  ${\rho}(\mu,w) w^2(\mu) - c(\mu)$  is divisible by $m(\mu)$ and moreover, the ratio
$$
\frac{{\rho}(\mu,w) w^2(\mu) - c(\mu)}{m(\mu)}
$$
is a polynomial of degree $\le 2N-1$\footnote{Notice that ${\rho}(\mu,w) w^2(\mu) - c(\mu)$ has degree $\le 2N+n -1$, since  by our assumptions ${\rho}(\mu,w) w^2(\mu)$ and $c(\mu)$ are monic polynomials of degree $2N+n$ so that, if we subtract one from the other, the highest degree terms cancel out.  This explains why the degree of the ratio ${\rho}(\mu,w) w^2(\mu) - c(\mu)$ and $m(\mu)$ is required to have degree $\le 2N-1$.}.  

Using this function $$\rho(\mu, w)=\mu^n + \rho_1(w)\mu^{n-1} + \dots + \rho_{n-1}(w) \mu + \rho_n(w),$$ we now define the rational  map
$$\mathcal R :\mathcal U\subseteq \R^N \to \R^n$$ 
by setting
$$
\rho(\mu, w) = \sigma(\mu, u) \overset{\mathrm{def}}{=} \det (\mu \Id - L(u)).
$$

If as coordinates $u_1,\dots, u_n$ on $\R^n$ we choose the coefficients of the characteristic polynomial $\sigma(\mu, u)$, then the $\mathcal R$-map is defined by
$$
u_i = \rho_i(w).
$$
(Note that in the companion coordinates used in Section~\ref{sec:2.4},
where
$\sigma(\mu,u)=\mu^n-u_1\mu^{n-1}-\cdots-u_n$,
the same relation takes the form
$u_i=-\rho_i(w)$.)

Now we have all the ingredients to formulate the reduction theorem for BKM systems that establishes the relationship between  the original BKM systems and the reduced BKM system.

\begin{Theorem}[\cite{BKMreduction}]\label{t2}
Let, as above, $\Bigl(w(t,x),p(t,x)\Bigr)$ be a common solution of two commuting Hamiltonian flows generated by $H$ and $F_\lambda$: 
$$
\bigl( w(t,x),p(t,x)\bigr) = \Phi_H^x \circ \Phi_{F_\lambda}^t (w_0, p_0),
$$
where $(w_0,p_0)$ is an initial condition lying on the zero level of the integrals
\begin{equation}
\label{eq:intsurface}
\mathcal X = \{ (w,\dot w) \in T\R^N~|~ F_\mu (w,\dot w) = 0,   \quad \mbox{for all $\mu \in\R$}\} \subset T\R^N.
\end{equation}

If $\lambda\in \R$, then
$$
\begin{aligned}
u(t,x) &= \mathcal R\bigl( w\left(   t/a,x\right)\bigr),\\
q(t,x) &= \frac{1}{a} \Bigl(\lambda^N +w_1\left(  t/a,x\right)\lambda^{N-1}  + \dots +  w_N\left(  t/a,x\right)\Bigr),
\end{aligned} 
$$
where $a=\sqrt{c(\lambda)}\ne 0$, is a solution of the BKM system  \eqref{eq:BKM} (if we want the solution to be real, we need $a=\sqrt{c(\lambda)}> 0$).

Similarly, if $\lambda = \infty$, then 
$$
u(t,x) = \mathcal R\bigl(w(t, x - \tfrac{c_1}{2} t)\bigr)\quad\mbox{and}\quad q(t,x) = w_1\left(t,x- \tfrac{c_1}{2}t\right) -\tfrac{c_{1}}{2},
$$
is a solution of the BKM system  \eqref{eq:BKM3}. Here $c_1$ is the first non-trivial coefficient of the polynomial $$c (\mu) =  \mu^{2N+n} + c_1\mu^{2N+n-1} + \dots + c_{2N+n}.$$
\end{Theorem}


\subsection{Step-by-step algorithm.} \label{sec:2.4}

Here we repeat once again the above reduction construction as a step-by-step algorithm that can be used for constructing solutions of \eqref{eq:BKM}  and \eqref{eq:BKM3} by solving a certain integrable Hamiltonian ODE system with $N$ degrees of freedom.

\begin{itemize}

\item[Step 0.]   Choosing parameters of the system.  These parameters are an arbitrary polynomial $m(\mu) = m_0 + m_1\mu + \dots + m_n \mu^n$ of degree $\le n$  (in particular,  the coefficient $m_n$ can equal zero) and  a number $\lambda\in \R$ for equation \eqref{eq:BKM} or $\lambda = \infty$ for equation \eqref{eq:BKM3}.

\item[Step 1.]  Explicit form of the equation.   

As already noticed,  system \eqref{eq:BKM} has an invariant meaning and can be written in any coordinate system.   The explanation below is given in the coordinate system $u^1,\dots, u^n$, in which $L$ has the following form (called first companion form \cite{nij3}): 
\begin{equation} \label{pcf}
L = \begin{pmatrix}
     u^1 & 1 & 0 & \dots & 0  \\
     u^2 & 0 & 1 & \dots & 0  \\
     \vdots & \vdots &\ddots & \ddots & \\
     u^{n - 1} & 0 & \dots & 0 & 1 \\
     u^n & 0 & \dots & 0 & 0
\end{pmatrix}
\end{equation}

  In this setting, the other ingredients of the BKM equations   \eqref{eq:BKM}  and \eqref{eq:BKM3} (of type I and III respectively)  are as follows:
$$
\begin{aligned}
m(\mu) &= m_0 + m_1\mu + \dots + m_n \mu^n, \\
\sigma(\mu, u) &= \mu^n - u_1\mu^{n-1} - \ldots - u_{n-1} \mu \ - \, u_n,   \\
\zeta & = - \begin{pmatrix} m_n u_1 + m_{n-1} \\ m_n u_2 + m_{n-2} \\ \dots \\ m_n u_n + m_0 \end{pmatrix}\\ 
\end{aligned}
$$

\item[Step 2.]   Choosing a monic polynomial of degree $2N + n$  (this polynomial serves as a parameter of the reduced  system \eqref{eq:BKMreduced})
$$
c(\mu) = \mu^{2N+n} + c_1 \mu^{2N+n-1} + c_2 \mu^{2N+n-2} + \dots   
$$
 
\item[Step 3.]   Define  the map  $$\mathcal R : \R^N(w_1, \dots, w_N) \to \R^n(u_1,\dots, u_n)$$ which will send solutions $w(t,x)$ of the reduced system to solutions $u(t,x)$ of the original PDE system \eqref{eq:BKM}. 
Consider the polynomial 
$$
\rho (\mu, w) = \mu^n + \rho_1(w) \mu^{n-1} + \dots + \rho_{n-1}(w)\mu + \rho_n(w),
$$
whose coefficients are certain functions of $w_1, \dots, w_N$ that can be explicitly found from the condition that the ratio
$$
Q(\mu, w) = \frac{\rho (\mu, w) w^2(\mu) - c(\mu)}{m(\mu)}
$$ 
is a polynomial in $\mu$ of degree $\le 2N-1$.
In the generic case, when $m(\mu)$ has $n$ distinct roots $\lambda_1,\dots, \lambda_n$ (and these roots are explicitly given), the above condition implies 
$$\rho (\lambda_i, w) w^2(\lambda_i) - c(\lambda_i)=0,$$ that is, $\rho (\lambda_i, w) = \dfrac{c(\lambda_i)}{w^2(\lambda_i)}$. Hence,  the polynomial $\rho_i$ can be reconstructed by using the Lagrange interpolating polynomial:
$$
\rho(\mu, w) = \mu^n + \sum_{i=1}^n \left( \dfrac{c(\lambda_i)}{w^2(\lambda_i)}  - \lambda_i^n \right) \prod_{s\ne i} \frac{\mu - \lambda_s}{\lambda_i - \lambda_s} 
$$
where $w(\mu) = \mu^N + w_1 \mu^{N-1} + \dots + w_N$.   
 
In the general case, when $m(\mu)$ has multiple roots  $\lambda_1,\dots, \lambda_s$ of multiplicities $k_1, \dots, k_s$   (with $k_1 + \dots + k_s = \deg m(\mu) \le n$),   the polynomial $\rho$ can be reconstructed from the following $d=\deg m(\mu)$ linear conditions  at the roots of $m(\mu)$:
 $$
 \rho(\lambda_i) = h(\lambda_i),  \dots, \rho'(\lambda_i)=  h'(\lambda_i),\dots , \rho^{(k_i-1)}(\lambda_i)=h^{(k_i-1)}(\lambda_i), \quad i=1,\dots, s,
 $$
 where $h(\mu)=\frac{c(\mu)}{w^2(\mu)}$ and the derivative is taken with respect to $\mu$.  If $d<n$, then we may interpret $\lambda=\infty$ as an additional root of $m$  of multiplicity $n-d$ to obtain $n-d$ additional conditions of the form
 $$
 \bar \rho'(0) = \bar h'(0), \ \bar \rho''(0) = \bar h''(0), \ \dots \ , \   \bar \rho^{(n-d)}(0) = \bar h^{(n-d)}(0), 
 $$
 where 
 $$
 \begin{aligned}
 \bar\rho(\mu) &= \mu^n \rho(\tfrac{1}{\mu}) = 1 + \rho_1\mu + \rho_2\mu^2 + \dots \\
 \bar h(\mu) &= \mu^n \frac{c(\frac{1}{\mu})}{ w^2 (\frac{1}{\mu})} = \frac{1 + c_1\mu +c_2\mu^2+\dots}{ (1 + w_1\mu + w_2\mu^2 + \dots)^2}
 \end{aligned}
 $$
 Solving these conditions leads to an explicit description of the coefficients $\rho_1(w), \dots, \rho_n(w)$ of the polynomial $\rho(\mu,w)$.
 
 Of course, one can also find these coefficients by formally resolving the polynomial relation
 $$
  \rho (\mu, w) w^2(\mu) - m(\mu) Q(\mu,w) =  c(\mu)   
 $$
 which amounts to a system of  $2N + n$ non-homogeneous linear equations on $2N + n$ unknown coefficients of the polynomials
 $$
 \rho  = \mu^n + \rho_1 \mu^{n-1} + \rho_2 \mu^{n-2} + \dots\quad\mbox{and}\quad Q= q_1 \mu^{2N-1} + q_2 \mu^{2N-2} + \dots
 $$
 
 The coefficients $\rho_1(w), \dots, \rho_n(w)$ of the polynomial $\rho(\mu, w)$ so obtained allow us to define the map
 $$
 \begin{aligned}
 \mathcal R   :   \R^N(w_1, \dots, w_N) & \to \R^n(u_1,\dots, u_n),    \\
 w = (w_1,\dots,w_N) & \overset{\mathcal R}{\mapsto} \mathcal (u_1, \dots, u_n) = \bigl(-\rho_1(w), \dots, -\rho_n(w)\bigr).
 \end{aligned}
 $$

Note that the polynomial $\rho(\mu, w)$  is uniquely defined by the parameters of the above constructions, i.e., polynomials $c(\mu)$ and $m(\mu)$.  However,  one can easily notice that $\rho(\mu, w)$ remains unchanged if we replace $c(\mu)$ by 
$c(\mu) +  m(\mu) q(\mu)$, where $q(\mu)$ is an arbitrary polynomial of degree $\le 2N-1$ with constant coefficients.

\item[Step 4.]  Commuting flows (ODE systems) on $T \R^N \simeq T^*\R^N$.    

Here our goal is to resolve  (analytically or numerically)  an integrable Hamiltonian system associated  with another Nijenhuis operator $M$.  This system can be written and then solved
in various coordinate systems. We will provide formulas related to the first companion coordinate system $(w_1,\dots, w_N)$ in which $M$ takes form \eqref{eq:matrixM}.

We consider the natural system on $\R^N(w_1,\dots, w_N)$ whose kinetic energy is defined by the (contravariant) metric $g$ given by \eqref{eq:gM}
and potential  $U_0=U_0(w)$  defined from the following matrix relation
$$
c(M) \bigl(m(M)\bigr)^{-1} = U_0(w) M^{N-1} + U_1(w) M^{N-2} + \dots + U_{N-1}(w) \Id
$$
Notice that in the coordinates $w_1,\dots, w_N$, the coefficients $U_0, \dots, U_{N-1}$ of the matrix polynomial on the right hand side are exactly the elements of the last column of the matrix $c(M) \bigl(m(M)\bigr)^{-1}$ on   the left hand side. 

Next we consider the Hamiltonian system with the Hamiltonian
\begin{equation}\label{eq:BKMreduced}
H(w, p) = \frac{1}{2} g_0^{-1}(p, p) + U_0(w) =\frac{1}{2} \sum_{i, j = 1}^N (g_0)^{ij} p_i p_j + U_0(w)  
\end{equation}
with $g$ and $U_0$ defined above  (note that $(g_0)^{ij}$ are the entries of $g^{-1}$ defined by \eqref{eq:gM} so that in the Hamiltonian setting we do not need to invert this matrix).

We will also need an integral $F(w,p)$  of this  Hamiltonian system. Namely, in the case of equation \eqref{eq:BKM} we take 
$$
F_\lambda (w,p) = \frac{1}{2} g_0^{-1}(M_\lambda^* p, p) - V_\lambda (w),   
$$
where $$M_\lambda = \det (\lambda\Id - M) \cdot (M - \lambda \Id)^{-1}$$ and $$V_\lambda = U_0 \lambda^{N-1} + U_1\lambda^{N-2} + \dots + U_{N-1}.$$   Recall that $$M_\lambda^* : T^*_w \R^N \to T^*_w \R^N$$ is the dual operator to $$M_\lambda : T_w \R^N \to  T_w \R^N$$ (which, in terms of matrices, is equivalent to transposition).

Similarly, in the case of equation \eqref{eq:BKM3} we formally set $\lambda = \infty$ and take
$$
F_\infty(w,p) =  \frac{1}{2} g_0^{-1}\bigl( (M - \tr M\cdot \Id )^* p, p\bigr) + U_1 (w).
$$

The functions $H(w,p)$ and $F_\lambda(w,p)$ ($\lambda\in\R \cup \{\infty\}$) commute with respect to the canonical Poisson bracket on $T^*\R^N$  and, therefore, generate  a Hamiltonian $\R^2$-action on $T^*\R^N$.

\item[Step 5.]   Finding orbits of the Hamiltonian $\R^2$-action generated by $H$ and $F_\lambda$ ($\lambda\in\R \cup \{\infty\}$).

This step is equivalent to simultaneously solving the Hamiltonian systems generated by  $H$ and $F_\lambda$:
\begin{equation}
\label{eq:twoflows}
\left\{ \begin{aligned}  
\frac{\dd w_\alpha}{\dd x} & = \frac{\partial H}{\partial p_\alpha} \\ 
\frac{\dd p_\alpha}{\dd x} & = - \frac{\partial H}{\partial w_\alpha}
\end{aligned}\right.
\quad \mbox{and} \quad 
\left\{ \begin{aligned}  
\frac{\dd w_\alpha}{\dd t} & = \frac{\partial F_\lambda}{\partial p_\alpha} \\ 
\frac{\dd p_\alpha}{\dd t} & = -\frac{\partial F_\lambda}{\partial w_\alpha}
\end{aligned}\right.
\end{equation}

This integration can be done in quadratures because $H$ and $F_\lambda$ are, in a natural sense, included into a Liouville integrable system with $N$ degrees of freedom (see Theorem \ref{t2}).  This way leads to 
 quite complicated analytic formulas involving hyperelliptic integrals.  For some good choice of parameters  (i.e., coefficients of the polynomials $m(\mu)$ and $c(\mu)$)   these formulas simplify and one can expect to find explicit solutions in elementary functions.  

On the other hand, since we deal with a system of ODEs,  one can integrate \eqref{eq:twoflows} numerically. To that end,  we 
\begin{itemize}
\item choose an arbitrary initial condition $(w_0, p_0)$ on the common level $H_0=\cdots =H_{N-1}=0$,  
\item solve (numerically) the first system to find the integral curve
$(w(x), p(x))$ of the Hamiltonian flow of $H$   (such that $w(0)= w_0$ and $p(0)= p_0$) and  
\item taking $(w(x), p(x))$ as an initial condition,  solve  (numerically) the second system to find the integral trajectory $(w(t,x), p(t,x))$ of the Hamiltonian flow of $F_\lambda$ (such that $w(0, x) = w(x)$, $p(0,x) = p(x)$).   
\end{itemize}

\item[Step 6.]   Final step.  Finding solutions to \eqref{eq:BKM}.

Take the solution $(w(t,x), p(t,x))$  of \eqref{eq:twoflows} found in Step 5.
The corresponding solution $u(t,x)$ of \eqref{eq:BKM} is obtained from it by applying the mapping $$\mathcal R : \R^N (w_1,\dots, w_N) \to \R^n(u_1,\dots,u_n)$$ from Step 3:
$$
u(t, x ) = (u_1(t,x), \dots, u_n(t,x)) =\mathcal R\bigl( w_1(t/a,x), \dots, w_N(t/a,x)  \bigr), \quad\mbox{with } a = \sqrt{c(\lambda)},
$$ 
for equation \eqref{eq:BKM}  (Type I and II) and
$$
u(t, x ) = (u_1(t,x), \dots, u_n(t,x)) =\mathcal R\bigl( w_1(t,x-\tfrac{c_1}{2} t), \dots, w_N(t,x-\tfrac{c_1}{2} t)  \bigr).
$$ 
for equation \eqref{eq:BKM3}  (Type III and IV).

\item[Step 7.]  If necessary, we can also reconstruct the function $q(t,x)$:
$$
q(t,x) = \frac{1}{a} \Bigl(\lambda^N + w_1\left( t/a,x\right)\lambda^{N-1} + w_2\left(  t/a,x\right)\lambda^{N-2} + \dots + w_N\left( t/a,x\right)\Bigr), \quad\mbox{where } a = \sqrt{c(\lambda)}.
$$
for equation \eqref{eq:BKM}  (Type I and II) and
$$
q(t,x)=w_1\left(t,x-\frac{c_1t}{2}\right)-\frac{c_1}{2}.
$$
for equation \eqref{eq:BKM3}  (Type III and IV).

\end{itemize}

\section{Topology of IHS associated with differentially non-degenerate Nijenhuis operators and finite-gap solutions of integrable PDEs}\label{ihs}

{\bf Construction.}   Let $M$ be a Nijenhuis operator differentially non-degenerate almost everywhere.  As the model example we consider
$M$ of the form \eqref{dnd1} defined on $\R^N(w_1,\dots,w_N)$, that is,
\begin{equation}
\label{eq:matrixM}
M = \begin{pmatrix}
-w_1 & 1 & 0 & \cdots  & 0 \\
-w_2 & 0 & 1 & \ddots &   \vdots\\
\vdots & \vdots & \ddots & \ddots & 0 \\
-w_{N-1} &  0 & \dots & 0 & 1\\
-w_N & 0  & \dots & 0 & 0 
\end{pmatrix}
\end{equation}

Consider the family of operators (with $\lambda\in\R$ as a parameter)    
$$
A_\lambda = \det (\lambda\cdot\Id-M) \, (\lambda\cdot\Id-M)^{-1} = \lambda^{N-1} A_0 + \lambda^{N-2} A_1 + \dots + \lambda A_{N-2} + A_{N-1}
$$
Note that the formula for $M_\mu$ used in Section \ref{sec:Rmap}, see e.g.,  \eqref{eq:Mmu},  is different from the one above by the sign, $A_\mu= -M_\mu$.  
The operators $A_i$  (as well as $A_\lambda$) satisfy one remarkable property: the evolutionary flows 
$$
w_t = A_i(w) w_x
$$
generated by these operators pairwise commute, i.e., are symmetries of each other in the sense of the theory of evolutionary PDEs. The same terminology is used for the operators themselves:  we say that $A_i$'s are symmetries of each other (see details and discussion in \cite{nij4, nijapp5}).    In terms of tensor calculus, the definition is as follows: an operator 
$B$ is said to be a {\it symmetry} of an operator $A$ if $AB=BA$ and for any vector field $\xi$ we have:
$$
A[\xi, B\xi] + B[A\xi,\xi] - [A\xi, B\xi]=0.
$$

Next, consider the contravariant 
metric $g^{-1}=\Bigl( g^{ij}\Bigr) $   and operator $M$ on   $\mathbb{R}^N(w_1,\dots, w_N)$ :
\begin{equation} \label{eq:gM}
    g^{-1}= 
\begin{pmatrix}    
0 & \cdots & 0 & 0 & \!\!\!\! 1\\ 
0 &  \cdots & 0  &\! 1 & w_1\\ 
\vdots & \iddots &\iddots & \iddots  &w_{{2}} \\ 
0&\!\!1&w_{{1}}& \iddots &   \vdots    \\
1&w_{{1}}&w_{{2}}  &\cdots & \!\!\!  w_{{N-1}}
  \end{pmatrix}.    
\end{equation}

The metric  $g$ and operator $M$ enjoy the following useful properties  (see \cite{nijapp2}):  $g$ is flat and geodesically compatible with $M$ in the sense of  \cite[\S 6.2]{nij1} or \cite[\S 1.1]{nijapp5}. The latter property allows us to construct Poisson commuting integrals for the geodesic flow of $g$ using the formula from \cite{MT97,BM2003}: 

\begin{Proposition} \label{fact:1}  
The functions 
\begin{equation}\label{eq:formulaintegrals}
I_\lambda =  \frac{1}{2}   g^{-1}\left(\det( \lambda\,\Id-M ) ( \lambda\,\Id- M^*)^{-1} p, p\right) : T^*\mathbb{R}^{N}  \to \R
\end{equation}
or, equivalently, $I_k=\frac{1}{2}g^{-1}(A_i^* p, p)$, $i=0,\dots,N-1$, are Poisson commuting integrals of the geodesic flow of $g$.
The integrals $I_0,\dots,I_{N-1}$ are functionally independent almost everywhere.
\end{Proposition}

Thus, starting from a Nijenhuis operator $M$, we have constructed an integrable system on $T^*\R^N$ with quadratic integrals.

\begin{Problem}
\label{prob:b13}
What about the topology of this integrable system?  What are the Lagrangian fibration associated with this integrable system and its singular fibers?   Are the singularities non-degenerate in the sense of
Eliasson \cite{Eliasson1990}? Unlike the more standard situation of the geodesic flows of Riemannian metrics on closed manifolds  (e.g., geodesics on an ellipsoid),  
this is a geodesic flow on a non-compact manifold of a  generally indefinite metric.  So the corresponding Lagrangian fibration certainly admits non-compact fibers.
Are there any compact fibers (Liouville tori) here?  In dimension $N=2,3$,  we can construct the bifurcation diagram. 
What does it look like?  Notice that the whole construction is based on the Nijenhuis operator $M$ so that it is expected that the final answer to all these questions will be given in terms of $M$.
\end{Problem}

Let us now explain how one can implement non-trivial potentials to the picture.  This construction is quite general and can be applied to any geodesically compatible pair $(g, M)$, if $M$ is $\gl$-regular.
 
Consider an arbitrary symmetry $S$ of the operator $M$.  Since $S$ commutes with $M$  (i.e., $SM=MS$),  $S$ is a linear combination of powers of $M$: 
$$
S = g_1 M^{N-1} + g_2 M^{N-2} + \dots + g_N \Id.
$$
(In the notation of Section \ref{sec:Rmap}, $g_j=U_{j-1}.$)

Now consider the Hamiltonian system on $T^*\mathsf M$ with the Hamiltonian of the form
\begin{equation}
\label{eq:12}
H = K + V = \tfrac{1}{2} \, g^{-1} (p, p) + g_1(x).
\end{equation}
In other words, we consider a natural dynamical system on a manifold with the metric $g$ and potential $V=g_1(x)$. 

\begin{Proposition} [\cite{BKMreduction}, \cite{Fomenko_80}] 
\label{thm:4}
The Hamiltonian system on  $T^*\mathsf M$ with the Hamiltonian \eqref{eq:12} is Liouville integrable. As commuting first integrals of it, we can take the functions
$H_\lambda =  I_\lambda + V_\lambda$, where 
$$
I_\lambda = \tfrac{1}{2} g^{-1} (A^*_\lambda p , p ) \quad\mbox{and}\quad V_\lambda (x) = g_1\lambda^{N-1} + g_2 \lambda^{N-2} + \dots + g_N,
$$
or, equivalently  (if we take the coefficients of the expansion of $H_\lambda$ in powers of $\lambda$)
$$
H_i = I_i + g_{i+1} = \frac{1}{2} g^{-1}( A_i^* p, p) + g_{i+1}(w)
$$
\end{Proposition}

In applications, the symmetry $S$ is always of the form $S=f(M)$, where $f$ is a real-analytic function (we admit potentials with singularities, so $f$ can also be a rational function).  In what follows, we set $S=f(M)$. The following observation can be verified by direct calculations. 

\begin{Proposition} \label{fact:2}
For almost every point,   $M$ has $N$ 
different eigenvalues  which we call $q_1,\dots, q_N$. As $M$ is differentially nondegenerate, $q_i$  are 
functionally independent and give a local coordinate system. In this coordinate system, the metric $g=\Bigl( g_{ij}\Bigr)$  and operator $M$ are as  follows: 
\begin{equation}\label{eq:g_and_M}
g = \sum_{i=1}^N \left(\prod_{s=1, s\ne i}^N (q_i-q_s)  \right) \dd q_i^2  \ , \  \  M = \operatorname{ diag}(q_1,q_2,\dots, q_N) .   
\end{equation}
Moreover, the functions $H_i= I_i+ g_{i+1}$ $(i=0,\dots, N-1)$  satisfy the following ``St\"ackel'' relation  (cf.~\cite{Staeckel1893}): 
{\begin{equation}\label{eq:stackel}
\begin{pmatrix} q_1^{N-1} & q_1^{N-2}& \cdots & 1 \\ 
                               q_2^{N-1} & q_2^{N-2}& \cdots & 1 \\
                                   \vdots & & & \vdots \\ 
                                   q_N^{N-1} & q_N^{N-2}& \cdots & 1\end{pmatrix} \begin{pmatrix} H_0 \\ 
                               H_1 \\
                                   \vdots   \\ 
                                  H_{N-1}\end{pmatrix} = 
\begin{pmatrix}  \frac{1}{2}p_1^2 +f(q_1)\\ 
                                \frac{1}{2}p_2^2 + f(q_2)\\
                                   \vdots   \\ 
                                  \frac{1}{2}p_N^2 + f(q_N)\end{pmatrix},
\end{equation}
}
where $p_i$ are the momenta corresponding to the coordinates $q_i$. 
In particular,     $H_0, \dots, H_{N-1}$ Poisson commute and are functionally independent almost everywhere.  
\end{Proposition}

\begin{Remark}{\rm
   The integrable systems constructed by \eqref{eq:stackel} are sometimes called {\it Benenti } or
   {\it Benenti-St\"ackel }
   systems, see e.g.  \cite{Benenti1992, Benenti2005, BM2006, IMM}. This is a well-studied class of finite-dimensional integrable systems. Finite dimensional reductions of various integrable PDEs are related to Benenti systems, see e.g. \cite{novikov74,VBM}.}   
\end{Remark}

For finite-dimensional reductions of integrable PDEs we are interested in, the function $f(\mu)$  is given by 
\begin{equation} \label{eq:Umu}
f(\mu) =  \frac{\mu^{2N+n} + c_{2}\mu^{2N+n-2}+ c_{3}\mu^{2N+n-3}+ \cdots + c_{2N+n}}{m(\mu)} =  \frac{c(\mu)}{m(\mu)}
\end{equation}
with arbitrary constant coefficients $c_2,\dots, c_{2N+n}$, where  $m(\mu)$ is a polynomial of degree at most  $n$. Observe that the second highest coefficient of the polynomial $c(\mu)$ is zero, $c_1=0$, which is a normalization we adopted here. 

One of the main results of \cite{BKMreduction} states that every solution of the natural system with Hamiltonian \eqref{eq:12} leads  (via a certain explicit formula) to a solution of a certain integrable PDE system from a large family of multicomponent PDEs including KdV, Camassa-Holm and other classical examples of integrable PDEs.   This is the reason why the topological analysis of \eqref{eq:12} is important.

\begin{Problem}
\label{prob:b14}
Describe the singular Lagrangian fibration corresponding to the integrable system generated by \eqref{eq:12}. More specifically:

\begin{itemize}

\item  Describe the singular points and singular fibers, as well as their symplectic invariants when appropriate.  

\item  Describe the equilibrium points. Which of them are non-degenerate? Are there any focus points?  

\item Describe different types of solutions (periodic, quasiperiodic, homoclinic, etc.).  From the physical point of view, the most important are solutions with fast decay as $x\to \infty$.  In terms of the reduced system,  this means that we are interested in  the solutions $w(t)$ such that $\lim_{t\to\pm\infty} w(t)= \mathsf p$,  where $\mathsf p$ is an equilibrium point.    This leads us to an 
analysis of equilibrium points  (of rank $0$) for the reduced system.

\end{itemize}

\end{Problem}

The simplest examples are the reduced equations for $n=1$ and $N=1,2$, which correspond to the KdV equation (we assume $m=-1$)
$$
u_t = \tfrac{1}{2} u_{xxx} + \tfrac{3}{2} u u_x 
$$

\begin{Ex}{\rm
The first of them explains how to find classical travelling waves, i.e. solutions of the form $u=u(x+\tfrac{c}{2}t)$, $c\in\R$. The substitution gives 
$$
cu' = u'''+ 3 u u'  ,
$$
which is equivalent to
$$
(u'' + \tfrac{3}{2} u^2 - cu)_x = 0\quad\mbox{or}\quad u_{xx} + \tfrac{3}{2} u^2 - cu = a, \quad a=\mathrm{const}\in\R
$$
This is an equation for a point moving on a line under some cubic potential, i.e., a  one degree of freedom Hamiltonian system.
This system can be easily solved (in quadratures).
$$
\begin{aligned}
u'' + \tfrac{3}{2} u^2 - cu - a &= 0 \\
u' (u'' + \tfrac{3}{2} u^2 - cu -a) &=0 \\
\bigl(\tfrac{1}{2}(u')^2 + \tfrac{1}{2} u^3 - \tfrac{1}{2}c u^2 - au\bigr)'&=0 \\
\tfrac{1}{2}p_u^2 + \tfrac{1}{2} u^3 - \tfrac{c}{2} u^2 -au &= h = \mathrm{const}
\end{aligned}
$$
Thus, the Hamiltonian of this system is $H(p_u,u) = \frac{1}{2} p_u^2 +  \tfrac{1}{2} u^3 - \tfrac{c}{2}u^2 - au$.  
For $a=0$, the level lines of this Hamiltonian are shown in Fig. \ref{fig0}.  This is exactly the singular Lagrangian fibration we are interested in.

\begin{figure}[ht!]
\centering
\includegraphics[width=0.70\textwidth]{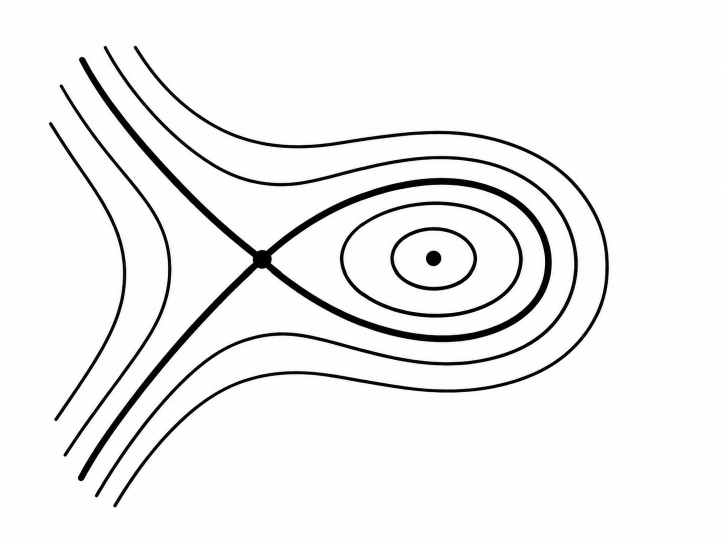}
\caption{Singular Lagrangian fibration for KdV}\label{fig0}
\end{figure}

Here we can clearly see the non-compact fiber (like a non-compact version of the atom $B$) that contains a non-degenerate saddle point. The solution we are interested in is the loop that starts and ends at this point.   In this case the integration of 
$$
p_u = u' = \sqrt{-u^3 + cu^2} 
$$

can be done explicitly. 

}\end{Ex}

\begin{Ex}{\rm
If $n=1$, $N=2$, the situation becomes more interesting.  
Now we deal with a natural system on $\R^2$.

We consider two coordinate systems:
\begin{itemize}

\item diagonal:  $x,y$
\item companion: $u=x+y, v=xy$

\end{itemize}
 
In the diagonal coordinates $x,y$, the metric is
\[
g_0=(x-y)(dx^2-dy^2)=-2\,du\,dv+u\,du^2.
\]
Let $R$ be a polynomial or rational function. For $N=2$, the
St\"ackel relation~\eqref{eq:stackel}, with separation potential $-R/2$, gives
\[
xH+F=\tfrac12\bigl(p_x^2-R(x)\bigr),\qquad
yH+F=\tfrac12\bigl(p_y^2-R(y)\bigr).
\]
Solving these relations for the Hamiltonian $H$ and its commuting
first integral $F$ gives
\begin{eqnarray}
H & =&  \frac{(p_x^2 - R(x)) - (p_y^2-R(y))}{2(x-y)}  \label{eq:Ha}\\ 
F & = & \frac{y(p_x^2 - R(x)) - x(p_y^2-R(y))}{2(y-x)} \label{eq:Fa} \end{eqnarray}
and we consider the ``Liouville torus'' $\mathcal X$  (i.e., fiber of the corresponding Lagrangian fibration) defined by
$$
\mathcal X = \{ H = 0, F=0\}.
$$ 
On the coordinate chart $x\ne y$, the separation relations show that
the fiber $\mathcal X$  is given by
\[
p_x^2=R(x),\qquad p_y^2=R(y).
\]
The Hamiltonian flows of $H$ and $F$ preserve this fiber.
More generally, on the common level $H=h$, $F=f$, the same equations
hold with $R(s)$ replaced by $R(s)+2hs+2f$. Thus the choice of zero
as the common level entails no loss of generality here.

In companion coordinates we have:
$$
(p_x \ \ p_y) = (p_u \ \ p_v) \begin{pmatrix}  u_x & u_y \\ v_x & v_y \end{pmatrix} =
(p_u \ \ p_v) \begin{pmatrix}  1 & 1 \\  y & x \end{pmatrix} = ( p_u + yp_v \ \ p_u + x p_v ). 
$$
Then, 
\[
\begin{aligned}
H={}&-p_up_v-\frac u2p_v^2
    +\frac{R(y)-R(x)}{2(x-y)}&  = \qquad  \eqref{eq:Ha} \textrm{ above},\\
F={}&\frac12\left(p_u^2+2up_up_v+
          (u^2-v)p_v^2\right) +\frac{xR(y)-yR(x)}{2(y-x)}& = \qquad  \eqref{eq:Fa} \textrm{ above}.
\end{aligned}
\]
Notice that the ``potentials'' are symmetric w.r.t. interchanging $x \leftrightarrow y$ and, therefore, can be expressed in terms of symmetric polynomials $u=x+y$ and $v=xy$.  Moreover, the numerators equal zero for $x=y$ and therefore are divisible by $y-x$ so that the potentials are smooth in $u$ and $v$.

Now we want to choose a quite special ``potential function'' $R$  (equivalently, one can say that we choose a very special separation curve or spectral curve).  Namely, we are interested in the potentials shown in Fig. \ref{fig2}.   More specifically, we assume that $R$ is a polynomial of degree 5 with two double roots.  Recall that in the general case $R$ is a rational function.  Our goal is to get a detailed description of this singular fiber $\mathcal X$. 
\begin{Problem}
\label{prob:b15}  
Describe

\begin{itemize}
\item[P1:] the topology of the fiber  (``singular Liouville torus'')
$\mathcal X = \{ H = 0, F=0\}$;

\item[P2:] the type (nondegeneracy) of the equilibrium point located on  $\mathcal X$;

\item[P3:]  the two-dimensional orbit(s) $\mathcal O$  of the natural Poisson $\R^2$-action located on $\mathcal X$ and having compact closure;

\item[P4:] one-dimensional orbits of this action;

\item[P5:] the dynamics defined by the Hamiltonians $H$ and $F$ on the two-dimensional orbit $\mathcal O$ with compact closure;

\item[P6:] projection of $\mathcal O$ onto the configuration space $\R^2(u,v)$. 

\end{itemize}
\end{Problem}

The ultimate goal is to see/observe the soliton-like dynamics.  More precisely,  the orbit $\mathcal O$ can be parametrised  by the times $x$ and $t$ corresponding to the Hamiltonian flows of $H$ and $F$:
$$
(t,x) \ \mapsto \ (u(t, x), v(t,x), p_u(t,x), p_v(t,x)) \in \mathcal O,
$$
so that for fixed $t$ we get an integral trajectory of $H$ and for fixed $x$ we get an integral trajectory of $F$.  In particular, 
due to compactness,  these flows are both complete and $(t,x)\in \R^2$ is arbitrary.
In the simplest case,  $u(t, x)$ is expected to be a soliton-like solution.  In general,   this solution is expressed in terms of $u(t, x), v(t,x)$ via a certain explicit (known) formula.

The conjectures related to the above problems P1 - P6 are as follows:

\begin{itemize}

\item[C1] The singular fiber  $\mathcal X$ is the direct product  of  two copies of the curve shown in Fig. \ref{fig1}  (this is the same curve as for the simplest KdV soliton).

\begin{figure}[ht!]
\centering

\includegraphics[width=0.5\textwidth]{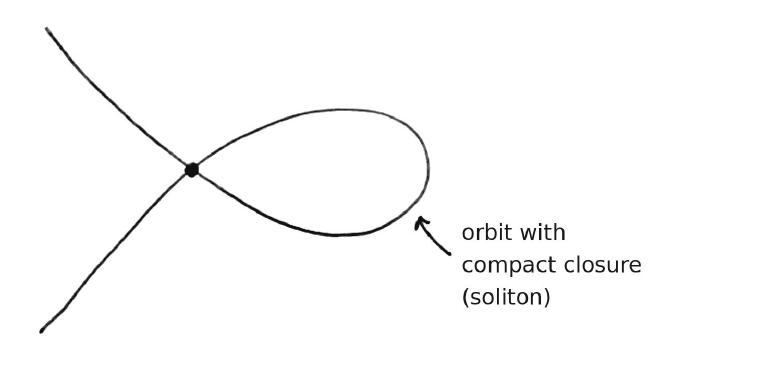}
\caption{Classical 1-soliton}\label{fig1}
\end{figure}

\item[C2]  This fiber contains exactly one equilibrium point  (point of rank zero). This point is of saddle-saddle type.  On the configuration space $\R^2(u,v)$ shown in Fig. \ref{fig3}, this point is $A$.

\item[C3]  There is one two-dimensional orbit $\mathcal O$ with compact closure. This orbit is homeomorphic to $\R^2$. The closure of $\mathcal O$ is homeomorphic to a 2-torus.
But this torus is not smooth. It should be understood as the direct product of two copies of the closure of the soliton orbit shown in Fig.\ref{fig1}  (a drop-like curve). 

\item[C4]  The closure of $\mathcal O$ consists of 4 orbits:   $\mathcal O$ itself, the equilibrium point, and two one-dimensional orbits that ``start'' and ``end'' at the singular point.  They look like 1-solitons.  The description of these singular solutions can be given in terms of the graph of $R$ shown in Fig. \ref{fig2} and eigenvalues of $M$.

Case 1. One of the eigenvalues $\lambda$ is constant and equals  $a$.   The other eigenvalue $\mu$ varies from $b$ to $c$ and then back from $c$ to $b$.  In Fig. \ref{fig3}  this corresponds to the motion along the edge  $AB$ of the triangle (from $A$ to $B$ and then back).

Case 2. One of the eigenvalues $\lambda$ is constant and equals  $b$.   The other eigenvalue $\mu$ varies from $a$ to $b$, then the eigenvalues become equal, but since $\lambda$ is constant,
the first eigenvalue remains unaffected and continues to move from $b$ to $c$,  after which it turns back and  returns to $a$ (again through $b$).    In Fig. \ref{fig3}  this corresponds to the motion along the edge  $AC$ of the triangle (from $A$ to $C$ and then back through the point $D$).

These orbits can be interesting  because the flows of $H$ and $F$ are proportional on them with a constant coefficient.  This basically means  that the common flow $u(x,t), v(x,t)$ should be given as $u(x+\alpha t), v(x+\alpha t)$.

\begin{figure}[ht!]
\centering
\includegraphics[width=0.5\textwidth]{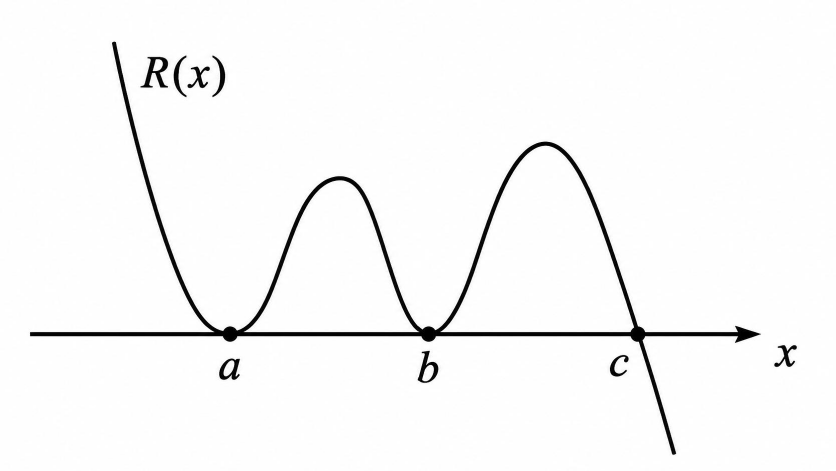}
\caption{Graph of the potential function $R(\cdot)$}\label{fig2}
\end{figure}

\item[C5]  In the topological sense, the dynamics  will be quite simple.   Typical trajectory for the Hamiltonian is shown in Fig.\ref{fig4}.  But we need more as the topology does not allow us to see the values of maxima of $u(x)$, which need to be studied.   

\item[C6]  The projection of the closure of $\mathcal O$ onto the configuration space $\R^2(u,v)$ is the triangle shown in Fig. \ref{fig3}.  It is something unusual for integrable systems.   A quite typical situation in integrable systems is that the projection is a square.  Interior points have 4 preimages,  points on the edges --- two, and vertices have just one preimage. In our case,  the situation is similar with the only exception.  Namely,  the main feature here is that the point $D$ has infinitely many preimages.   This point corresponds to collision of eigenvalues  (which happens exactly once for each solution located on $\mathcal O$).  On the (singular) torus $T^2 = \bar{\mathcal O}$ in Fig. \ref{fig4},  the preimages of $D$ are shown as the black line connecting two points $D$.      

\end{itemize}

If these conjectures are correct, then we obtain an interesting conclusion.  Since our system is invertible, together with each solution $u(x)$  this orbit will contain the solution $u(-x)$.   On the other hand these solutions must have the form 
$u(t_1,x)$ and $u(t_2, x)$.   In particular this means that if we have two solitons (two maxima for $u(x)$), then under evolution they will interchange.  Also this implies that for one of the solutions we should have $u(x) = u(-x)$. We think that this solution is characterised by the following property:   when this solution passes through $D$  (intersecting  the black line on the singular orbit $\mathcal O$), the velocity  $\dot u, \dot v$ will be zero.

\begin{figure}[ht!]
\centering
\includegraphics[width=0.9\textwidth]{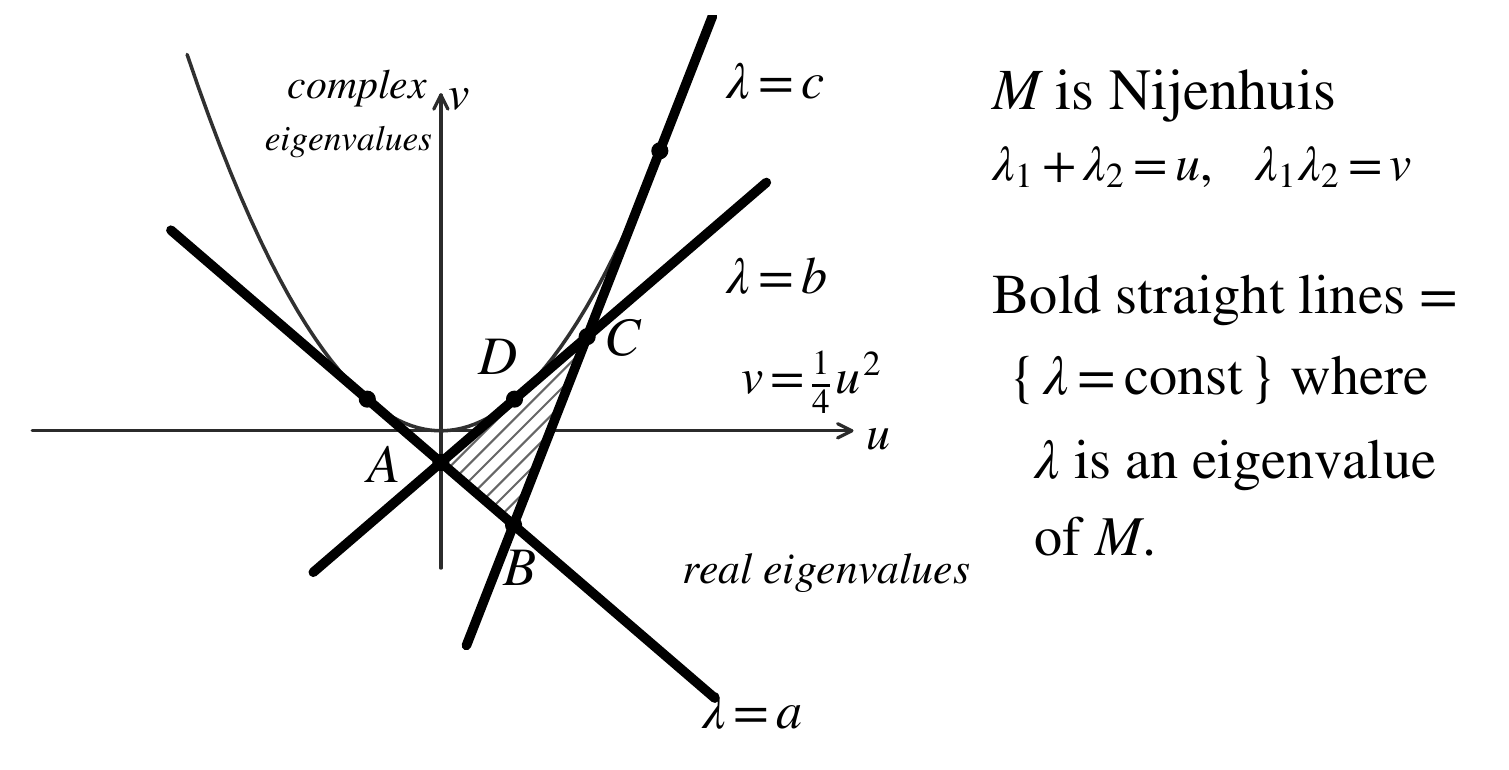}
\caption{Projection of the orbit $\mathcal O$ carrying soliton solutions onto $\R^2(u,v)$}\label{fig3}
\end{figure}

\begin{figure}[ht!]
\centering
\includegraphics[width=0.7\textwidth]{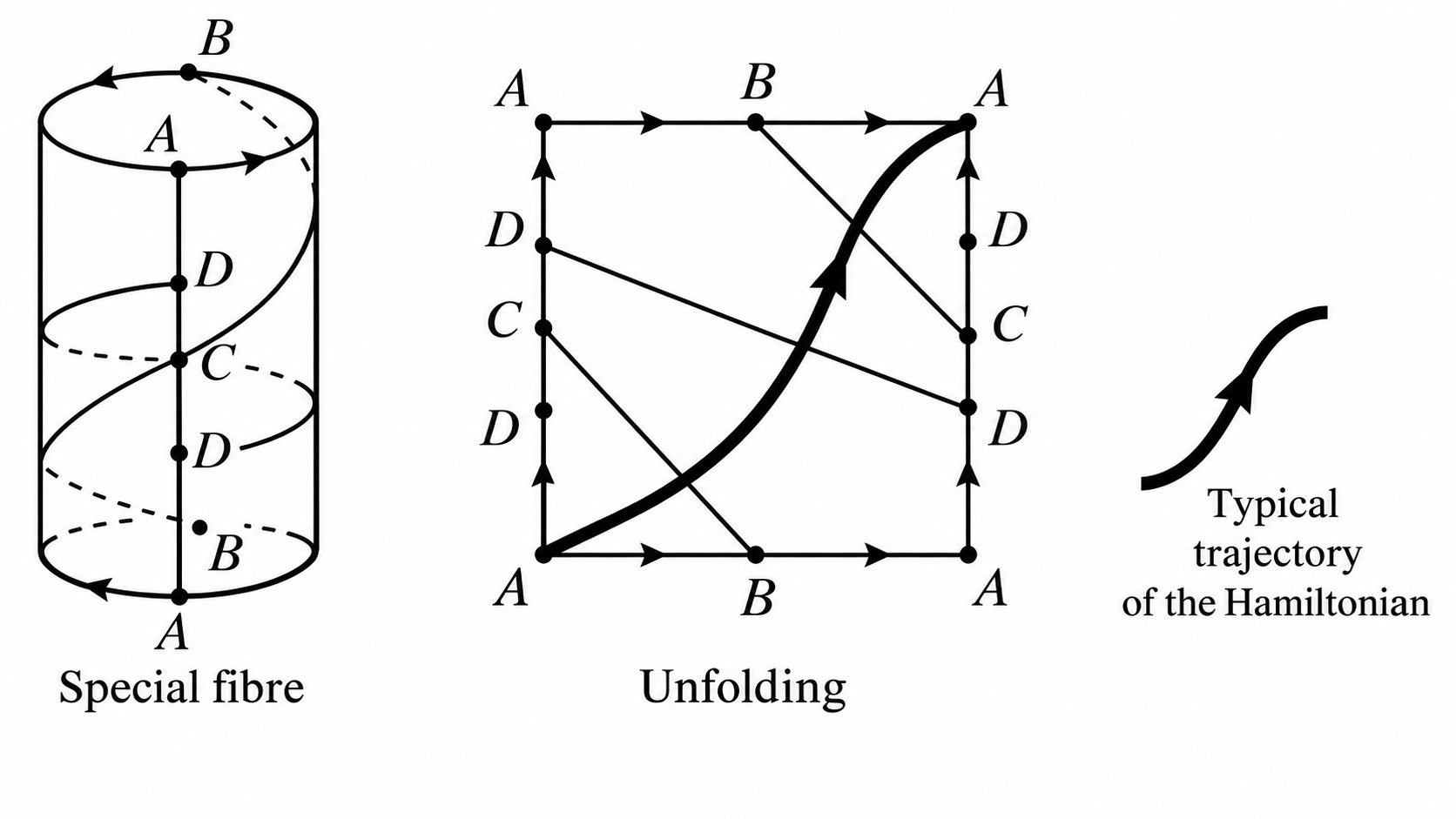}
\caption{Topology of the closure of $\mathcal O$ with its partition into 4 domains each of which projects into the triangle shown in Fig. \ref{fig3}.  }\label{fig4}
\end{figure}

}\end{Ex}

This example will give us a topological interpretation of 2-soliton solutions of KdV.

For $k$-soliton solutions we will have to set $N=k$ and consider the potential of the form
$$
R(t) = (t-a) \prod_{j=1}^k (t-a_j)^2. 
$$
Problem \ref{prob:b15} can be naturally generalised to this case.

Problem  \ref{prob:b15} is related to the KdV case which has been studied in detail  by many experts.   The new story appears if we consider the multicomponent case $n>1$, where the potential $R$ is not polynomial but rational. 

\begin{Problem}
\label{prob:b16}
How to choose the potential $R$ in such a way that the Hamiltonian system with $g$ given by \eqref{eq:gM} and potential $R$ admits  soliton-like solutions?   Study the topology of the corresponding singular fibers following the list of questions from Problem \ref{prob:b15}.  When does this system have non-degenerate saddle point(s) of rank $0$?  What is the topology of the corresponding singular fibers?
\end{Problem}


\section{Benenti type systems and infinite-dimensional integrable systems}

In this section we construct particular systems from the class of integrable systems, described in section \ref{ihs}. In all these examples $N = 2$. We change the notation and use the names $w, z$ for coordinates instead of   $w^1, w^2$.


\subsection{A pair of commuting quadratic functions with potentials}

Let $\mathbb R^2$ be a two-dimensional affine space with coordinates $w, z$. Consider  a differentially non-degenerate Nijenhuis operator $M$ in canonical coordinates:
$$
M = \left( \begin{array}{cc}
     - w  & 1 \\
     - z & 0
\end{array}\right).
$$
The operators $A_0, A_1$ are given by the formulas
\begin{equation}\label{operators}
A_0 =  \Id = \left( \begin{array}{cc}
      1 & 0  \\
     0 &  1 
\end{array}\right), \quad A_1 = M-\tr M \Id = \left( \begin{array}{cc}
     0 & 1  \\
     -z & w 
\end{array}\right).    
\end{equation}
The contravariant metric $g^{-1}$ is in the form
$$
g^{-1} = \left( \begin{array}{cc}
     0 & 1  \\
     1 & w 
\end{array}\right).
$$
  Let $h_0$ denote the Hamiltonian of the geodesic flow of $g$, and
let $f_0$ denote its quadratic first integral corresponding to $-A_1$.
In the notation of Proposition~\ref{fact:1}, $h_0=I_0$ and $f_0=-I_1$.
Explicitly,
\[
\begin{aligned}
h_0 &= \tfrac12\bigl(2p_1p_2+w p_2^2\bigr),\\
f_0 &= \tfrac12\bigl(-p_1^2-2w p_1p_2+(z-w^2)p_2^2\bigr).
\end{aligned}
\]  
It is a simple exercise to check  that the differentials   $\ddd f_0$ and $\ddd h_0$  become linearly dependent on the two-dimensional plane $p_1 = p_2 = 0$ in $T^* \mathbb R^2$. 

This implies  that this system is highly degenerate: all its singular points are stationary, that is  $\xi_{f_0} = \xi_{h_0} = 0$  for corresponding Hamiltonian vector fields. However, the system possesses  some peculiar properties. First, notice  that it admits the separation of coordinates. The separating coordinates in this case are the eigenfunctions of the Nijenhuis operator $M$:
\begin{equation}\label{separation}
\lambda_1 = \frac{1}{2} \big(- w + \sqrt{w^2 - 4 z}\big), \quad \lambda_2 = \frac{1}{2} \big(- w - \sqrt{w^2 - 4 z}\big).    
\end{equation}
However, these coordinates are real-valued for $w^2 - 4 z \geq 0$ and complex-valued for $w^2 - 4 z < 0$. 

\begin{Problem}
Describe the Liouville foliation for this system. What is the difference (if any) in Liouville foliations for two regions: $w^2 - 4 z \geq 0$ and $w^2 - 4 z < 0$ respectively?
Is there  some sort of bifurcation?   
\end{Problem}

\begin{Problem}
Using computer algebra construct the projections of particular solutions for given values of the  integrals  $h_0, f_0$  onto $\mathbb R^2$. 
\end{Problem}

We now add the potentials as in Proposition~\ref{thm:4}. For $N=2$, that
proposition gives $H_0=I_0+g_1$ and $H_1=I_1+g_2$. Thus the
Hamiltonian and first integral with potentials are
\[
h:=H_0=h_0+g_1,\qquad f:=-H_1=f_0-g_2.
\]
The potentials $g_1,g_2$ are determined by the identity
\[
S=c(M)\bigl(m(M)\bigr)^{-1}=g_1M+g_2\mathrm{Id},
\]
where $c(\lambda)$ is a monic polynomial of degree $5$ or $6$ and $m(\lambda)$ is a polynomial of degree $0, 1, 2$. We will denote the degree of $m(\lambda)$ as $d$. The degree of $c(\lambda)$ is  $ 4 + n$. Note
that the given potentials are separable as well in the same coordinates.

Notice,  that section \ref{ihs} hints  that the most interesting cases are the ones related to the singular level surfaces. The algebraic formulas that  provide the  correspondence are the same in all the cases. Thus, the behavior of the interesting solutions is given by the dynamics of the Hamiltonian vector fields on the corresponding singular surface.  


\subsection{Case $n = 1, d = 0$, KdV equation}\label{kdvsec}

This case corresponds to the case $m(\lambda) = \mathrm{const}$ and symmetry $S$ yields polynomial potential of degree $5$. W.l.o.g. we assume that $m(\lambda) = 1$ (note that in section \ref{ihs} we used $m(\lambda) = 1$ and $m(\lambda) = -1$ both). Consider separation coordinates in the form \eqref{separation}. 

\textbf{The real case.} We consider the case $  w^2 - 4 z   \geq 0$, that is the eigenvalues of $M$ are real. Notice  that taking $\lambda_1, \lambda_2$ as coordinates, we map the set $w^2 - 4z > 0$ to the set $\lambda_1 > \lambda_2$.

\begin{figure}[ht!]
\centering
\includegraphics[width=0.50\textwidth]{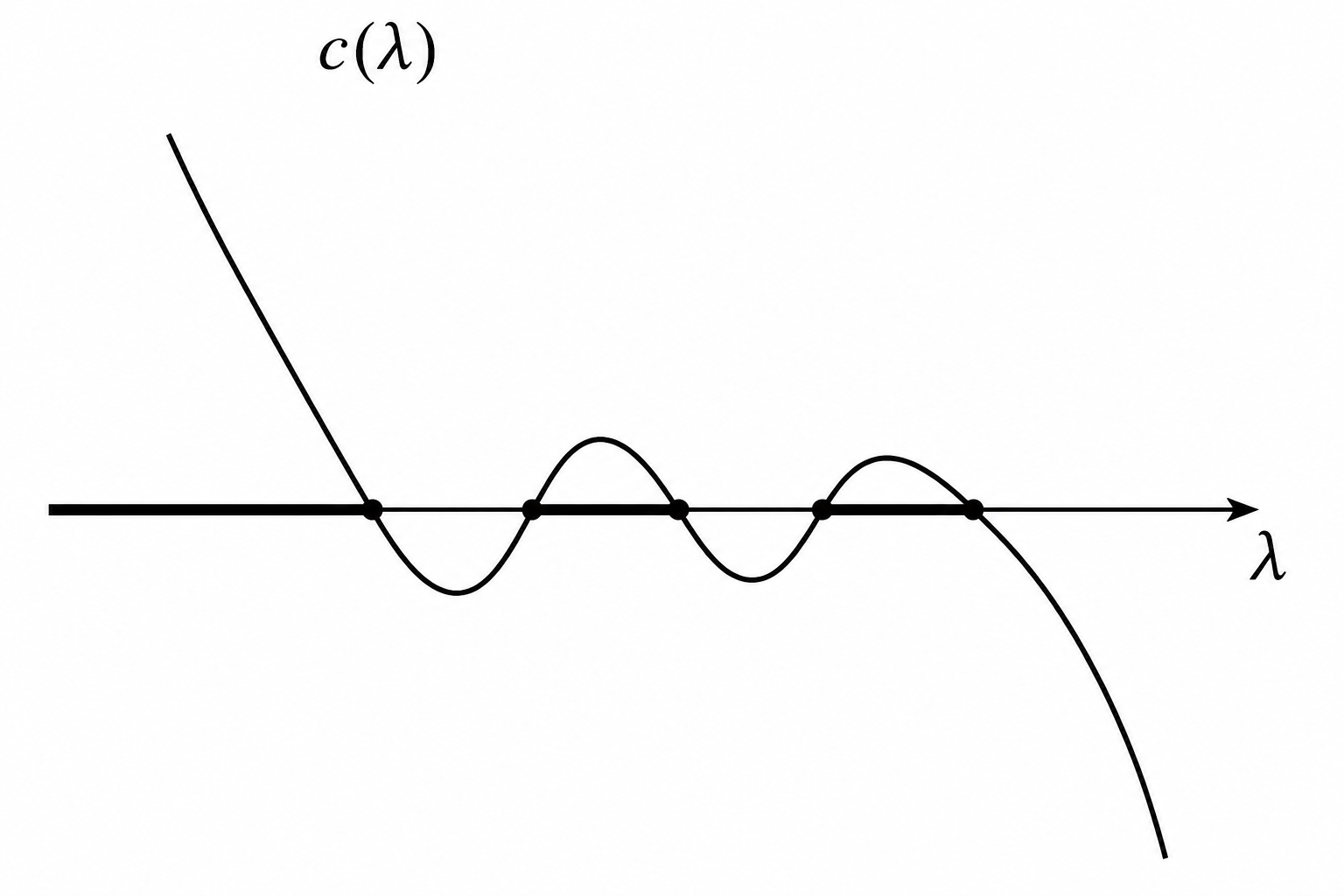}
\caption{Graph of $- c(\lambda)$. The bold sections indicate regions of admissible motion.}\label{antic}
\end{figure}

The differentials of the coordinates are
$$
\begin{aligned}
\ddd \lambda_1 = - \frac{1}{2} \ddd w + \frac{w \ddd w}{2 \sqrt{w^2 - 4z}} - \frac{\ddd z}{\sqrt{w^2 - 4z}}, \quad \ddd \lambda_2 = - \frac{1}{2} \ddd w - \frac{w \ddd w}{2 \sqrt{w^2 - 4z}} + \frac{\ddd z}{\sqrt{w^2 - 4z}}.    
\end{aligned}
$$
The matrix for the contravariant metric $g$ in these coordinates is given by $g^{-1} (\ddd \lambda_i, \ddd \lambda_j)$. At the same time the components of $f^{ij}$ are related to the components $g^{ij}$ via operator field $A_1$ from \eqref{operators}. By direct computation we get 
$$
A_1 = \left( \begin{array}{cc}
     -\lambda_2 & 0  \\
     0 & -\lambda_1 
\end{array}\right), \quad g^{ij} = \left( \begin{array}{cc}
     \frac{1}{\lambda_1 - \lambda_2} & 0  \\
     0 & \frac{1}{\lambda_2 - \lambda_1}
\end{array}\right), \quad f^{ij} = \left( \begin{array}{cc}
      \frac{\lambda_2}{\lambda_1 - \lambda_2} & 0  \\
     0 &   \frac{\lambda_1}{\lambda_2 - \lambda_1}
\end{array}\right)
$$
Now let us calculate the potential. We assume that $c(\lambda)$ has five distinct real eigenvalues. We draw $ - c(\lambda)$ because we are interested in the zones where $- c(\lambda) > 0$ (see Fig.  \ref{antic}).

\begin{figure}[ht!]
\centering
\includegraphics[width=0.50\textwidth]{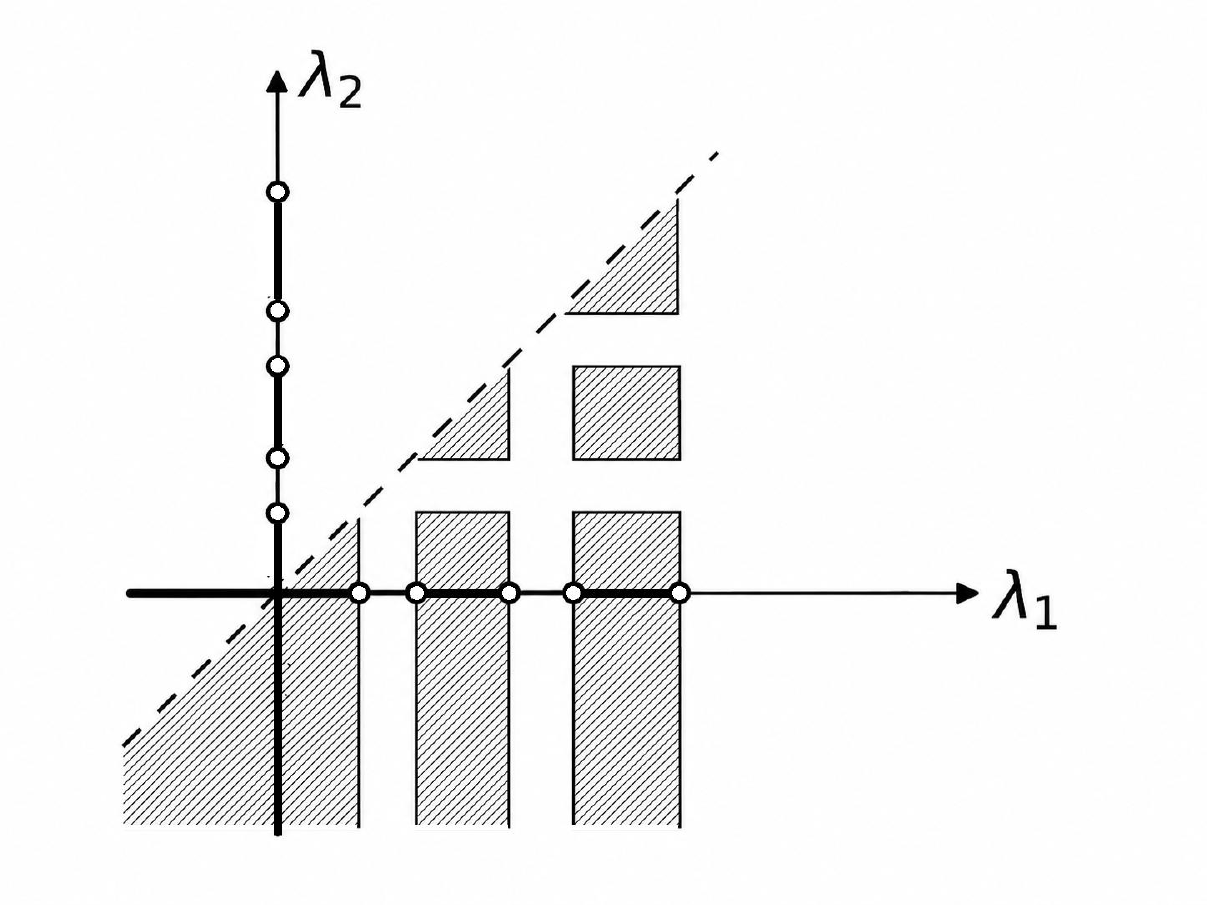}
\caption{$\lambda_1 > \lambda_2$, diagonal $\lambda_1 = \lambda_2$ is the image of the parabola $w^2 - 4z = 0$. The gray  zones are the images of the projections of Liouville surfaces (admissible motion of points). Some projections are compact.}\label{realdiag}
\end{figure}

We have
$$
S = \left( \begin{array}{cc}
     c(\lambda_1) & 0  \\
     0 & c(\lambda_2) 
\end{array}\right) = g_1 \left( \begin{array}{cc}
     \lambda_1 & 0  \\
     0 & \lambda_2 
\end{array}\right) + g_2 \left( \begin{array}{cc}
     1 & 0  \\
     0 & 1 
\end{array}\right).
$$
Solving this system with respect to $g_1, g_2$ we get
$$
g_1 = \frac{c(\lambda_1) - c(\lambda_2)}{\lambda_1 - \lambda_2}, \quad g_2 = \frac{- \lambda_2 c(\lambda_1) + \lambda_1 c(\lambda_2)}{\lambda_1 - \lambda_2}. 
$$
Now pick the level surface as (we assume now that $f, h$ stand for constants, not the functions)
$$
\begin{aligned}
h & = \tfrac{1}{2} \frac{p_1^2 - p_2^2}{\lambda_1 - \lambda_2}  + \frac{c(\lambda_1) - c(\lambda_2)}{\lambda_1 - \lambda_2}, \\
f & = \tfrac{1}{2} \frac{p_1^2 \lambda_2  - p_2^2\lambda_1}{\lambda_1 - \lambda_2} - \frac{- \lambda_2 c(\lambda_1) + \lambda_1 c(\lambda_2)}{\lambda_1 - \lambda_2}, \\
\end{aligned}
$$
This is a linear system on $p_1^2, p_2^2$. Solving it with respect to $p_1^2, p_2^2$, we get
$$
p_1^2 = 2(\lambda_1 h- f - c(\lambda_1)), \quad p_2^2 = 2(\lambda_2 h - f - c(\lambda_2)).
$$
  Along the Hamiltonian flow of $h$, we have
$\dot{\lambda}_i=\partial h/\partial p_i$. Substituting the
expressions for $p_1^2$ and $p_2^2$ gives the following equations
for the squared velocities on the chosen common level: 
$$
(\dot \lambda_1)^2 = 2\frac{\lambda_1 h - f - c(\lambda_1)}{(\lambda_1 - \lambda_2)^2}, \quad (\dot \lambda_2)^2 =2 \frac{\lambda_2 h-f - c(\lambda_2)}{(\lambda_1 - \lambda_2)^2}. 
$$
  In the real region $\lambda_1>\lambda_2$, both right-hand sides
must be nonnegative. These inequalities determine the regions of
admissible motion shown in Fig.~\ref{realdiag}.  

We are interested in the following problems.

\begin{Problem} \label{problem7}
What are the conditions on $f, h$, that yield singular surface levels? Is it true, that they correspond to the cases, when some of the roots of the polynomial $\lambda h - f - c(\lambda)$ coincide?
\end{Problem}

\begin{Problem}
What are the bifurcations, which correspond to the coinciding roots? How are  tori/cylinders glued together?
\end{Problem}

The next problem deals with the border of our coordinate system.

\begin{Problem}
What happens if a trajectory lies in the triangular zone? In particular, what happens if it hits the diagonal? In other words: we see the blow-up, but it might be an  artifact of the coordinate system, 
what happens to the  foliation?
\end{Problem}

To answer the previous question, some computer simulation might help. It should be done in the  original coordinates (they were discussed earlier, see Fig.  \ref{fig4})


\textbf{The complex case. }The main question here is rather general: what happens in the case $  w^2  - 4 z  < 0$? The separation coordinates become complex-valued, yet conjugate.

\begin{Problem}
For fixed values of the integrals $h$ and $f$, compute and plot the
projections of the connected components of their common level set
onto the $(w,z)$-plane in the region $w^2-4z<0$. Compare these
projections with those obtained for the same values of the integrals
in the region $w^2-4z>0$.
\end{Problem}
To make this comparison in the $(w,z)$-plane, use
$w=-(\lambda_1+\lambda_2)$ and $z=\lambda_1\lambda_2$ to transform
the regions shown in Fig.~\ref{realdiag}. For the numerical computations one may
also use the flat coordinates $\mu_1=w$, $\mu_2=z-\tfrac14 w^2$,
in which the boundary between the two regions is the line $\mu_2=0$.

\begin{Problem}
Determine the critical points of the momentum map $(h,f)$ lying
over the region $w^2-4z<0$, and find their critical values.
Can these values be characterized in terms of multiple nonreal
roots of the separation polynomial
\[
P(\lambda)=\lambda h-f-c(\lambda)?
\]
\end{Problem}

\begin{Problem}
Describe the topology of the singular fibers containing the critical
points found in the preceding problem. How do nearby regular fibers
degenerate as their values under the momentum map approach the
corresponding critical value?
\end{Problem}

There are two commuting vector fields on the two-dimensional Liouville surface. Solving the corresponding systems together, we get the function $w(t, x)$ (we denote the  times by  $x$ and $t$). Setting  
$$
u(t, x) =  2 w(t, x - \frac{c_1}{2}t)-c_1 , 
$$ where $c_1$ is the second coefficient of the monic polynomial $c(\lambda)$, 
we get that the corresponding function satisfies the KdV equation
$$
u_t = -\frac{1}{2} u_{xxx} + \frac{3}{2} u u_x.
$$
The answer to the next problem might be just given by computer algebra.

\begin{Problem}
What solutions of KdV are given by regular solutions (in the complex and in the real region) of the geodesic flows? What  solutions  are given by flows on the singular surfaces (in the complex and in the real region)?  
\end{Problem}

  To compare the KdV solutions arising from regular and singular fibers,
one may plot representative spatial profiles $x\mapsto u(t,x)$ at
several fixed values of $t$. 


\subsection{Case $n = 1, d = 1$, Harry Dym and Camassa-Holm equations, peakons} \label{s:HD}

\begin{figure}[ht!]
\centering
\includegraphics[width=0.50\textwidth]{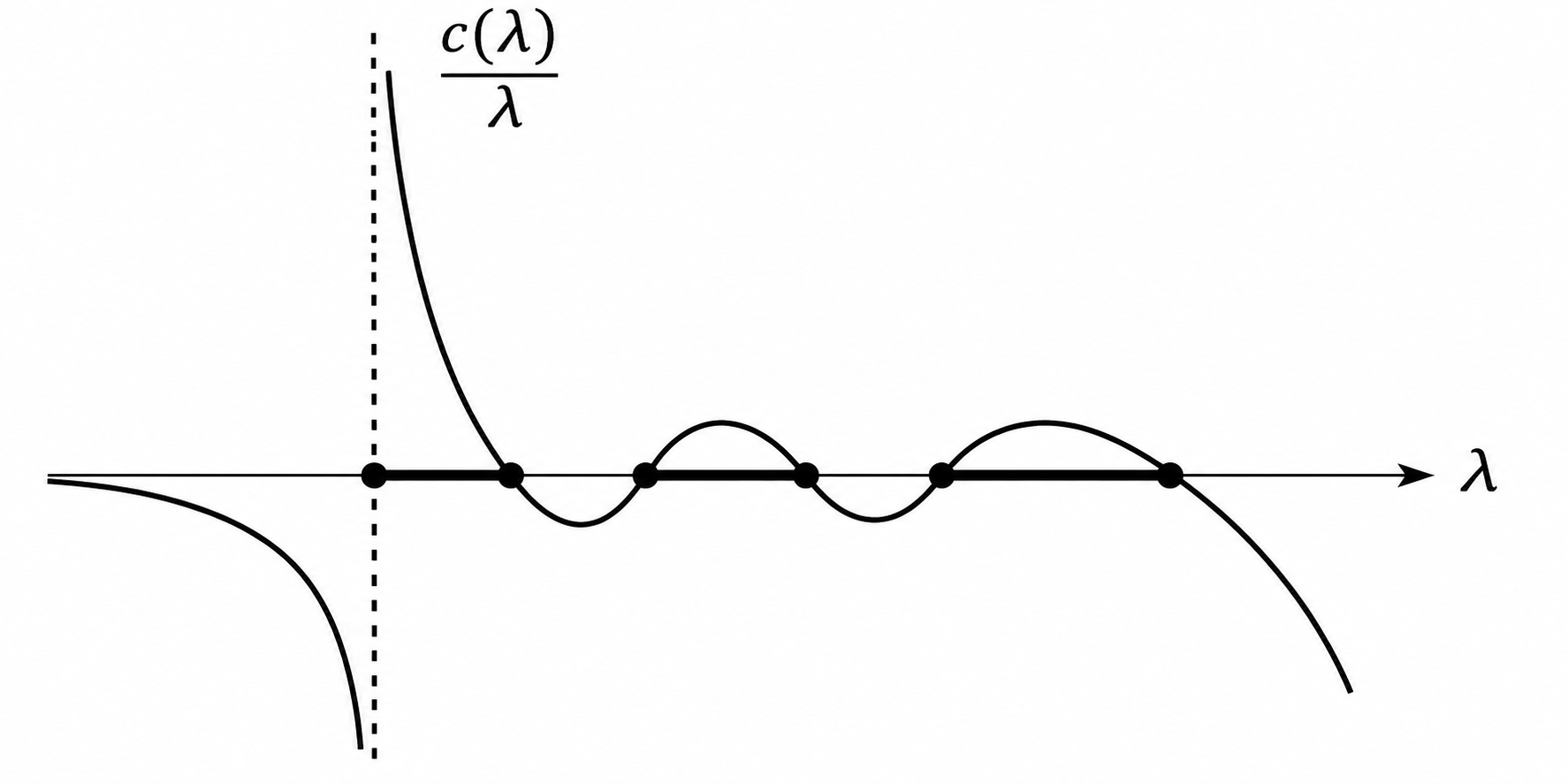}
\caption{Potential energy  for $N=2$, $m(\lambda)=m_1\lambda$. The asymptote is $\lambda = 0$. The bold regions are admissible regions of motion. } \label{fig:p}
\end{figure}

This case corresponds to the case $m(\lambda) = m_1 \lambda$ and symmetry $S$ yields polynomial potential of degree $5$. We take  
$m(\lambda) = \lambda$. The main difference from the previous case is the fact that the potential is not polynomial, but rational, that is
$$
S = c(M) M^{-1} = g_1 M + g_2 \Id.
$$
As before we assume that $c(\lambda)$ has five distinct real eigenvalues and $c(0) > 0$. A possible  diagram  of the corresponding  potential is on  Fig. \ref{fig:p}. We again start with a real case.

\begin{figure}[ht!]
\centering
\includegraphics[width=0.50\textwidth]{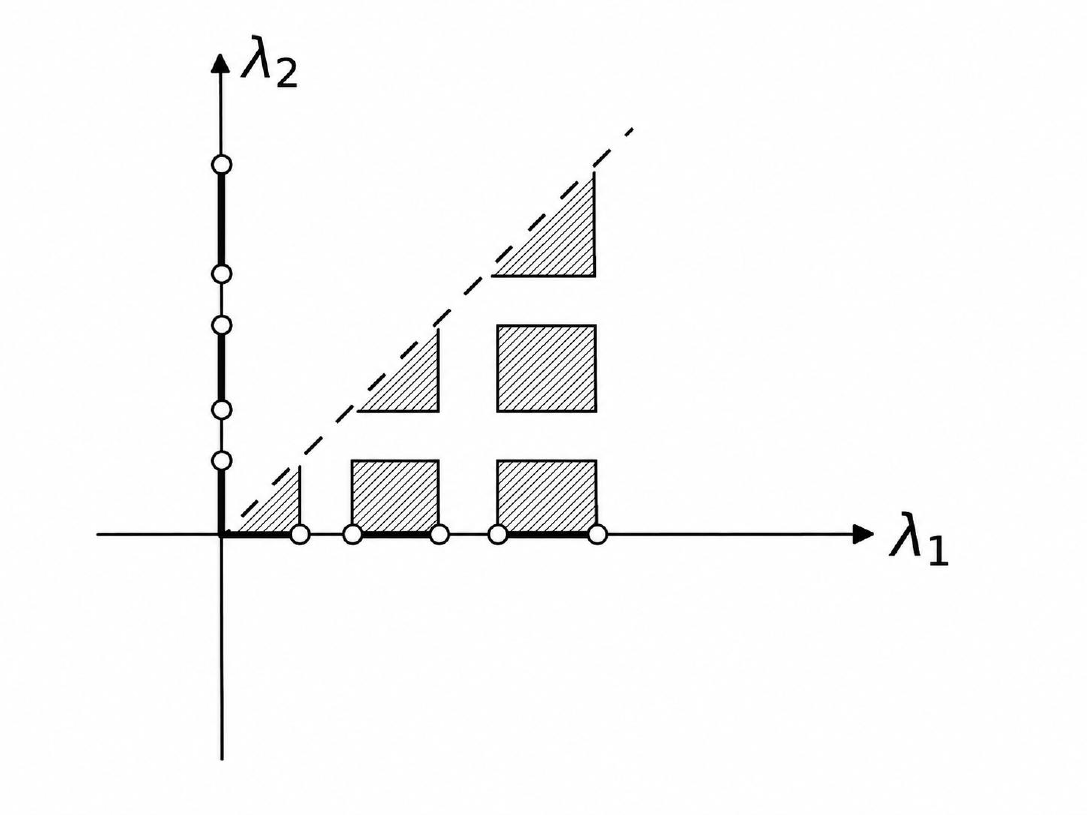}
\caption{$\lambda_1 > \lambda_2$, diagonal $\lambda_1 = \lambda_2$ is the image of the parabola $w^2 - 4z = 0$. The gray  regions are the admissible zones for the motion of points.} \label{fig:adm}
\end{figure}

\textbf{The real case. } In separation coordinates we have:
$$
S = \left( \begin{array}{cc}
     \frac{c(\lambda_1)}{\lambda_1} & 0  \\
     0 & \frac{c(\lambda_2)}{\lambda_2} 
\end{array}\right) = g_1 \left( \begin{array}{cc}
     \lambda_1 & 0  \\
     0 & \lambda_2 
\end{array}\right) + g_2 \left( \begin{array}{cc}
     1 & 0  \\
     0 & 1 
\end{array}\right).
$$
Solving this system with respect to $g_1, g_2$ we get
$$
g_1 = \frac{\lambda_2 c(\lambda_1) - \lambda_1 c(\lambda_2)}{\lambda_1 \lambda_2(\lambda_1 - \lambda_2)}, \quad g_2 = \frac{- \lambda^2_2 c(\lambda_1) + \lambda^2_1 c(\lambda_2)}{\lambda_1 \lambda_2(\lambda_1 - \lambda_2)}. 
$$
  Substituting the separation relations into Hamilton's equations,
as in Section~\ref{kdvsec} immediately before Problem~\ref{problem7}, gives for
$m(\lambda)=\lambda$  
$$
(\dot\lambda_i)^2=
\frac{2(\lambda_i^2h-\lambda_i f-c(\lambda_i))}
{\lambda_i(\lambda_1-\lambda_2)^2},
\qquad i=1,2. 
$$
The corresponding admissible regions are indicated by the bold line on the Fig. \ref{fig:adm}.

\begin{Problem}
What are the conditions on $f, h$ that yield singular surface levels? Is it true that they correspond to the cases, when some of the roots of the polynomial $\lambda^2 h - \lambda f - c(\lambda)$ coincide?
\end{Problem}

\begin{Problem}
What are the bifurcations, which correspond to the coinciding roots? How are tori/cylinders glued together?
\end{Problem}

\begin{Problem}
What happens  when the solution hits the border $\lambda_2 = 0$? The conjecture is that it reflects, like in billiards.  
\end{Problem}

\textbf{Complex case. }The problems here are the same as before.

\begin{Problem}
Determine the critical points of the momentum map $(h,f)$ lying
over the region $w^2-4z<0$, and find their critical values.
Can these values be characterized in terms of multiple nonreal
roots of the separation polynomial
\[
P(\lambda)=\lambda^2 h-\lambda f-c(\lambda)?
\]
\end{Problem}

\begin{Problem}
Describe the topology of the singular fibers containing the critical
points found in the preceding problem. How do nearby regular fibers
degenerate as their values under the momentum map approach the
corresponding critical value? 
\end{Problem}
The potential has a pole along $z=0$. Since $w^2-4z=w^2\geq 0$
on this line, the pole lies outside the region considered in these
two problems. 

We proceed  as before: we solve the pair of commuting systems to obtain functions $w(t, x), z(t, x)$. Following the general construction, we get that the function 
$$
u(t, x) = - \frac{c(0)}{z^2(t/\sqrt{c(0)}, x)}
$$
satisfies 
the reduction of the  coupled Harry Dym equation (assuming the positive square root of $c(0)$)
$$
u_t = -\Bigg( \frac{1}{\sqrt{-u}}\Bigg)_{xxx} - \frac{u_x}{-u^{\frac{3}{2}}}.
$$
The difference from the original Harry Dym equation is in the existence of the second summand. At the same time, the  functions
$$
q(t, x) = w\left(t, x- \frac{c_1t}{2}\right) - \frac{c_1}{2}, \quad  u(t, x) = - \frac{c(0)}{z^2(t, x - \frac{c_1}{2} t)}
$$
satisfy the Camassa-Holm equation
$$
\begin{aligned}
u_t & = 2 u q_x + q u_x, \\
u & =  q_{xx} + 2 q.    
\end{aligned}
$$
The next problem can be solved by computer algebra.

\begin{Problem}
What solutions of Camassa-Holm are given by regular solutions (in the complex and in the real region) of the geodesic flows? Which  solutions  are given by flows on the singular surfaces (in the complex and in the real region)?    
\end{Problem}

\begin{Problem} \label{prob:20}
Is it true that,  in case of the real region,  the solution hitting the line $\lambda_2 = 0$, yields a cuspon (that is a special type of peakon with a cusp at the vertex) solution of Camassa-Holm equation?    
\end{Problem}



\subsection{Case $n = 2, d = 0$, Kaup-Boussinesq system}

This case corresponds to the case $m(\lambda) = \mathrm{const}$ and symmetry $S$ yields polynomial potential of degree $6$. W.l.o.g. we assume that $m(\lambda) = 1$. We assume that $c(\lambda)$ has six distinct real eigenvalues.

We still have separation coordinates and the formulas are the same as in Section \ref{kdvsec}. The only difference here is that $c(\lambda)$ is the polynomial of degree six.  

\begin{Problem}
Describe the singularities and the bifurcations  of the corresponding system. What happens if the solution hits the curve $w^2 - 4z = 0$?    
\end{Problem}

\begin{Problem}
Same questions for the complex case: what are the singularities, what are the bifurcations?    
\end{Problem}

The following reconstruction formulas \weg{We normalize  $c_2=0$.  The formulas
}$$\begin{aligned}
u(t, x) & = 2 w(t, x - \frac{c_1}{2}t)-c_1,  \\ v(t, x) & = - 3 w^2 (t, x -\frac{c_1}{2}t) + 2 z(t, x - \frac{c_1}{2}t)+ 2 c_1 w(t, x-\frac{c_1}{2}t )-c_2\end{aligned}
$$
give solutions of the Kaup--Boussinesq system
$$
\begin{aligned}
u_t & = v_x  + \frac{3}{2} u u_x, \\
v_t & =-\tfrac{1}{2}  u_{xxx}+ v u_x + \frac{1}{2} u v_x.
\end{aligned}
$$
The question remains as before.
 
\begin{Problem}Characterize the spatial profiles $x\mapsto (u(t,x),v(t,x))$ of the
Kaup--Boussinesq solutions obtained by the reconstruction formulas
above. Which joint orbits on regular fibers give periodic profiles,
and which give quasiperiodic profiles? For joint orbits on singular
fibers, describe the asymptotic behavior of the profiles as
$x\to\pm\infty$, whenever they are defined on the whole real line.
Consider both the region $w^2-4z>0$ and the region $w^2-4z<0$. 
\end{Problem} 

\begin{Problem}
To which  singular level surfaces do solitons correspond?
\end{Problem}


\subsection{Case $n = 2, d = 2$, two-component Camassa-Holm system}

This corresponds to the case $m(\lambda) = \lambda^2 + 1$. The main difference from the previous case is the fact that the potential is not polynomial, but rational, that is
$$
S = c(M) (M^2 + \Id)^{-1} = g_1 M + g_2 \Id.
$$
As before we assume that $c(\lambda)$ has six distinct real eigenvalues. Notice that $S$ has a pole at $z = 1, w = 0$. This differs from the previous case with singularity, where the entire line $z = 0$ was singular.

\textbf{The real case. } In separation coordinates we have:
$$
S = \left( \begin{array}{cc}
     \frac{c(\lambda_1)}{\lambda^2_1 + 1} & 0  \\
     0 & \frac{c(\lambda_2)}{\lambda^2_2 + 1} 
\end{array}\right) = g_1 \left( \begin{array}{cc}
     \lambda_1 & 0  \\
     0 & \lambda_2 
\end{array}\right) + g_2 \left( \begin{array}{cc}
     1 & 0  \\
     0 & 1 
\end{array}\right).
$$
Solving this system with respect to $g_1, g_2$ we get
$$
g_1 = \frac{(\lambda_2^2 + 1) c(\lambda_1) - (\lambda_1^2 + 1) c(\lambda_2)}{(\lambda_2^2 + 1) (\lambda_1^2 + 1) (\lambda_1 - \lambda_2)}, \quad g_2 = \frac{- \lambda_2(\lambda_2^2 + 1) c(\lambda_1) + \lambda_1(\lambda_1^2 + 1) c(\lambda_2)}{(\lambda_2^2 + 1) (\lambda_1^2 + 1)(\lambda_1 - \lambda_2)}. 
$$
  Substituting the separation relations into Hamilton's equations,
as in Section~\ref{kdvsec}  immediately before Problem~\ref{problem7}, gives for
$m(\lambda)=\lambda^2+1$  
\[
\dot\lambda_i^2=
\frac{
2\left((\lambda_i^2+1)(\lambda_i h-f)-c(\lambda_i)\right)}
{(\lambda_i^2+1)(\lambda_1-\lambda_2)^2},
\qquad i=1,2.
\]

\begin{Problem}
Describe the singularities and the bifurcations of this system. What happens if the solution hits the curve $w^2 - 4z = 0$?    
\end{Problem}

\begin{Problem}
Are there any special solutions  which can be written explicitly?    
\end{Problem}

\begin{Problem}
Same questions for the complex case: what are the singularities, what are the bifurcations?    
\end{Problem}

In the complex case the situation might be different, as we have a meromorphic operator field with two poles. One might expect some interesting effects. The relation between $u, v$ and $w, z$ is given by the complex relation
$$
1+v(t,x)-\ii u(t,x)
=-\frac{c(-\ii)}{(1-z(t,x- \tfrac{c_1t}{2})+\ii w(t,x- \tfrac{c_1t}{2 }))^2},
\quad
q(t,x)=w(t,x-\tfrac{c_1t}{2})-\frac{c_1}{2}.   
$$
The corresponding integrable system with dispersion is a two-component Camassa-Holm equation
(recall the  formulas in the  introduction for the change of $x, t$).
$$
\begin{aligned}
u_t & = 2 u q_x + q u_x + v_x, \\
v_t & = 2v q_x + q v_x + 2 q_x - u_x, \\
u & = 2 q + q_{xx}. 
\end{aligned}
$$
We have the following problems.

\begin{Problem}
What solutions of the two-component Camassa-Holm are given by regular solutions (in the complex and in the real region) of the geodesic flows?     
\end{Problem}

\begin{Problem}
Which solutions  are given by flows on the singular surfaces (in the complex and in the real region)?   What role do  the poles play? 
\end{Problem}

\begin{Problem}
Are there any solutions  that can be found explicitly?    
\end{Problem}


\section{Possible analogues of famous results for BKM systems}
\subsection{Exact solutions of BKM systems}

A popular topic in the theory of integrable systems is the construction of exact solutions. ``Exact'' means that the solution is given by certain formulas. In the best case, these formulas include only elementary functions. This ``best case'' is quite rare in the literature. Exact solutions may be implicit, include special functions, and results of integration.

We suggest searching for exact solutions of BKM systems using the finite-dimensional reduction methods obtained in \cite{BKMreduction} and discussed above. Recall that the method allows one to find certain solutions of BKM systems by first solving certain explicit finite-dimensional integrable systems and then applying explicit algebraic formulas which transform solutions of the finite-dimensional system into solutions of BKM systems.

Let us first explain a naive approach. We will consider the BKM system (of any type) corresponding to parameters $n,  m(t)$ and the corresponding $N$-dimensional Benenti system. We know that in diagonal coordinates $q_1,\dots, q_N$, the Benenti system  is a Stäckel system with commuting integrals $H_0,\dots, H_{N-1}$ satisfying
\begin{equation}\label{eq:stackel1}
\begin{pmatrix} q_1^{N-1} & q_1^{N-2}& \cdots & 1 \\ 
                               q_2^{N-1} & q_2^{N-2}& \cdots & 1 \\
                                   \vdots & & & \vdots \\ 
                                   q_N^{N-1} & q_N^{N-2}& \cdots & 1\end{pmatrix} \begin{pmatrix} H_0 \\ 
                               H_1 \\
                                   \vdots   \\ 
                                  H_{N-1}\end{pmatrix} = 
\begin{pmatrix}   \frac{1}{2}p_1^2 + f(q_1)\\ 
                                 \frac{1}{2}p_2^2 + f(q_2)\\
                                   \vdots   \\ 
                                   \frac{1}{2}p_N^2 + f(q_N)\end{pmatrix},
\end{equation}
where the function $f(\mu)$ is given by 
\begin{equation} \label{eq:Umu1}
f(\mu) =  \frac{\mu^{2N+n} + c_{1}\mu^{2N+n-1}+ c_{2}\mu^{2N+n-2}+ c_{3}\mu^{2N+n-3}+ \cdots + c_{2N+n}}{m(\mu)} =  \frac{c(\mu)}{m(\mu)}
\end{equation}
with arbitrary constant coefficients $c_1, c_2,\dots, c_{2N+n}$ in BKM I and BKM II systems. In the BKM III and BKM IV systems, the essential coefficients are $c_2,\dots, c_{2N+n}$ as  we assume    $c_1=0$ for simplicity.  The general case differs by the Galilean shift described in Theorem \ref{t2}.  No such normalisation is needed for   BKM I and BKM II. 

Such systems can be integrated by the method of separation of variables: if we denote by $t_0$ the time corresponding to the Hamiltonian system generated by $H_0$, $t_1$ the time corresponding to the Hamiltonian system generated by $H_1$, and so on, we obtain  
\begin{equation} \label{eq:separation}
    \begin{array}{ccc}
\sqrt{2}\,t_0 & =  & \sum_{i=1}^N \pm \int^{q_i} \frac{s^{N-1}} {\sqrt{-c(s)/m(s)}}  ds \\
\sqrt{2}\, t_1  & = &   \sum_{i=1}^N \pm \int^{q_i} \frac{s^{N-2}} {\sqrt{-c(s)/m(s)}} ds \\
&\vdots& \\ 
\sqrt{2}\, t_{N-1}  & = &  \sum_{i=1}^N \pm \int^{q_i} \frac{1} {\sqrt{-c(s)/m(s)}} ds .
\end{array} 
\end{equation}
This is a system of explicit functional equations on $q_1, \dots, q_N$, depending on $t_0,\dots, t_{N-1}$ as parameters. Solving this system with respect to $q_i$, we obtain $q_i(t_0,\dots,t_{N-1})$.
    For BKM III and IV, assuming $c_1=0$ and
    $N\geq 2$, we set
    \[
    t_0=x,\qquad t_1=t,
    \]
    and keep the remaining commuting times fixed.

    For BKM I and II, let
    $a=\sqrt{c(\lambda)}$.
    Since the corresponding reduced flow is generated
    by
    \[
    F_\lambda
    =-\sum_{j=0}^{N-1}
      \lambda^{N-1-j}H_j,
    \]
    we set
    \[
    t_j=x\delta_{j0}
        -\tfrac{1}{a}\lambda^{N-1-j}t,
    \qquad j=0,\ldots,N-1.
    \]
    In particular, for $\lambda=0$ and $N\geq2$,
    only $t_{N-1}$ varies in addition to $t_0$.

    The resulting functions $q_i(x,t)$ are then
    transformed into solutions of the BKM system
    by the reconstruction map $\mathcal R$ described in
    Theorem~\ref{t2}.

\begin{Remark}
    One of the most famous exact solutions of the KdV equation, the two-soliton solution, can be obtained by this method (Multi-soliton solutions of KdV were first obtained by the inverse
scattering method \cite{GGKM1967} and by Hirota's direct method
\cite{Hirota1971}.) 
In this case, the polynomial $m(\mu)$ is constant, and the polynomial $c(\mu)$ has degree $5$ and has two double roots  and one simple root. Then, one can solve the integrals appearing in \eqref{eq:separation} explicitly, and then solve the system \eqref{eq:separation} with respect to $q_1, q_2$ to obtain exact solutions of this finite-dimensional system, and finally of the KdV system.
\end{Remark}

Note that the complexity of the system increases with increasing $N$. The integrals in \eqref{eq:separation} become more complicated, and even if one can evaluate the integrals (for example, using special functions), one obtains a complicated system of algebraic equations. For certain carefully chosen polynomials $c(\mu), m(\mu)$, one can still find exact solutions; see, e.g., \cite[Example 1.5]{BKMreduction}, but for arbitrary $c$ and $m$, this naive approach does not work. 

We suggest approaching this problem via Riemann surfaces. The ``correct'' Riemann surface in this setup is $$\frac{1}{2}p^2 = -\frac{c(q)}{m(q)}$$ (we assume that both $p$ and $q$ are complex-valued). We may think that the variables $q_1,\dots, p_N$ evolve on this Riemann surface. The integrands in the formula \eqref{eq:separation} can be viewed as differentials on this Riemann surface. It is known that, on a compact Riemann surface, the set of differentials with   fixed numbers  and types of zeros and poles is finite-dimensional. This allows one to solve the system \eqref{eq:separation} in special functions.   This approach is very classical: the system \eqref{eq:separation} corresponds to  an Abel map
written in the coordinates $q_1,\dots,q_N$, and solving it with respect to
the $q_i$ is the Jacobi inversion problem.  For St\"ackel systems this goes
back to Jacobi and St\"ackel \cite{Staeckel1893}. As far as we know, it is
essentially the only approach producing quasi-periodic (finite-gap)
solutions of integrable infinite-dimensional systems. 

\begin{Problem}
   Apply this method, or its naive version, to other BKM systems. For BKM IV and BKM III systems with constant polynomial $m$, the methods should work; the case when $m$ is not constant may be more complicated. 
\end{Problem}
 
In the context of integrable infinite-dimensional systems, this procedure
was first applied to the KdV equation in the works of Novikov \cite{novikov74},
Dubrovin \cite{Dubrovin1975}, Its and Matveev \cite{ItsMatveev1975},
Lax \cite{Lax1975}, and McKean and van Moerbeke \cite{McKeanVanMoerbeke1975};
see the survey \cite{DubrovinMatveevNovikov1976}. A general algebro-geometric
scheme, applicable to the KP equation and to general Zakharov--Shabat
equations, was developed by Krichever \cite{Krichever1976, Krichever1977}.

\subsection{Uhlenbeck-Veselov-Moser for Benenti systems corresponding to BKM systems} \label{sec:moser}

We consider the Benenti system corresponding to the $N$-gap reduction of an $n$-component BKM system, and concentrate on the regions where the eigenvalues of the corresponding $N$-dimensional operator $M$ are real and distinct. Then, one can take them as coordinates. We denote them $q_1 < \dots < q_N$. As discussed above, in this coordinate system the operator $M$ is diagonal, $M = \textrm{diag}(q_1,\dots, q_N),$
and the corresponding Benenti system has the Stäckel form 

\begin{equation}\label{eq:st3}
\begin{pmatrix} (q_1)^{N-1} & (q_1)^{N-2}& \cdots & 1 \\ 
                               (q_2)^{N-1} & (q_2)^{N-2}& \cdots & 1 \\
                                   \vdots & & & \vdots \\ 
                                   (q_N)^{N-1} & (q_N)^{N-2}& \cdots & 1\end{pmatrix} \begin{pmatrix}  I_0 \\ 
                                I_1 \\
                                   \vdots   \\ 
                                    I_{N-1}\end{pmatrix} = 
\begin{pmatrix} \tfrac{1}{2}f(q_1) (p_1)^2 + U(q_1)\\ 
                      \tfrac{1}{2} f(q_2) (p_2)^2 + U(q_2)\\
                                   \vdots   \\ 
                      \tfrac{1}{2} f(q_N) (p_N)^2 + U(q_N)\end{pmatrix}
\end{equation}
with explicit functions $f$ and $U$. Actually, the function $f$ for the finite-dimensional reductions of BKM systems is $f \equiv 1$, and the function $U(t) = c(t)/m(t)$, where $c$ is a monic polynomial of degree $n + 2N$ with arbitrary coefficients (for BKM IV, the second leading term is zero), and $m(t)$ is the polynomial used for the construction of the BKM system.

Let us now recall the following relatively easy result:

\begin{Lemma}[Lemma 2 of \cite{BKMmoser}] \label{lem:2}
Consider two families of Benenti systems sharing the same $L$: one constructed by $U(t)$ and $f(t)$, another by $\widetilde U(t)$ and $\widetilde f(t)$. Suppose for certain values of the integrals $I_0 = h_0, \dots,  I_{N-1} = h_{N-1}$ of the first system and for certain values of the integrals  
$\widetilde I_0 = \widetilde h_0, \dots, \widetilde I_{N-1} = \widetilde h_{N-1}$ of the second system the products 
satisfy
\begin{equation} \label{eq:condition} \begin{array}{l}
(-U(t) +  h_0 t^{N-1} + \cdots + h_{N-2} t + h_{N-1}) f(t)\\ \qquad   = (-\widetilde U(t) + \widetilde h_0 t^{N-1} + \cdots + \widetilde h_{N-2} t + \widetilde  h_{N-1}) \widetilde f(t).\end{array}
\end{equation} 
Then, every solution of the first Benenti system, viewed as a curve on $\mathbb{R}^N$, corresponding to the values of the integrals $h_0, \dots, h_{N-1}$, is a solution of the second Benenti system corresponding to the values of the integrals $\widetilde h_0, \dots, \widetilde h_{N-1}$.  
\end{Lemma}

Using Lemma \ref{lem:2}, we can change the function $f$ in \eqref{eq:st3} to make it a polynomial of degree $N+1$ with $N+1$ different real roots and with leading coefficient $-4$. The kinetic energy of the integral $I_0$  
is then a (contravariant) metric of constant curvature $1$, and if we correctly choose the ranges of the coordinates $q_1, \dots, q_N$, it is positive definite.  
The corresponding Stäckel system then gives an integrable system on the $N$-sphere, such that the integral 
$I_0$ is the sum of the kinetic energy of the standard metric on the sphere and a certain potential. 

\begin{Problem}\label{prob:M5}
Find the expression for the potential in the ambient coordinates in $\mathbb{R}^{N+1}.$ Find also the expression for the integrals in ambient coordinates.
\end{Problem}

In the KdV case (BKM IV, $n=1$), the corresponding integrable system on the sphere is the famous Neumann system~\cite{Neumann1859}. This was observed by Ju.~Moser~\cite{moser80,moser81, Moser1980Chern}  and A. Veselov \cite{veselov}. The expression for the integrals in ambient coordinates can again be found in \cite{moser81}, and according to Moser, is due to Uhlenbeck ~\cite{Uhlenbeck1982}. 

Note that for the Riemannian metric of constant curvature 1, the coordinates $q_1, \dots, q_n$ are essentially the famous sphero-ellipsoidal coordinates. The formula relating the ambient coordinates and the sphero-ellipsoidal coordinates is of course well known (and can be found in e.g.,  \cite{separation,DullinRichterVeselovWaalkens2001},  or via search engines) which  makes Problem \ref{prob:M5}, for small dimensions, an exercise in computer algebra.   

Let us recall the formula. Let  $a_0<\cdots<a_N$ be the roots of the polynomial $f$ chosen
above, and let $(X_0,\ldots,X_N)$ be the ambient Cartesian
coordinates on the unit sphere $S^N\subset\mathbb{R}^{N+1}$.
The sphero-ellipsoidal coordinates $q_1,\ldots,q_N$ are related
to the ambient coordinates by the identity
\[
\sum_{\alpha=0}^{N}\frac{X_\alpha^2}{\lambda-a_\alpha}
=
\frac{\prod_{i=1}^{N}(\lambda-q_i)}
     {\prod_{\alpha=0}^{N}(\lambda-a_\alpha)}.
\]
Taking the residue at $\lambda=a_\alpha$ gives
\[
X_\alpha^2=
\frac{\prod_{i=1}^{N}(a_\alpha-q_i)}
     {\prod_{\beta\ne\alpha}(a_\alpha-a_\beta)},
\qquad \alpha=0,\ldots,N.
\]
See \cite[eqs.~(15)--(16)]{DullinRichterVeselovWaalkens2001}.
These formulas give a direct starting point for the computation
in Problem~\ref{prob:M5}.

\begin{Problem}
If $m(t) = \const$ or the zeros of $m$ lie outside of the ranges of the coordinates $q_i$, the isoenergy hypersurface is compact, so the Liouville foliation has compact fibers. Describe the topology and singular points of this Liouville foliation in small dimensions. Make visualizations of the trajectories. 
\end{Problem}

Note that the solution of this problem will give a new, interesting class of natural integrable systems. These systems generalize the famous Neumann system and are also separable in sphero-ellipsoidal coordinates. One can then study them separately, asking and answering standard questions about integrable systems (i.e., construct the  marked molecule, symplectic invariants, action-angle variables, etc., see~\cite{BolsinovFomenko2004}). 

\begin{Problem}
If the zeros of $m$ lie inside the range of a coordinate $q_i$, the potential is unbounded and the isoenergy hypersurface is not compact. Make visualizations of the trajectories. Pay special attention to the case $m(t)=t$,
which is related to the Camassa--Holm equation
and may be relevant to the appearance of
peakon-type singular solutions.
\end{Problem}

As discussed above, the finite-gap solutions of the BKM systems are related to the solutions of the Benenti systems via a transformation of the form $u = u(q)$. 

\begin{Problem}
Find the transformation $u(q)$ in ambient coordinates (i.e., find the functions 
$u_1, \dots, u_n$ as functions of standard coordinates in $\mathbb{R}^{N+1}$). 
\end{Problem}

For the KdV case, the answer is known and is very nice: up to a factor and an additive constant, 
the function $u$ coincides with the potential energy of the Neumann system.

Related finite-dimensional constructions were
studied by Qiao using the nonlinearisation of
spectral problems and Neumann-type constraints.
In particular, generalized C.~Neumann systems
associated with several spectral problems,
including energy-dependent ones, were constructed
in \cite{Qiao1997GeneralizedNeumann};
see also \cite{Qiao2003CamassaHolm} for a
Neumann-like constrained system associated with
the Camassa--Holm hierarchy and its
algebro-geometric solutions.

\subsection{The peakon metric: Lax representation and polynomial Killing tensors}
\label{s:peakon_polynomial}

The Camassa-Holm (CH) equations in the classical article \cite{CH93} differ from the one we construct in section \ref{s:HD} in the form of constraint: it has $u = q - q_{xx}$, not $u = q_{xx} + 2 q$. This corresponds to the different choice of the sign of $m(\lambda)$ in our construction. In this section, we discuss the classical version.

General peakon solutions for the CH equation are related to different finite-dimensional integrable systems, found in the literature; see e.g. \cite{CI10, HI10, CI10, CK}. In this case, the $N$-peakon ansatz
\[
u(x,t)=\sum_{a=1}^N p_a(t)\,e^{-|x-q_a(t)|}
\]
defines a finite-dimensional invariant manifold on which the PDE reduces to an integrable Hamiltonian system. In a
chamber $q_1<\cdots<q_N$ of the configuration space its Hamiltonian is
\begin{equation}\label{eq:CHpeakonHamiltonian}
 H(q,p)=\frac12\sum_{i,j=1}^N p_i p_j e^{-|q_i-q_j|}
       =\frac12\,g^{ij}(q)p_i p_j,
 \qquad g^{ij}(q)=e^{-|q_i-q_j|}.
\end{equation}
In particular, the corresponding geodesic flow is defined by a
Riemannian metric as long as the peak positions are pairwise distinct.
For $p_i>0$, one convenient form of its finite-dimensional Lax pair is
\begin{equation}\label{eq:CHpeakonLax}
 \mathcal L_{ij}=\sqrt{p_i p_j}\,e^{-|q_i-q_j|/2},\qquad
 \mathcal A_{ij}=\begin{cases}
 \frac12\operatorname{sgn}(i-j)\mathcal L_{ij},&i\ne j,\\
 0,&i=j,
 \end{cases}
 \qquad \dot{\mathcal L}=[\mathcal A,\mathcal L].
\end{equation}
The version of the Lax matrix used in \cite{CH93} is
$L_{ij}=p_j e^{-|q_i-q_j|/2}$; it is conjugate to $\mathcal L$
when all $p_i>0$. Its characteristic polynomial and trace invariants
extend as polynomials in the momenta without this positivity assumption.

The coefficients $E_k$ of
\begin{equation}\label{eq:CHpeakonSpectral}
 \det(I+t\mathcal L)=\sum_{k=0}^N E_k t^k,\qquad E_0=1,
\end{equation}
are homogeneous polynomials in the momenta, of degrees $k=1,\ldots,N$.
They Poisson-commute and are functionally independent on an open dense
subset. In particular, $E_1=\sum_i p_i$ and
$2H=E_1^2-2E_2$. Thus, the geodesic flow of the peakon metric is
Liouville integrable \emph{within the class of polynomial first integrals},
with one integral of each degree $1,\ldots,N$; this is classical
\cite{CH93}. 

There is also superintegrability in the class of general
smooth first integrals: using spectral/action--angle variables,
Holm and Ivanov \cite{HI10} exhibit, on the regular part of the
peakon phase space, $N-1$ additional first integrals, giving $2N-1$
functionally independent integrals. The additional integrals need not
be polynomial in momenta.

\begin{Problem}
Determine the graded algebra of polynomial first
integrals (equivalently, of Killing tensors) of the peakon metric
\eqref{eq:CHpeakonHamiltonian}. Is every polynomial first integral a
polynomial in $E_1,\ldots,E_N$? In particular, are there any Killing
tensors not generated by these spectral invariants? Investigate this
question also for general $N$ and for higher polynomial degrees.    
\end{Problem}

Preliminary brute-force computer-algebra calculations of the prolonged
Killing equations support a negative answer to the question about
additional polynomial integrals. They give no extra integrals for
$N=3$ in degrees $2$ and $3$, and for $N=4$ in degree $2$.

For $N=4$ in degree $3$ the same conclusion has been checked in the
translation-invariant class, but a complete exclusion of
non-invariant integrals is not asserted here. These are computational
observations, rather than a general classification theorem.

Notice that the linear integral $E_1=\sum_i p_i$ generates simultaneous
translations of the $q_i$. Its Poisson bracket acts on each
finite-dimensional space of Killing tensors of a fixed degree.
Consequently, one can study its (possibly complex) eigenvectors and
generalised eigenvectors, as in the method discussed in
Section~\ref{sec:tannery}. This symmetry can substantially simplify
the relevant Killing equations.

There is a useful coordinate system in which the \emph{contravariant}
metric has polynomial entries and the \emph{covariant} metric is
tridiagonal. On the chamber $q_1<\cdots<q_N$ set
\begin{equation}\label{eq:CHpeakonExpChange}
 x_i=e^{q_i},\qquad \pi_i=e^{-q_i}p_i,
 \qquad 0<x_1<\cdots<x_N.
\end{equation}
This is a canonical change of coordinates and momenta. In these
coordinates
\begin{equation}\label{eq:CHpeakonPolyMetric}
 \widetilde g^{ij}(x)=x_{\min(i,j)}^2,\qquad
 H=\frac12\sum_{i,j=1}^N x_{\min(i,j)}^2\pi_i\pi_j,
\end{equation}
while the covariant metric is
\begin{equation}\label{eq:CHpeakonTriMetric}
 g=\frac{dx_1^2}{x_1^2}
   +\sum_{i=2}^N\frac{(dx_i-dx_{i-1})^2}{x_i^2-x_{i-1}^2}.
\end{equation}
In particular, its covariant coefficient matrix is tridiagonal.

The Lax pair also acquires a polynomial form. Set
\begin{equation}\label{eq:CHpeakonPolyLax}
 B_{ij}=x_{\min(i,j)}\pi_j,
 \qquad
 A_{ij}=\begin{cases}
  \frac12\operatorname{sgn}(i-j) B_{ij},&i\ne j,\\[2pt]
  \frac12 B_{ii}+\displaystyle\sum_{k=i+1}^N B_{ik},&i=j.
 \end{cases}
\end{equation}
Then, for the Hamiltonian \eqref{eq:CHpeakonPolyMetric},
\begin{equation}\label{eq:CHpeakonPolyLaxEquation}
 \dot B=[A,B].
\end{equation}
Indeed, when $p_i>0$, the matrices $B$ and $\mathcal L$ are conjugate
by the diagonal matrix
$D=\operatorname{diag}(\sqrt{x_i/p_i})$; the Lax equation
\eqref{eq:CHpeakonPolyLaxEquation} and its polynomial coefficients
extend to arbitrary momenta. In dimensions $N=3,4$ this formula can
also be verified directly by substituting Hamilton's equations and
checking each matrix entry.

The spectral invariants of $B$ now admit particularly simple expressions.
Write
\begin{equation}\label{eq:CHpeakonPolyInvariants}
 \det(I+tB)=\sum_{k=0}^N E_k t^k,\qquad
 E_k=\sum_{1\le i_1<\cdots<i_k\le N}
 \left[x_{i_1}\prod_{r=2}^k
      (x_{i_r}-x_{i_{r-1}})\right]
 \pi_{i_1}\cdots\pi_{i_k}.
\end{equation}
In particular,
\[
 E_1=\sum_{i=1}^N x_i\pi_i,\qquad
 E_2=\sum_{i<j}x_i(x_j-x_i)\pi_i\pi_j,
 \qquad
 E_N=x_1\prod_{i=2}^N(x_i-x_{i-1})\prod_{i=1}^N\pi_i.
\]
To see \eqref{eq:CHpeakonPolyInvariants}, write
$B=M\operatorname{diag}(\pi_1,\ldots,\pi_N)$ with
$M_{ij}=\min(x_i,x_j)$. For every ordered subset
$I=\{i_1<\cdots<i_k\}$, the corresponding principal minor satisfies
\[
 \det M_I=x_{i_1}\prod_{r=2}^k
            (x_{i_r}-x_{i_{r-1}}),
\]
and the general formula follows by expanding $\det(I+tB)$ in
principal minors. These formulas make the whole known family of
commuting integrals polynomial both in $x_i$ and in $\pi_i$.

\begin{Remark}\label{peakon_remark}
The finite-dimensional reduction discussed in Section~\ref{s:HD}
seems to produce, in the peakon regime considered there, only a single
peakon of the cuspon type. In contrast, the construction of \cite{CH93} gives genuine multipeakon solutions with an arbitrary number of peaks (which are not cusps). Thus the two finite-dimensional constructions are different; understanding their precise relationship, if any, would be of interest.
\end{Remark}

\subsection{Novikov--Dubrovin--Moser for BKM systems.}

S. Novikov has  shown in  \cite{novikov74} that stationary solutions of (higher) KdV  give quasiperiodic   potentials such that the corresponding  Schr\"odinger-Hill operators  have  finitely many connected components  of spectrum.  Slightly later, B. Dubrovin \cite{Dubrovin1975} has shown that all quasiperiodic potentials for which  the corresponding  Schr\"odinger--Hill operators  have  finitely many connected components  of spectrum correspond to stationary solutions of KdV. This correspondence, and also its relation to the Neumann problem discussed in \S \ref{sec:moser} 
were partially rediscovered and popularised by J. Moser \cite{moser80,moser81}.  We suggest   generalizing  these  results for all BKM systems.

Let us first recall the result of \cite{Dubrovin1975} using our terminology. We consider the solution $u(x,t)$  of the  KdV equation coming from our finite-dimensional reduction corresponding to the parameter $N$ and to  the polynomial $c$ with $2N+1$ real roots which we denote by $r_1\le r_2\le \cdots \le r_{2N+1}$.
 
Now, fix $t$ and consider the  linear second order operator  $-\tfrac{\partial^2 }{\partial x^2} + \tfrac{1}{2} u$ and the
corresponding Schrödinger  equation 

\begin{equation}\label{eq:schrodinger}
   (-\tfrac{\partial^2 }{\partial x^2} + \tfrac{1}{2} u)\Psi= \tfrac{\lambda}{2} \Psi .
\end{equation}
Note that \eqref{eq:schrodinger} is an ordinary differential equation in $x$,  and $t$ should be viewed as a parameter. 

We say that $\lambda$ lies in the spectrum of this equation if there exists a bounded nontrivial solution $\Psi$.  
\begin{Lemma}[Theorem 4 of \cite{BKMmoser}]\label{lem:p3} 
 In the notation above, $\tfrac{1}{2} \lambda$ lies in  the spectrum of $ -\tfrac{\partial^2 }{\partial x^2} + \tfrac{1}{2} u$ if and only if   
\begin{equation} \label{eq:bands} \lambda \in [r_1, r_2] \cup [r_3, r_4]\cup \cdots \cup [r_{2N-1}, r_{2N}] \cup [r_{2N+1}, +\infty). \end{equation} 
\end{Lemma}

We see that the spectrum consists  of finitely many closed intervals and one half-line $\left[\tfrac{1}{2}r_{2N+1}, +\infty\right)$. This is very untypical, as for a generic periodic or quasiperiodic potential  the spectrum consists of infinitely many intervals and one half-line.

In jargon, the intervals are called {\it bands} and the  half-line $[r_{2N+1}, +\infty)$ is called {\it the  infinite band}. The connected components of the
complement of the spectrum are called {\it gaps}.   This terminology explains the terminology ``finite-gap solution of KdV''. 

Note that the number of (finite)  bands  is at most $N$.  Since $r_i$ may be equal to $r_{i+1}$, the number of bands can be smaller than $N$ and some of the bands could   actually  be points.

\begin{Problem} \label{prob:M6}
    Generalize Lemma \ref{lem:p3} for other BKM systems. 
\end{Problem}

Note that the case $n=1$ can be done straightforwardly, and possibly has been done before, as all $n=1$  BKM systems  were known and studied from different aspects before. 
In the case $n>1$,    we have $n$ unknown functions $u_1(x,t),\dots, u_n(x,t). $
This implies that  the  spectral problem   is expected to look very different, as the equation is multicomponent, so instead of one function $u(x,t)$ we should somehow employ, in the spectral problem,  $n$ functions $u_1(x,t),\dots, u_n(x,t)$. 

  We encode the functions $u^1,\ldots,u^n$ in the characteristic
polynomial of the Nijenhuis operator $L$ defining the BKM system,
\begin{equation}\label{eq:sigma} 
  \sigma(\lambda,u)=\det(\lambda\,\mathrm{Id}-L(u)).
\end{equation}
Indeed, in the companion coordinates of~\eqref{pcf}, we have
\[
  \sigma(\lambda,u)
  =\lambda^n-u^1\lambda^{n-1}-\cdots-u^n.
\]
Thus the coefficients of $\sigma$ recover all $n$ components of $u$.

The analog of the 
equation \eqref{eq:schrodinger} in this case is the equation 
\begin{equation}\label{eq:SL2} \psi_{xx}+ \frac{\sigma(\lambda, u(x)) }{2m(\lambda)}\psi    =0.  \end{equation}

We see that for $n\ge 2$ the dependence on the $\lambda$ is not linear, unlike in  the 
Schr\"odinger equation.  Indeed, the function   $$
  \frac{\sigma(\lambda,u(x))}{2m(\lambda)}
 $$  is a  rational function in $\lambda$ whose coefficients  depend on  $x$. 

The special case $m(\lambda) = \textrm{const}=m_0$ of this equation   appeared in the literature, see  e.g. \cite{af,alonso,JaulentMiodek1976}, under the name ``energy-dependent potentials'', in the relation to multicomponent integrable systems. The scattering problem for such systems was 
studied  quite recently, see e.g.  \cite{nab,HrMa}.

We say that $\lambda$ lies in the spectrum of \eqref{eq:SL2}  if there exists a bounded nonzero solution $\Psi$. We expect that for  certain $u$ coming from finite-dimensional reductions of the BKM systems, the spectrum is the 
finite-band one, in the sense that  it  contains finitely many connected components, and Problem \ref{prob:M6} is precisely about it.   Of course, it will be  interesting to understand whether in the class of quasi-periodic coefficients of  the equation \eqref{eq:SL2}  the finite-bandness of the spectrum is equivalent to the property that the coefficients came from a BKM system; the next section will be devoted to this question. Moreover, as mentioned above, it will generally  be interesting to develop the inverse scattering approach to BKM systems; the analogue of the corresponding Sturm-Liouville equation is expected to be \eqref{eq:SL2}. 

\subsection{Rational linear problem and energy-dependent potentials}

We start with a general formulation of the energy-dependent linear problem. Consider a differential operator $H_\lambda : C^\infty(\mathbb R) \to C^\infty (\mathbb R)$ on the real line in the form
$$
H_\lambda = D^2 + \sigma(\lambda) \operatorname{Id}.
$$
Here $\sigma(x,\lambda)$ is a rational function
of the spectral parameter $\lambda$, whose
coefficients are smooth functions of $x$. 
  We refer to $\sigma(x, \lambda)$ as an energy-dependent potential. For a given $\lambda$, consider the ODE in the form
\begin{equation}\label{lin1}
H_\lambda \psi = 0.
\end{equation}
We say that $\lambda$ belongs to the spectrum of the energy-dependent operator $H_\lambda$ if there exists a solution $\psi$ of \eqref{lin1} which is bounded by some constant on the entire line $\mathbb R$. Here are some examples.

\begin{Ex}
\rm{
Consider $m = 1$, $\sigma(x, \lambda) = \lambda + u(x)$. In this case, the equation \eqref{lin1} takes the form
$$
H_\lambda \psi = \psi_{xx} + u(x) \psi + \lambda \psi = 0.
$$
We arrive at the standard spectral problem of the Sturm–Liouville operator, well known in the literature. 
}
\end{Ex}

\begin{Ex}
\rm{
Consider $$m = 1 \, , \   \sigma(x, \lambda) = \lambda^n + u_1(x) \lambda^{n - 1} + \dots + u_n(x).$$ In this case, the equation \eqref{lin1} takes the form
$$
H_\lambda \psi = \psi_{xx} + \sum_{i = 1}^n \lambda^{n - i} u_i(x) \psi + \lambda^n \psi = 0.
$$
In this case, we arrive at the problem of finding the spectrum of a potential that depends on the energy in a polynomial sense. These potentials were introduced by Jaulent and Miodek~\cite{JaulentMiodek1976} and Mart\'{\i}nez Alonso~\cite{alonso} and later studied by Antonowicz and Fordy \cite{af}. The scattering problem for such potentials came to attention only quite recently (see \cite{nab}, \cite{HrMa}). 
}    
\end{Ex}

We start with the easiest problem.

\begin{Problem}
Given the formulation, construct finite-gap type solutions in the rational case.    
\end{Problem}

The examples from the previous section might be useful. 

\begin{Problem}
Is it true that all finite-gap type solutions come from a Benenti-type system for an appropriate choice of symmetry $S$ (and functions $g_1$, $g_2$, respectively)?
\end{Problem}

The next problem seems to be hard.

\begin{Problem}
What is the solution of the rational problem in the general case? If the roots of the denominator are complex, is it true that the spectrum in general consists of an infinite number of gaps?    
\end{Problem}

\section{Do finite-dimensional reductions determine the BKM system?}

BKM systems  form
a family of multicomponent PDE systems introduced in \cite{nijapp4}; we described it in Section \ref{s:descriptionBKM}. Their special cases include famous and actively studied systems such as KdV, coupled KdV, Harry Dym, coupled Harry Dym, Camassa--Holm, multicomponent Camassa--Holm \cite{HI10}, Dullin--Gottwald--Holm, and Kaup--Boussinesq equations. 

One can understand a BKM system  as a dynamical system on the space of curves on an $n$-dimensional manifold: for a local curve $\hat u(x)$, the solution of a BKM can be viewed as a family $u_t(x) = u(x,t)$ of curves depending on $t$.  

The paper \cite{BKMreduction} suggests a construction of special solutions of this system. The general structure is as follows: one chooses a number $N$ (``number of gaps'') and constructs, by a certain explicit formula, a finite-dimensional Hamiltonian system on the cotangent bundle to an $N$-dimensional manifold, whose Hamiltonian is the sum of a kinetic energy corresponding to a flat (pseudo-Riemannian) metric and a certain potential energy. The system is integrable; it admits $N$ functionally independent quadratic integrals in involution. We take one of the integrals and denote by $t$ the time corresponding to the flow along the integral. Next, consider a solution of the system $(w_1(x,t),\dots, w_N(x,t))$. The paper \cite{BKMreduction} constructs, from this solution, a solution $u(t,x)$ of the BKM system. Such solutions will be called {\it finite-gap} solutions. 

\begin{Problem}
\label{prob:M1}
Is it possible, for a generic initial curve $\hat u(x)$ and for any $k$, to find a finite-gap solution $u(x,t)$ such that at the point $x = 0$, the curve $u(x,0)$ approximates the curve $\hat u(x)$ up to order $k$?
\end{Problem}

For the KdV system, the answer is classically known and is positive. 
New nontrivial advances on  this topic are due to  Quaschner and Steneker
\cite{QuaschnerSteneker2026}. 
They proved that for BKM systems with $\deg m=0$,
every prescribed finite jet of initial data can be
realised by a finite-gap solution. In particular,
this applies to the  already known  KdV case and to the   Kaup--Boussinesq equations.
For one-component BKM systems with $\deg m=1$,
including the Camassa--Holm equation, they proved
the same result for an open set of real jets and
for a Zariski-open dense set of complex jets.

\begin{Problem}
Extend the jet-approximation results of
\cite{QuaschnerSteneker2026} to general BKM systems.
In particular, does every prescribed finite jet of
initial data arise from a finite-gap solution?
Can the results for $\deg m=1$ be extended to
multicomponent systems, and can the restrictions
on the admissible jets be removed?
\end{Problem}

The problem can be tackled using computer algebra. Namely, the formulas for transforming the solutions $w(x,t)$ to $u(x,t)$ are explicit, though complicated; see \S \ref{sec:Rmap}.  

\begin{Problem}
\label{prob:M2}
Can one obtain BKM systems as limits of finite-dimensional systems? That is, can one approximate a generic curve (or a certain class of curves) by finite-gap solutions and obtain a system of ODEs which converges to the BKM PDE system?
\end{Problem}

For special cases of BKM systems, the answer is positive; see \cite{BM2006, BM2008, BS2023}. 

Let us give more details. The finite-dimensional reduction from \cite{BKMreduction} can be geometrically viewed as follows. As mentioned above, one can understand BKM systems as a dynamical system on the space of (local, real-analytic, parameterised) curves on the $n$-dimensional manifold. Finite-dimensional reduction can be viewed as the restriction of the dynamical system 
to a finite-dimensional subset of the set of curves (for $N$ ``gaps'', the subset is $2N$-dimensional). We will call the subset an $N$-gap family. The dynamics with respect to the $t$-variable then gives a vector field along the curves from the finite-dimensional subset.   Problem \ref{prob:M2} is as follows:  when we approach an arbitrary curve by a sequence of curves from $N$-gap families with $N \to \infty$, do  we obtain the generator of the PDE?

    Of course, one needs to check whether other structures related to the finite-dimensional systems converge to those of the BKM systems. These other structures include quadratic-in-momenta integrals (whose infinite-dimensional analogs are called conservation laws), symplectic and Poisson structures (whose infinite-dimensional analog is called a Hamiltonian operator).

\begin{Problem}
\label{prob:M3}
The finite-dimensional systems appearing in the procedure of finite-dimensional reduction are quantum integrable. That is, the Carter quantization of the integrals yields commuting  second-order operators. Will the limit of these quantizations give a quantization of the BKM system?
\end{Problem}

Let us discuss the terms used in the problem. In the finite-dimensional systems, the procedure of Carter quantization~\cite{Carter1977} is essentially the replacement of the commuting integrals $K_i + U_i$ (where $K_i$ are the quadratic-in-momenta parts and $U_i$ are the potential parts) by second-order differential operators of the form  
$\hat K_i + U_i$, with  
$$
\hat K_i = \nabla_s K_i^{sr} \nabla_r,
$$
where $\nabla$ is the Levi-Civita connection of the metric. It is known \cite{Bbook, KM2025, MT01, Matveev2020, DuvalValent2005}   that for the finite-dimensional systems appearing in the procedure of finite-dimensional reduction, the operators $\hat K_i + U_i$ commute. 

Note that in our setup the metric is flat. That implies that in a certain coordinate system, namely in the coordinate system in which the metric is given by constant entries, the Carter quantization coincides with the Weyl quantization. Recall that Weyl quantization is essentially the replacement of momenta $p_i$ by the differential operators $i\hbar \partial_{x_i}$ (plus a certain completely symmetrized ordering procedure to make the output of the construction self-adjoint). This allows the definition of quantization for real-analytic functions. 

Another equivalent viewpoint on this procedure is that one replaces the Poisson brackets of functions by commutators of the corresponding differential operators.  

The main difficulty of Problem \ref{prob:M3} is that it is not clear at all what a quantization of a system of PDEs is. There are different viewpoints on what quantization may mean in this context. Quantization of certain integrable systems of PDEs has been studied in the literature; there exist several quite complicated and seemingly non-equivalent approaches, in which the authors of the present note are not experts. 

Generally, quantization may involve:
\begin{itemize} 
\item Replacing Poisson brackets with commutators.
\item Finding commuting quantum analogs of classical conserved quantities.
\item Possibly, understanding the spectrum of the quantum Hamiltonian and its implications. 
\end{itemize} 

The idea is to consider the sequence of quantizations of the finite-dimensional systems appearing in the procedure of finite-dimensional reductions and to understand whether, and in which sense, they may converge to a quantization of a BKM system. Clearly, this problem is related to Problem \ref{prob:M2}, which is a ``classical'' analog of Problem \ref{prob:M3}. As in Problem \ref{prob:M2}, one should also try to understand how and whether other mathematical objects related to the quantization, in particular the common set of eigenfunctions, converge to something reasonable. So a partial goal is to describe the set of common eigenfunctions of the quantizations of the finite-dimensional systems appearing in the procedure of finite-dimensional reductions. An additional difficulty may arise from the fact that the background metric is non-Riemannian and, in fact, has split signature. 

Note that in the KdV case, there are also other quantization methods suggested in the literature; as mentioned above, the authors of this paper are not experts. 
According to the literature, there exist three main approaches:

\begin{itemize} 
\item \textbf{Quantum Inverse Scattering Method (QISM)}
    \begin{itemize} 
    \item Extends the classical inverse scattering method to quantum systems.
    \item Used to quantize many integrable systems, including KdV, especially in 1+1 dimensions.
    \end{itemize}
    
\item \textbf{Virasoro Algebra and Conformal Field Theory (CFT)}
    \begin{itemize} 
    \item In certain contexts (e.g., Gelfand–Dikii theory), the quantum KdV hierarchy is related to the Virasoro algebra in 2D conformal field theory.
    \item The quantum KdV Hamiltonians can be constructed using modes of the energy-momentum tensor in a CFT ~\cite{BLZ1996}.
    \end{itemize} 
    
\item \textbf{Free Fermion / Bosonisation Techniques}
    \begin{itemize} 
    \item In some formulations, the KdV system can be expressed using bosonised fermionic fields, allowing for an explicit quantization.
    \end{itemize}
\end{itemize}

Note that links with Nijenhuis geometry are expected.

\section{ Haantjes geometry }

Let $L$ be an operator and $\mathcal{N}_L$ its Nijenhuis torsion. 
The following recursion formula defines the {\it Haantjes torsion} of level $m$ (see \cite{Kosmann,TT2021, TempestaTondo2022higher}):
\begin{equation}\label{haam}
\begin{aligned}
\mathcal T^{(1)}_L(\xi, \eta) = & \mathcal N_L(\xi, \eta), \\
\mathcal T^{(m)}_L(\xi, \eta) = & L^2 \mathcal T^{(m - 1)}_L(\xi, \eta) + \mathcal T^{(m - 1)}_L(L\xi, L\eta) - \\
& - L\mathcal T^{(m - 1)}_L(L\xi, \eta) - L\mathcal T^{(m - 1)}_L(\xi, L\eta), \quad m = 2, 3, \dots.   
\end{aligned}    
\end{equation}
For $m=2$, this formula coincides with the classical definition of Haantjes torsion $\mathcal H_L = \mathcal T^{(2)}_L$ introduced in \cite{nijenhuis,Nijenhuis1955, Haantjes}. That is, the Haantjes torsion is given by  
\begin{equation}\label{haam1}
\begin{aligned}
\mathcal H_L(\xi, \eta) = & L^2 \mathcal N_L(\xi, \eta) + \mathcal N_L(L\xi, L\eta) - \\
& - L\mathcal N_L(L\xi, \eta) - L\mathcal N_L(\xi, L\eta).    
\end{aligned}    
\end{equation}
Since $\mathcal N_L$ is a tensorial object, $\mathcal H_L$ is also a tensorial object. We see that the formula for the Haantjes torsion is linear in the derivatives of $L$, with coefficients which are homogeneous polynomials of degree three in $L$.

 We suggest defining Haantjes geometry by mimicking the definition of Nijenhuis geometry: a Haantjes structure is an operator with vanishing Haantjes torsion. Haantjes geometry studies such operators and their natural relations to other objects.   

There exists, however, an earlier definition of Haantjes geometry, due to F. Magri \cite{Magri2014}, as follows.
We give it for completeness, but will not use it. This definition is not equivalent to the above one. 

Let $M$ be an $n$-dimensional smooth manifold and let $A\in C^\infty(M)$ be a function. We consider its differential $dA$. 
Assume we are given $n$ operators  
\[
K_1=\Id,\quad K_2,\dots,K_n,
\]
which commute pairwise,
\[
K_jK_\ell = K_\ell K_j,\qquad j,\ell=1,\dots,n.
\]
The triple $(M,\{K_j\}_{j=1}^n,dA)$ is called a \emph{Haantjes manifold} if, for all
$j,\ell=1,\dots,n$, the $1$-forms
\[
K_jK_\ell\, dA
\]
are closed (equivalently, locally exact). In other words, there exist local functions
$A_{j\ell}\in C^\infty(M)$ such that
\[
K_jK_\ell\, dA = dA_{j\ell},\qquad j,\ell=1,\dots,n.
\]
The symmetric matrix $(A_{j\ell})$ is then called the \emph{matrix potential}.

We note that this definition of Haantjes structure does not require the Haantjes torsion of the operators to vanish at all. Later, in his paper \cite{Magri2018}, F. Magri often assumed the vanishing of Haantjes torsion as an additional assumption. Though a Haantjes structure in the sense of F. Magri is an interesting notion, we will stick to our definition: a Haantjes structure is a choice of an operator with vanishing Haantjes torsion.

Haantjes operators in our definition enjoy the following nice property: any polynomial, with coefficients which may depend on the position, in a Haantjes operator is again a Haantjes operator, by \cite{Bogo2004}.   
Starting with one gl-regular Haantjes operator, one then produces the whole algebra (with respect to pointwise multiplications and taking linear combinations with functional coefficients). Such algebras are called \emph{cyclic}. They are automatically commutative. We refer to \cite{NTT2021,TT2021} for the definition (e.g., \cite[Definition 7]{NTT2021}) and for a discussion of Haantjes algebras, in particular of semi-simple commutative Haantjes algebras, and their interplay with (multi)separation of variables.

It was shown by us in 
\cite{BKMHaantjes2026} that, locally near every
algebraically generic point, a gl-regular Haantjes
operator is a polynomial, with smooth functional
coefficients, in a gl-regular Nijenhuis operator.
Moreover, the latter can be chosen to have constant
coefficients in suitable local coordinates.
Consequently, every cyclic Haantjes algebra generated
by a gl-regular operator is, locally near almost
every point, generated by a gl-regular Nijenhuis
operator.

 It would be interesting to extend the above statement to the case when all elements $K$ of the Haantjes 
algebra are compatible with the symplectic form in the sense that for all vectors $X,Y\in T_pM$ we have 
\begin{equation}\label{eq:compatible}
          \omega(X, K Y)=\omega(KX, Y).
        \end{equation}

\begin{Problem} \label{prob:h1}
   Let $(M^{2n}, \omega)$ be a symplectic manifold. 
   Suppose we have a Haantjes operator $K$, compatible with $\omega$ in the sense of   \eqref{eq:compatible},  such that the operators $\Id, K, K^2, \dots, K^{n-1}$ are linearly independent at almost every point. 

   Is it true that in a neighborhood of almost every point there exists a Nijenhuis operator $L$, compatible with $\omega$, such that $K$ is a polynomial in $L$ (with coefficients which may depend on the position).
\end{Problem}

The semisimple case of Problem~\ref{prob:h1} is essentially covered by
results of Tempesta and Tondo
\cite[Propositions~37--38]{TT2022}.
They showed that a semisimple Abelian Haantjes algebra
admits, locally, a Nijenhuis generator. In the
symplectic-Haantjes setting, this generator can also
be chosen compatible with the symplectic form.
Thus, Problem~\ref{prob:h1} is answered positively when $K$
has $n$ distinct eigenvalues. The main open question
concerns the non-semisimple case.

Consider a Haantjes algebra generated by $K$ as in Problem \ref{prob:h1}: we assume that the operator is compatible with the symplectic form $\omega$, and that 
$$\Id, K, K^2, \dots, K^{n-1}$$ are linearly independent at almost every point. Next, suppose there exists a function $f$ such that it is a common conservation law for all the operators $K_i$, which, as we know, locally means that there exist functions $f_1=f, f_2,\dots, f_n$ such that  
\begin{equation} \label{eq:co} 
df_i= K^{*i-1}df.  
\end{equation}
Then, from the compatibility condition we immediately have that the functions $f_i$ Poisson-commute with respect to the symplectic structure (see also \cite{TT2022}).  

\begin{Problem} \label{prob:h2}
Solve Problem \ref{prob:h1} under the additional assumption that there exist functions $f=f_1, f_2,\dots, f_n$ such that their differentials are linearly independent and such that the condition \eqref{eq:co} is fulfilled. 
\end{Problem}
\weg{
\begin{Problem} \label{prob:h3}
   Let $(M^{2n}, \omega)$ be a symplectic manifold. 
   Suppose we have a Haantijes operator $K$ such that the operators $\Id, K, K^2, \dots, K^{n-1}$ are linearly independent at almost every point, and which is compatible with $\omega$ in the sense \eqref{eq:compatible}. 
Assume in addition that there exist functions $f=f_1, f_2,\dots, f_n$ such that their differentials are linearly independent and such that the condition \eqref{eq:co} is fulfilled.   Then, prove or  disprove the following statement:  the operators $K_i$ are mutual symmetries (the definition of symmetry is e.g. in \cite{nij4, nijapp5}. Symmetries are PDE-analogs of commuting vector fields.) 
 
\end{Problem} 

In relation to the last problem,}
We would like to attract attention to 
the construction of non-semisimple cyclic Haantjes algebras such that their generators have a common conservation law, see 
\cite{BKM26a}. In this construction,  symmetries and conservation laws are put in a relation.   

In our opinion, the study of families of commutative operators with vanishing Nijenhuis or Haantjes torsion is an interesting and promising  topic, generally and in the symplectic or Poisson setting.

\section{Maximally superintegrable systems}  
\subsection{Algebras of maximally superintegrable systems}

We work on an $n$-dimensional manifold and consider the Hamiltonian system on the cotangent bundle 
such that the Hamiltonian is the sum of kinetic, magnetic, and potential terms: 
\begin{equation} \label{eq:ham}  H=K+ V + U. \end{equation}
Recall that the kinetic term $K$ is quadratic in momenta. The magnetic term $V$ is linear in momenta. The coefficients in the quadratic and linear terms may depend on the point $x$. 
The potential $U$ is a function on $M$. We will always assume that the bilinear form corresponding to $K$ 
is nondegenerate, so $K$ is the kinetic energy corresponding to a metric. The magnetic and potential terms can vanish and are assumed to be zero for some problems.

In dimension $n$, we say that the Hamiltonian system is \emph{maximally superintegrable} if its flow admits $2n-1$ functionally independent integrals which are polynomial in momenta.

Let us first consider the case when $V$ and $U$ vanish. In this case, the following result holds.  
\begin{Theorem}[\cite{Matveev2025}] \label{thm:222} Let $(M^n, g)$ be a connected $C^\infty$-smooth pseudo-Riemannian manifold. Assume its geodesic flow is maximally superintegrable and denote by $F_1=H, F_2, \dots, F_{2n-1}$ its functionally independent integrals polynomial in momenta. 

Then, the Poisson bracket of any two of them is algebraically dependent on $F_1, \dots, F_{2n-1}$, in the sense that for every $i,j$
there exists a polynomial $P$ of $2n$ variables which nontrivially depends on the last variable such that 
$$P( F_1,F_2,\dots, F_{2n-1}, \{F_i,F_j\})\equiv 0.$$    
\end{Theorem}

Note that the functional dependence of $\{F_i ,F_j\}$ on the functions 
$
F_1,F_2,\dots, F_{2n-1}
$ is almost trivial; the difficulty is to show that functional dependence is in fact an algebraic one.  
Note also that in dimension $2$, superintegrability is necessarily maximal, as $2n-1= 4-1 =3$, so Theorem \ref{thm:222} can  be applied in the two-dimensional superintegrable setup, where many problems are not solved. 

\begin{Problem} Generalize Theorem \ref{thm:222} to Hamiltonians of the form \eqref{eq:ham}.      
\end{Problem}
Of course, the terms of highest order in momenta of the integrals are integrals for the kinetic energy. It is not clear (and is possibly wrong) whether the polynomial relation survives for the magnetic terms and potential energies. We expect, though, a positive answer.

The case with no magnetic terms could be easier and is definitely better studied. Let us consider the answer in the case 
when the dimension is $2$ (and $V=0$), in the simplest situation where one of the 
integrals is linear in momenta and the other is linear, quadratic, or cubic.

The simplest case, of course, is when both integrals are linear. In this case the metric has constant curvature and the potential energy is constant.   

The case when one integral is linear and the other is quadratic is known under the name ``Darboux-superintegrable case'' and was studied, e.g., in \cite{KKMW,MZZ}. In this case, the polynomial relation for the higher-order parts (i.e., the polynomial $P$ from Theorem \ref{thm:222}) also works for systems with potentials.

In the case when one integral is linear and the other is cubic, the situation is as follows: from the results of \cite{sewa} it follows that in the nontrivial case the corresponding potential energy must be constant. Note that in this case we view the case when the cubic integral is the product of a quadratic integral and a linear one as a trivial case.

The situation in the presence of a nontrivial magnetic
field is less developed. Nevertheless, a number of
interesting examples and classification results are known.

Let us first consider a Riemannian surface of constant
Gaussian curvature $K$, equipped with the magnetic
form proportional to the volume form with a constant coefficient. 

The corresponding magnetic geodesic flow is maximally
superintegrable. Indeed, the three Killing vector
fields of the metric give rise to three first integrals
linear in momenta.  This construction works locally
for all three constant-curvature geometries: the
sphere, the Euclidean plane and the hyperbolic plane.

By
\cite{AgapovDemina2027},  if  the magnetic geodesic flow of a
Riemannian metric $g$ on a two-dimensional manifold
admits two additional  first integrals     linear in momenta, 
  then  we are in the situation as above: 
 the Gaussian curvature of $g$ is constant
and the magnetic form is proportional to the volume form with a constant coefficient. 

A higher dimensional  example of superintegrable magnetic geodesic flows can be constructed   
as follows: take 
the product of two Euclidean planes with the same
constant magnetic field $B$. Denote
\[
    P_i=p_{x_i}+\frac B2 y_i,\qquad
    Q_i=p_{y_i}-\frac B2 x_i,\qquad i=1,2.
\]
The Hamiltonian
\[
    H=\frac12(P_1^2+Q_1^2+P_2^2+Q_2^2)
\]
has seven functionally independent polynomial first
integrals:
\[
\begin{gathered}
    x_i+\frac{Q_i}{B},\qquad
    y_i-\frac{P_i}{B},\qquad i=1,2,\\
    P_1^2+Q_1^2,\qquad
    P_2^2+Q_2^2,\qquad
    P_1P_2+Q_1Q_2.
\end{gathered}
\]
The additional integral in the second line exists
because the two cyclotron motions have the same
angular frequency. Thus the system is maximally
superintegrable in dimension four. The construction
extends to products of arbitrarily many identical
magnetic Euclidean planes.

In contrast, direct products of the analogous magnetic
flows on spheres or hyperbolic planes are generally not
maximally superintegrable. Indeed, the angular frequency
of a closed magnetic orbit on a surface of curvature $K$
depends on its energy when $K\ne 0$. Hence the motions
on different factors are not generally synchronised,
and, when both individual motions are periodic, the
closure of a generic product orbit is a two-torus.

More complicated superintegrable systems with magnetic
fields have been constructed and studied by Marchesiello,
\v{S}nobl, Winternitz, Kub\r{u} and their collaborators;
see, for example,
\cite{MarchesielloSnobl2018,MarchesielloSnobl2020,
KubuMarchesielloSnobl2023,MarchesielloReyesSnobl2024}.
These examples are three-dimensional. Most have a flat
background metric, but examples on manifolds of
nonconstant curvature are also known. In some, but
not all, examples the magnetic vector potential is
locally gauge-equivalent to the metric dual of a
Killing vector field (see also \cite{BKM_magnetic}
for a discussion of integrability in this setting).

The known constructions of maximally superintegrable
systems with polynomial first integrals typically
involve a nontrivial scalar potential. On the other
hand, examples of magnetic geodesic flows without
a scalar potential are known, but the three-dimensional
examples discussed above are either only minimally
superintegrable in the class of polynomial integrals,
or require nonpolynomial integrals for maximal
superintegrability. This motivates the following
problem.

\begin{Problem}
Construct, in dimension three, maximally
superintegrable magnetic geodesic flows with
nontrivial magnetic fields, without a scalar
potential, and with  first integrals polynomial in momenta.
\end{Problem}

When discussing magnetic geodesic flows without
a scalar potential, we mean Hamiltonians which,
locally, can be written in the form
\[
H=\frac12 g^{ij}(p_i+A_i)(p_j+A_j),
\]
up to an additive constant. This condition
is invariant under magnetic gauge transformations.

\subsection{Superintegrable Tannery metrics on the $2$-sphere} \label{sec:tannery}

If a Riemannian  manifold is closed (=compact without boundary), 
geodesics of a maximally superintegrable metric are closed.  In dimension two 
this implies that the manifold is    the projective  plane with metric of constant curvature or is topologically the sphere.    
A Riemannian metric on $S^2$ with all geodesics closed is called a \emph{Zoll metric} \cite{Zoll1903}. All 
geodesics of a Zoll metric necessarily have the same length; see, e.g., \cite{GG81,Besse}.

Under the assumption that one of the integrals is linear, metrics on the sphere with all geodesics closed are called Tannery metrics and are described in \cite{Besse}. 
We repeat the description below. 
The freedom of the description is a function of one variable. 

\begin{Problem} Which of the Tannery metrics are superintegrable, with the additional integral of degree $\ge 4$? 
\end{Problem}

Local examples of superintegrable metrics admitting a Killing vector field 
 and  an  irreducible  polynomial first integral of arbitrarily high
degree were constructed by Valent \cite{Valent2017}.
For odd degrees, some of these metrics extend
to globally defined Zoll metrics on the sphere
\cite{Valent2021}.

A smooth rotationally symmetric (i.e.\ $\mathrm{SO}(2)$-invariant) Riemannian metric on the
$2$-sphere can be written in polar coordinates $(r,\theta)\in[0,\pi]\times\mathbb{R}/2\pi\mathbb{Z}$ as
\begin{equation}\label{eq:rot}
g \;=\; A(r)^2\,dr^2 \;+\; B(r)^2\,d\theta^2,
\end{equation}
where $A(r)>0$ for all $r\in[0,\pi]$.  Smoothness at the poles $r=0,\pi$
implies the following first-order conditions:   
\begin{equation}\label{eq:smooth}
B(0)=B(\pi)=0,\qquad B'(0)=A(0),\qquad B'(\pi)=-A(\pi).
\end{equation}
(Equivalently, after a reparametrisation of $r$ one may assume $A\equiv 1$ and write
$$g=dr^2+f(r)^2d\theta^2$$ with the corresponding pole conditions on $f$.)

A standard normal form for a \emph{rotationally symmetric Zoll metric} on $S^2$ is
\begin{equation}\label{eq:zoll}
g \;=\; \bigl(1+h(\cos r)\bigr)^2\,dr^2 \;+\; \sin^2 r\,d\theta^2,
\qquad (r,\theta)\in[0,\pi]\times\mathbb{R}/2\pi\mathbb{Z},
\end{equation}
where $h:[-1,1]\to\mathbb{R}$ is smooth and satisfies
\begin{equation}\label{eq:hconds}
h(-x)=-h(x)\ \text{ for all }x\in[-1,1],\qquad h(\pm 1)=0,\qquad 1+h(\cos r)>0\ \text{ for all }r.
\end{equation}
For the normalization \eqref{eq:zoll}, one typically has the common length of prime geodesics equal
to $2\pi$; multiplying the metric by a constant rescales this common length.
The round metric corresponds to $h\equiv 0$, i.e.
\[
g_{\mathrm{round}} = dr^2 + \sin^2 r\,d\theta^2.
\]

Next, in the case when one of the integrals is linear, the procedure described in the previous section can be simplified. Indeed, in this case the Poisson bracket of the additional integral of degree $m$ with the linear integral is again an integral of degree $m$. This gives an endomorphism of the (necessarily) finite-dimensional 
space of integrals of degree $m$. Consider an eigenvalue of this endomorphism; this corresponds to an integral such that its Poisson bracket with the linear integral is proportional to it. This gives an ansatz for the system of PDEs corresponding to an integral of higher degree in the following form:
the coefficients of the integral are essentially functions of one variable (multiplied by a certain fixed function of the other variables). This reduces the system of PDEs for the existence of the integral to a system of ODEs.

This approach is due to  \cite{sewa}. In this paper, a series of examples is constructed. One should combine the approach of \cite{sewa} with the form of Tannery metrics above. The $m=3$ case was partially done in  \cite{sewa}  and finished  in \cite{DMPS2026}, using the approach described above.

\subsection{Non-degenerate quadratic abundant  superintegrable systems and separation of variables on complex projective spaces of different signature and scalar curvature}

A non-degenerate quadratic  abundant (super)integrable system  (see e.g. \cite{KressSchoebelVollmer2023, KressSchoebelVollmer2024}) in dimension $n$ is a system with Hamiltonian of the form
\begin{equation} \label{eq:ham1}  H=K+ U(\alpha_1,\dots, \alpha_{n+2}) \end{equation}
characterised by the following properties:
\begin{enumerate}
\item The potential energy $U$ is a linear combination of $n+2$ linearly independent functions with coefficients $\alpha_1,\dots, \alpha_{n+2}$:
\begin{equation}
U=U_1 \alpha_1 + \dots  + U_{n+2} \alpha_{n+2}
\end{equation}
(one of which can be constant).  
\item There exist $n(n+1)/2$ functions of the form $K_i + U_i(\alpha_1,\dots, \alpha_{n+2})$ such that $K_i$
are quadratic in momenta and linearly independent, and $U_i$ are linear combinations of linearly independent functions with coefficients depending on $\alpha_1,\dots, \alpha_{n+2}$, and such that $K_i + U_i(\alpha_1,\dots, \alpha_{n+2})$ are  integrals of \eqref{eq:ham1} for each choice  of the parameters $\alpha$. 
\end{enumerate}

An example of such a situation comes from separation of variables of the standard metric of $\mathbb{C}P(n)$. In this case, $n$ of the $2n$ separating coordinates are ignorable, so the metric locally looks as follows (see, e.g., \cite{BDMN}):
\begin{equation} \label{eq:coords}
        ds^2 = \sum_{\alpha, \beta= 1}^{n} g_{\alpha \beta }(x) \, dx_\alpha dx_\beta + \sum_{i,j=1}^{n} h_{ij}(x) \, dt_i dt_j.
    \end{equation}
Above, $t_i$ are ignorable coordinates and the $n$-dimensional metric $g$ has constant sectional curvature.

In this case the functions $h^{ij}(x)$  (i.e., the entries of the inverse matrix to $h_{ij}$)
can be taken as the components of the potentials $U_1, \dots, U_{n+2}$. One can show there are $n+2$ linearly independent functions among these functions.

The ignorable coordinates determine a maximal abelian subalgebra of the corresponding isometry algebra 
($SU(n+1)$ in the case of $\mathbb{C}P(n)$). In the case of $CP(n)$, the only maximal abelian subalgebra is the Cartan subalgebra, and the functions $U_1,\dots, U_{n+2}$ are essentially unique. 
The situation in other signatures and curvatures is different, as $SU(p,q)$ has other maximal abelian subalgebras. Moreover, if both $p$ and $q$ are large, the corresponding projective space has abelian subalgebras of higher dimension.

\begin{Problem} Does every non-degenerate quadratic superintegrable system such that the kinetic energy corresponds to a metric of constant curvature ``come'' from a complex projective space and its analogs in other signatures and curvatures?
\end{Problem}

Results in this direction include \cite{BKW1983,BKW1985,WOR1994,ORW1993,ORWZ1990}. 

Note that in the Riemannian signature the complex hyperbolic space $\mathbb{C}H^n$  has $n+2$  maximal abelian subalgebras. If one knows the  number of different nondegenerate quadratic  integrable systems such that the kinetic energy is the hyperbolic metric, a simple comparison  of the numbers can give an answer.

\section{Killing tensors, Frobenius structures, and integrable systems of hydrodynamic type}
\label{sec:killing-frobenius}

Killing tensors, Frobenius structures, and integrable systems of hydrodynamic type
are closely related, although these connections are far from being completely
understood. Let $(\mathsf M, g)$ be a (pseudo-)Riemannian manifold and $\nabla$ be the Levi-Civita connection of the metric $g$. A symmetric tensor $K$ of degree $d$ is a Killing tensor if
$\nabla_{(i_0}K_{i_1\ldots i_d)}=0$; equivalently, it defines a first integral
of the geodesic flow of $g$ which is homogeneous of degree $d$ in momenta.
On the other hand, families of commuting $(1,1)$-tensor fields and
Frobenius multiplications naturally arise in the description of commuting
systems of hydrodynamic type and of quadratic integrals of geodesic flows.
We formulate several problems suggested by recent progress in these directions.

\subsection{Frobenius structures and Hamiltonian operators}

\emph{Duality of operator Frobenius algebras.}
A duality construction for Frobenius algebras of operator fields was introduced
in \cite{KTF_BKMDuality2026}.  If the operators in the original algebra are
mutual symmetries, the operators in the dual algebra have the same property.
This gives a method of constructing new commuting systems of hydrodynamic
type.  The same work \cite{KTF_BKMDuality2026} solves the Eisenhart--St\"ackel problem for arbitrary
Segre characteristic, under the assumption that the corresponding
$(1,1)$-tensor fields are mutual symmetries.

\begin{Problem}
Study the duality of operator Frobenius algebras in concrete integrable
hierarchies.  Which hydrodynamic-type systems admit a nontrivial duality,
and how are their conservation laws, Hamiltonian structures, and solutions
related?  Classify, at least in small dimensions, operator Frobenius
algebras which are isomorphic to their duals.
\end{Problem}

\emph{Hesse--Frobenius structures and superintegrability.}
Commutative associative multiplications  compatible with a metric occur in
the theory of nondegenerate abundant second-order superintegrable systems.
Vollmer \cite{KTF_Vollmer2026} established a correspondence, under 
appropriate additional compatibility conditions, between abundant
structures on flat manifolds of dimension at least three and certain
Frobenius-type structures.  This brings the classification of finite-dimensional
commutative algebras into the study of superintegrability.

\begin{Problem}
Use the classification of commutative Frobenius algebras to classify the
corresponding Hesse--Frobenius structures and abundant superintegrable
systems in low dimensions, in particular dimensions $3$, $4$, and $5$.
Can analogous correspondences be established for metrics of nonzero or
nonconstant curvature?  Which algebraic structures can actually arise
from superintegrable Hamiltonians?
\end{Problem}

\emph{Compatible Hamiltonian operators.}
In the theory of non-homogeneous Poisson brackets of the form
$\mathcal P_3+\mathcal P_1$, the compatibility equations are closely related
to Frobenius algebras and pencils of such algebras.  For a nondegenerate
Darboux--Poisson third-order part, a substantial classification was obtained
in \cite{nijapp3}.

\begin{Problem}
Extend the algebraic description of compatible Hamiltonian operators
beyond the nondegenerate Darboux--Poisson setting.  Describe the
corresponding Frobenius and Hesse--Frobenius structures, especially in
low dimensions, and determine which of the resulting compatible pairs
generate new integrable hierarchies.
\end{Problem}

\emph{Generalised Frobenius manifolds.}
New families of generalised Frobenius manifold structures on orbit spaces
of affine Weyl groups have been constructed in
\cite{KTF_JiangEtAlI2025,KTF_JiangEtAlII2026}.  Their relationship with
other Frobenius-type structures, Killing tensors, and Hamiltonian
integrable systems is still largely unexplored.

\begin{Problem}
Determine which of these generalised Frobenius manifold structures admit
an associated Hesse--Frobenius structure or a natural family of Killing
tensors.  Can the resulting geometric data be used to construct or
understand hierarchies such as Volterra, Ablowitz--Ladik, or
$q$-deformed KdV?
\end{Problem}

\subsection{Killing tensors, hydrodynamic-type systems, and separation}

\emph{Killing tensors and hydrodynamic-type equations.}
The Killing equation is an overdetermined system for the metric and the
coefficients of a polynomial first integral.  In dimension two, Bialy and
Mironov \cite{KTF_BialyMironov2011} interpreted the existence of a polynomial
integral of arbitrary degree as a rich quasilinear system of hydrodynamic
type.  This approach suggests studying both the equation itself and its
duality and singularities using Nijenhuis geometry.

\begin{Problem}
Study hydrodynamic-type systems obtained from Killing equations in
higher dimensions or for several commuting Killing tensors.  When do
these systems carry a natural Nijenhuis or operator Frobenius structure?
Describe the behaviour near singular points and the implications for
the existence of polynomial first integrals.
\end{Problem}

\emph{Third-order Hamiltonian operators.}
Homogeneous third-order Hamiltonian operators admit a projectively
invariant description in terms of quadratic line complexes and special
Killing tensors associated with the flat affine connection.  The
classification of such operators has been carried out in dimensions
up to four \cite{KTF_FerapontovPavlovVitolo2016}.

\begin{Problem}
Describe the special Killing tensors arising from homogeneous third-order
Hamiltonian operators in higher dimensions, and classify the associated
operators under projective equivalence.  Can the methods of prolongation
and representation theory be used to decide which Killing tensors belong
to this class?
\end{Problem}

\emph{Benenti--St\"ackel systems.}
A Benenti structure consists of a metric and a compatible self-adjoint
Nijenhuis operator, and gives a family of commuting quadratic Killing
tensors.  At regular points, the eigenvalues of this operator provide separation
coordinates.  The relation between non-diagonal operator fields and
quadratic integrability was substantially clarified in
\cite{KTF_BKMDuality2026}.

\begin{Problem}
Describe local normal forms of Benenti-type pairs $(g,L)$ near singular
points of $L$, especially where the Segre characteristic changes.
What happens to separation coordinates, the Killing tensors, and the
Liouville foliation across such singularities?  Investigate the
corresponding questions in dimensions three and four.
\end{Problem}

\emph{Finite-dimensional reductions beyond BKM.}
As discussed in the previous sections, certain finite-dimensional
Benenti--St\"ackel systems provide finite-gap-type solutions of BKM
hierarchies \cite{BKMreduction}.  Their singular invariant
fibres may encode soliton or other special limits of the associated PDEs.

\begin{Problem}
Compare the known finite-gap and finite-dimensional reductions for
other classical integrable hierarchies, such as NLS, mKdV, and
sine--Gordon, with the Benenti--St\"ackel picture above.  Which of
these reductions admit a description by compatible Nijenhuis operators
and Killing tensors?  Relate singular points of such operators and
singular Liouville fibres to special solutions of the PDEs.
\end{Problem}

\subsection{Killing tensors on symmetric spaces and Lie groups}

A Killing tensor is called \emph{decomposable} if it is a sum of symmetric
products of Killing vector fields; otherwise it is indecomposable.
On spaces of constant curvature \cite{Takeuchi1983, Thompson1986}, and on complex projective spaces \cite{KTF_Eastwood2023}
every Killing tensor, in
all degrees,  is decomposable.  However, Matveev and Nikolayevsky
\cite{KTF_MatveevNikolayevsky2024} found indecomposable quadratic
Killing tensors on quaternionic projective spaces and on the Cayley
plane.  Their later work describes all quadratic Killing tensors on
rank-one Riemannian symmetric spaces
\cite{KTF_MatveevNikolayevsky2026}; see also the prolongation approach
of Eastwood and Leistner \cite{KTF_EastwoodLeistner2026}, which gives
further examples, including $E_6/F_4$.  More recently,
Dearricott and Nikolayevsky \cite{KTF_DearricottNikolayevsky2026}
showed that on quaternionic projective spaces the entire algebra of
Killing tensors is generated by Killing vector fields and certain
quadratic Killing tensors. Note that in view of  \cite{KTF_MatveevNikolayevsky2026,MN25} it is sufficient to study this problem for irreducible spaces only. 

\begin{Problem}
Describe the graded algebra of Killing tensors on the remaining
irreducible Riemannian symmetric spaces.  In particular, determine
which spaces possess indecomposable Killing tensors, identify the
minimal degrees of generators, and find invariant or prolongation-based
criteria for decomposability.  Extend the known quadratic results to
higher degrees, especially on higher-rank and exceptional symmetric
spaces.
\end{Problem}

The situation is also interesting for Lie groups with invariant metrics.
For compact classical Lie groups with a bi-invariant Riemannian metric,
quadratic Killing tensors are decomposable
\cite{KTF_MatveevNguyenNikolayevsky2026}.  A classification of
left-invariant Killing tensors on almost abelian Lie groups has been
obtained by Berto-Cuevas, del Barco and Moroianu
\cite{KTF_BertoCuevasDelBarcoMoroianu2026}.  Further higher-rank
symmetric-space cases were treated by Nguyen and Nikolayevsky
\cite{KTF_NguyenNikolayevsky2026}.

\begin{Problem}
Describe Killing tensors of arbitrary degree on Lie groups with
left-invariant metrics, with particular attention to nilpotent and
solvable groups.  Determine when invariant Killing tensors are generated
by Killing vector fields and the metric, and when genuinely new
higher-degree generators occur.  Compare these questions with the
corresponding classification problems for homogeneous spaces.
\end{Problem}

\subsection{Haantjes operators  and integrability of hydrodynamic-type systems.}
A series of papers by Lorenzoni and collaborators develops the description
of integrable systems of hydrodynamic type
$u_t=A(u)u_x$ by means of $F$-manifolds, i.e., manifolds equipped
with a commutative associative multiplication satisfying the
Hertling--Manin condition \cite{HertlingManin1999}. Lorenzoni, Perletti and van Gemst studied
regular, possibly non-semisimple $F$-manifolds with compatible
connections, and extended semi-Hamiltonian techniques and the
generalised hodograph method \cite{Tsarev1991}
to non-diagonalisable systems
\cite{KTF_LorenzoniPerlettiVanGemst2025,
KTF_LorenzoniPerlettiVanGemstSemiHam2026,
KTF_LorenzoniPerlettiVanGemstHodograph2026}.
Arsie and Lorenzoni subsequently developed the theory of cyclic
$F$-manifolds, which also includes certain non-regular systems
\cite{KTF_ArsieLorenzoniCyclic2026}.

Most recently, Arsie and Lorenzoni
\cite{KTF_ArsieLorenzoniHaantjes2026} proved a conjecture of
Bolsinov, Konyaev and Matveev
\cite{KTF_BKMRiemannInvariants2026}. Namely, if
$K_1,\ldots,K_n$ are pointwise linearly independent mutual symmetries
of a hydrodynamic-type system and a constant linear combination of
them is $\gl$-regular at a point, then the Haantjes torsion
of every $K_i$ vanishes in a neighbourhood of that point. As a
consequence, near an algebraically generic point, an integrable system
with a $\gl$-regular operator can be represented as
$u_t=X\circ u_x$, where $\circ$ is an $F$-manifold multiplication and $X$ is a vector field.
No invariant metric is required in this construction, in contrast to
Frobenius or Hesse--Frobenius structures.

\begin{Problem}
Extend the results above to the  non-cyclic case. What are  the conditions on the system of hydrodynamic type that imply its integrability in the sense that one  can reduce it to a linear PDE system and then use the generalized hodograph method? 
\end{Problem}

\section{Towards the classification of $\gl$-regular Nijenhuis operators in dimension 3 (and higher)}

Let $L$ be an operator  (= field of endomorphisms) on a manifold $\mathsf M$.  Recall that a point $\mathsf p \in\mathsf M$ is called {\it algebraically generic} if the algebraic type of $L\mathsf (\mathsf q)$ is the same for all points $\mathsf q$  from a certain neighbourhood of $\mathsf p$. Otherwise, $\mathsf p$ is {\it singular}.

The splitting theorem  \cite[Theorem~3.1]{nij1} says that locally near every point $\mathsf p$, generic or singular,  a Nijenhuis operator $L$ splits into a direct product of operators having either a single real eigenvalue or a single pair of complex conjugate eigenvalues at $\mathsf p\in\mathsf M$. 

In the context of local normal forms for Nijenhuis operators, there are two major and quite different problems.  The first is a local description of an algebraically generic operator $L$
whose spectrum consists of one real number or two complex conjugate numbers. This problem is solved for $\gl$-regular operators, i.e. for real or complex Jordan blocks. In particular, a description of nilpotent Nijenhuis operators with a fixed structure of Jordan blocks is an open and quite difficult problem. It is known to be related to the  classification of flags of distributions with certain properties, but the knowledge on this subject is very limited.  
Some progress towards understanding these questions was made
in \cite{BKM-JordanChevalley}. In particular, in dimensions
three and four, tensorial criteria were obtained for the
existence of local coordinates in which an operator field
with a single eigenvalue and a single Jordan block takes
an upper triangular form. These results apply, in
particular, to Nijenhuis operators of this algebraic
type. However, they do not provide a classification
of such operators up to local diffeomorphisms,
especially near singular points.

The other problem is a description of Nijenhuis operators at singular points. This problem seems to be much harder and out of reach in full generality. However, there are  two special cases where we might expect some progress.

Let $$\chi (t) = \det (t\cdot \Id - L)=t^n +\sigma_1(x) t^{n-1} + \sigma_2(x) t^{n-2} + \dots + \sigma_n(x)$$ be the characteristic polynomial of a Nijenhuis operator $L$. Then the following relation holds  \cite{nij1}:
\begin{equation}
\label{eq10.1}
J L  = \begin{pmatrix}  \!\!\! -\sigma_1 & 1 & & & \\
 \!\!\! -\sigma_2 & 0 & 1 &  & \\
\vdots &  & \ddots & \ddots &  \\
-\sigma_{n-1} & & & 0 & 1 \\
 \!\!\! - \sigma_n & & & & 0   \end{pmatrix} J, \qquad\mbox{where } J =
 \begin{pmatrix}
 \frac{\partial \sigma_1}{\partial x_1} & \dots &  \frac{\partial \sigma_1}{\partial x_n} \\
 \vdots & \ddots & \vdots \\
 \frac{\partial \sigma_n}{\partial x_1} & \dots &  \frac{\partial \sigma_n}{\partial x_n}
 \end{pmatrix}
 \end{equation}

The other assumption which is expected to simplify the {\it local normal form} problem is $\gl$-regularity.  We deal with the real-analytic case.

\begin{Problem}
\label{prob:b17}
Describe local normal forms for $\gl$-regular operators in dimension 3 (and higher).
\end{Problem}

The classification is known in dimension 2. The most essential part of the classification is related to the case when the coefficients of the characteristic polynomial are functionally independent almost everywhere  (i.e. $L$ is differentially non-degenerate at a generic point).  In this situation, the operator $L$ can be naturally reconstructed from the coefficients of the characteristic polynomial of $L$. In dimension 2, the answer is given in terms of two functions $f(x,y)=\tr L$ and $g(x,y)=\det L$ in a certain coordinate system.  This coordinate system is chosen in such a way that $f$ and $g$ take the simplest possible simultaneous form (this form is then considered, by definition, as a normal form)

One  possible approach in higher dimensions is based on one important fact and one important principle.

{\bf Fact.} 
Consider a real-analytic $\gl$-regular operator $L$ with characteristic polynomial 
$$
\chi_L (\lambda) = \det (\lambda{\cdot}\operatorname{Id} - L) = \lambda^n - f_1 \lambda^{n - 1} - \dots - f_n
$$ 
for $n \geq 2$ in a sufficiently small neighbourhood of a point $\mathsf{p}\in \mathsf{M}$.  Then the following are equivalent:
\begin{enumerate}
\item[{\rm (i)}] $L$ is Nijenhuis.
    \item[{\rm (ii)}] There exists a local coordinate system $x=(x^1, \dots, x^n)$ in which $L$ takes the following form
    \begin{equation}\label{first}
  \Lfirst (x) =  \left(\begin{array}{ccccc}
     f_1 & 1 & 0 & \dots & 0  \\
     f_2 & 0 & 1 & \ddots & \vdots \\
     \vdots & \vdots & \ddots & \ddots & 0 \\
     f_{n - 1} & 0 & \dots & 0 & 1 \\
     f_n & 0 & \dots &  0 & 0 \\
    \end{array}\right),
    \end{equation}
    where $f_i=f_i(x)$ are coefficients of the characteristic polynomial in this coordinate system. 
    These coefficients satisfy the following system of PDEs: 
    \begin{equation}\label{first_set}
    \begin{aligned}
        & \frac {\partial f_i} {\partial x^{j}} = f_i  \frac{\partial f_1}{\partial x^{j+1}} + \frac{\partial f_{i + 1}} {\partial x^{j+1}} , \\
        & \frac{\partial f_n}{\partial x^{j}}  = f_n \frac{\partial f_1}{\partial x^{j+1}}.
        \end{aligned}
    \end{equation}
for $1 \leq i, j \leq n-1$. 

This system can be rewritten in an elegant matrix form...
     
       \item[{\rm (iii)}] There exists a local coordinate system $x=(x^1, \dots, x^n)$ in which $L$ takes the following form
    \begin{equation}\label{second}
    \Lsecond(x) = \left(\begin{array}{ccccc}
     0 & 1 & 0 & \dots & 0  \\
     0 & 0 & 1 & \ddots & \vdots \\
     \vdots & \vdots & \ddots &\ddots &  0\\
     0 & 0 & \dots & 0 & 1 \\
     f_n & f_{n - 1} & \dots  & f_2 & f_1 \\
    \end{array}\right),
    \end{equation}
    where $f_i=f_i(x)$ are coefficients of the characteristic polynomial in this coordinate system.  These coefficients satisfy a system of PDEs that can be written in the form
    \begin{equation}\label{second_set}
        \ddd \omega = 0,  \quad \ddd \bigl( L^*\omega \bigr)= 0,
    \end{equation}
where $\omega = f_n \ddd x^1 + \dots + f_1 \ddd x^n$.
    
\end{enumerate}

The system \eqref{first_set} is integrable in the sense that for every initial condition (for the Cauchy problem)
$$
f_i(0,\dots,0,x_n) = v(x_n), \quad i=1,\dots,n,
$$
there exists a unique solution of  \eqref{first_set}.

We continue  the discussion in dimension 3.  We denote the coordinates by $x,y,z$  and the functions $f_1,  f_2, f_3$ by $f$, $g$, $h$. We also assume that at a given point $\mathsf p\in \R^3$, we have $$f(\mathsf p)=g(\mathsf p)=h(\mathsf p)=0,$$ i.e. $L(\mathsf p)$ is a nilpotent Jordan block.

Then the solution can be found in the form of a power series in $x$ and $y$:
\begin{equation}
\label{eq:approx}
\begin{aligned}
f(x,y,z) &= f_0(z) + x f_1(z) + y f_2(z) + x^2 f_{11}(z) + 2xy\, f_{12}(z) + y^2 f_{22}(z) + \dots \\ 
g(x,y,z) &= g_0(z) + x g_1(z) + y g_2(z) + x^2 g_{11}(z) + 2xy\, g_{12}(z) + y^2 g_{22}(z) + \dots \\ 
h(x,y,z) &= h_0(z) + x h_1(z) + y h_2(z) + x^2 h_{11}(z) + 2xy \,h_{12}(z) + y^2 h_{22}(z) + \dots \\  
\end{aligned}
\end{equation}

The coefficients of these power series can be found step by step from the initial conditions $f_0(z)$, $g_0(z)$, $h_0(z)$.  

{\bf Principle.}  Given  a real analytic map $\mathcal F: \R^3 \to \R^3$,  consider its Taylor expansion and truncate it at a certain step. Then the truncation $\widetilde{\mathcal F}$ and the original map $\mathcal F$  are related by some coordinate transformation.  In other words,  the type of a singular point of a real-analytic map $\mathcal F$ is completely defined by its Taylor polynomial of a sufficiently large degree.

If this principle works in our problem, then we can proceed as follows. 

\begin{itemize}
\item[Step 1.]  Solve \eqref{first_set} approximately, i.e.  compute the Taylor expansion \eqref{eq:approx} and truncate it at a certain step.  In this way we obtain three functions $\widehat f(x,y,z)$, $\widehat g(x,y,z)$, $\widehat h(x,y,z)$ that approximate the true solutions 
$f(x,y,z)$, $g(x,y,z)$, $h(x,y,z)$.

\item[Step 2.]  Check that $\widehat f$, $\widehat g$, $\widehat h$ are good enough by constructing the Nijenhuis operator $L$ having $\widehat f$, $\widehat g$, $\widehat h$ as the coefficients of the characteristic polynomial. To that end, we use the standard formula \eqref{eq:fundidentity}.  Here there are two possibilities:

\begin{itemize}

\item  $L$ is smooth;

\item $L$ is not smooth (because of ``division by zero'').

\end{itemize} 

In the second case,  we will have to take a better approximation by truncating the Taylor expansion at a higher step. In the first case, we obtain a non-trivial example of a Nijenhuis operator  (it should be checked that $L$ so obtained is $\gl$-regular).    This will be, at least, a good example of a $\gl$-regular operator with the required properties.  However, 
this $L$ will not be in the first companion form.

\item[Step 3.] Given three good functions $\widehat f$, $\widehat g$, $\widehat h$ explicitly written in local coordinates, we can try to change coordinates to simplify them.  It is expected that one will be able to substantially simplify these functions so that the final {\it normal} form for 
the coefficients of the characteristic polynomial of $L$ will have reasonably simple normal form. 

\end{itemize}

Important hint:  without loss of generality we may assume that $h_0(z) = z^k$ and moreover $h(z)=z^k$

In fact, the above approach should be considered as a kind of experimental work which will lead us to a number of examples. Almost for sure, as a result we will observe some new important phenomena which will be then used for classification.

\weg{
Some experiments made by Steven Dickerson (final year MMath project) led, in particular, to the following conjecture.

{\bf Conjecture.}  In the above setting assume that the initial conditions for $f$, $g$ and $h$ are as follows
$$
h_0(z) = z^k, \quad g_0(z) =  z + \mbox{higher order terms}, \quad f_0(z) = \mbox{arbitrary}
$$    
Then the corresponding Nijenhuis operators are equivalent to each other and  the normal forms for the coefficients of the characteristic polynomial of $L$ are   (not companion form!)
$$
h(x,y,z)=z^k, \quad g(x,y,z)=z^{k-1}y + z, \quad f(x,y,z)=x.
$$ 

Notice that different initial conditions $f_0, g_0, h_0$ and $\widetilde f_0, \widetilde g_0, \widetilde h_0$ may produce operators $L$ and $\widetilde L$ of the form \eqref{first} that are equivalent to each other (i.e. can be transforms to each other by a suitable change of coordinates).  }

Some progress towards Problem~\ref{prob:b17} was recently made
by Akpan and Manenti \cite{AkpanManenti2026}.
In particular, they proved a conjecture of Bolsinov,
suggested by experiments of Steven Dickerson.

More precisely, suppose that the initial data in
\eqref{eq:approx} satisfy
\[
f_0(z)\text{ arbitrary},\qquad
g_0(z)=z\widehat g(z),\quad \widehat g(0)\ne 0,
\]
\[
h_0(z)=z^k\widehat h(z),\quad
k\ge 2,\quad \widehat h(0)\ne 0.
\]
Then there exist local coordinates in which the
coefficients of the characteristic polynomial take
the normal form
\[
f=x,\qquad
g=z(\alpha+y),\qquad
h=\varepsilon z^k,
\]
where $\alpha\ne 0$ and $\varepsilon=\pm1$.
The corresponding normal form of the Nijenhuis
operator is also given explicitly in
\cite{AkpanManenti2026}.

The same paper obtains normal forms in several
additional degenerate cases, including certain
cases where $\det L\equiv 0$.
A complete classification of three-dimensional
$\mathrm{gl}$-regular Nijenhuis operators,
however, remains open.

\begin{Problem}
\label{prob:b18}
Under what necessary and/or sufficient conditions do two sets of functions $f_0(z), g_0(z), h_0(z)$ and $\widetilde f_0(z), \widetilde g_0(z), \widetilde h_0(z)$ lead to equivalent Nijenhuis operators?
\end{Problem}

Such conditions definitely exist.  This means that the classification of $\gl$-regular operators can be done in terms of three functions of one variable by means of partitioning the set of such triples $\{ f_0(z), g_0(z), h_0(z)\}$ into some special equivalence classes.  The point is that solving \eqref{first_set} for given initial conditions is a very difficult problem, and partition into equivalence classes might be easier.

\weg{
\begin{Problem}
\label{prob:b19}  
Let $\mathsf{Sing}\subset \mathsf M$ be the set of singular points of a Nijenhuis operator $L$ and $\mathsf{Crit}\subset \mathsf M$ the set of critical points of the characteristic map $\mathcal F = (f,g,h)$
(i.e. such that $\operatorname{rank} \dd \mathcal F(x)<3$).  What is the relationship between these sets?   

Let $L$ be $\gl$-regular and suppose $L(\mathsf p)$ is a nilpotent Jordan block. In many examples of this kind (all?),  the critical set  $\mathsf{Crit}$ is a smooth hypersurface defined by the equation $\det L=0$ (as already mentioned, without loss of generality, one can assume that in a suitable coordinate system $\det L = z^n$).  Is it a general fact? 
\end{Problem}
}

For some reason, $h=\det L$ seems to be the most important coefficient of the  characteristic polynomial. In particular, if $\dd \det L\ne 0$ at a singular point $\mathsf p\in \mathsf M$, then $\mathsf p$ is {\it differentially non-degenerate} and $L$ reduces to the standard normal form
$$
L = \begin{pmatrix}
x & 1 & 0 \\
y & 0 & 1\\
z & 0 & 0
\end{pmatrix}
$$

Since $\det L = z^k$ in a suitable coordinate system, we can consider $k$ as the order of zero of $\det L$ at a singular point and this number is well defined.   The above matrix gives a normal form of $L$ in the case of zero of order 1.

It seems to be a good idea to start classification with $k=2$  (then $k=3,4,\dots$).  In dimension 2, the case $k=2$ is indeed quite simple.

\begin{Problem}
\label{prob:b20}
Classify all $\gl$-regular operators for $k=2$.
\end{Problem}

In terms of the 3-step approach discussed above,  $k=2$ means that the initial condition for $h$ is $h(0,0,z)=h_0(z)=z^2$.  Consider the next function $g$ and the initial condition $g(0,0,z)=g_0(z)$.  It is natural to assume that the type of singular point will essentially depend on the order of zero of $g_0$ at $z=0$.  If $g_0= z + \mbox{higher order terms}$, i.e., $g_0$ has zero of order 1, then see   the conjecture of Bolsinov, discussed above, according to which the answer will be quite simple. The opposite extremal case is when $g_0(z) \equiv 0$  (zero of infinite order).

\begin{Problem}
\label{prob:b21}
Let $$h_0(z)=z^2\, , \ g_0(z)\equiv 0\,  ,  \ f_0(z)=\mbox{some function}.$$  Treat this case.  How many terms in the Taylor expansion of $\mathcal F$ do we need to get a good approximation?  What is the normal form in this case? How does it depend on $f_0(z)$? 
\end{Problem}   

\begin{Problem}
\label{prob:b22}
Let $$h_0(z)=z^k\ , \ g_0(z)\equiv 0\, ,  \ f_0(z)\equiv 0.$$  Treat this case.  The coefficients of the Taylor series for $\mathcal F$ will depend only on one function $h_0(z)=z^k$.  Can one get a reasonable formula for them? How many terms in the Taylor expansion of $\mathcal F$ do we need to get a good approximation?  What is the normal form in this case? 
\end{Problem}

\section{Classification of left-symmetric algebras and related problems: approach based on Nijenhuis operators}

\begin{Definition}
 A finite-dimensional algebra $(\mathfrak a, \ast)$  is called left-symmetric if the multiplication satisfies the following identity:
 $$
 (u * v) * w - u* (v* w ) = (v * u) * w - v* (u* w ), \quad \forall u,v,w\in\mathfrak a.  
 $$ 
\end{Definition}

The main property of such algebras is that the commutator
$$
[u,v] = u*v - v*u
$$
satisfies the Jacobi identity. Hence $(\mathfrak a,[\, ,\,])$ is a Lie algebra; this is the reason why left-symmetric algebras are also called pre-Lie algebras,  see~\cite{Vinberg1963, Burde2006}.

Geometrically, let $\mathfrak g=(\mathfrak a,[\, ,\,])$ be the associated Lie algebra and let $G$ be a Lie group of $\mathfrak g$. Define a left-invariant connection $\nabla$ on $G$ by prescribing its value on left-invariant vector fields:
$$
\nabla_{u^L} v^L = (u*v)^L .
$$
Here $v^L$ stands for the unique left-invariant vector field, defined by $v$ in the tangent space of the unity of $G$.

The $\nabla$ is torsion-free and flat: the vanishing of the torsion is equivalent to $[u,v] = u*v-v*u$, and the vanishing of the curvature is precisely the left-symmetric identity. Conversely, every left-invariant torsion-free flat affine connection on $G$ arises in this way.

\begin{Problem}
\label{prob:b1}
The goal is to classify LSAs in dimensions  $\ge 3$, over  $\mathbb C$  and/or $\R$. 
\end{Problem}

A classification over $\mathbb{C}$  in dimension 3 was claimed  in the recent  paper \cite{Abdurasulov}, see also Bai~\cite{Bai2009}. The paper \cite{Abdurasulov} contains an  answer  but does not really provide a proof or details of calculations. A classification of three-dimensional real left-symmetric algebras, formulated equivalently as a classification of left-invariant flat torsion-free connections on
three-dimensional real Lie groups, was recently proposed
in \cite{AitAissaElBourkadiMansouri2026}.

Obviously, the approach suggested in the article can be extended to other small dimensions, where the classification of Lie algebras exists. However, it would be interesting to verify this classification
and compare it with the classification over
$\mathbb C$ and $\R$ of the corresponding linear
Nijenhuis operators.

We suggest  addressing this problem, and in partiucular reproving/checking \cite{Abdurasulov,AitAissaElBourkadiMansouri2026},
using   the following natural connection between left-symmetric algebras and Nijenhuis operators whose components depend linearly  on local coordinates $x^1,\dots, x^n$:
\begin{equation}
\label{eq:lsa1}
L = \Bigl(  L^i_j(x) \Bigr), \quad L^i_j =\sum_\alpha  \ell^i_{\alpha j} x^\alpha.
\end{equation}
 
\begin{Theorem}[\cite{nij1,nij2,Winterhalder1997}] \label{t:above}
An operator \eqref{eq:lsa1} is Nijenhuis if and only if  $\Bigl(   \ell^i_{\alpha j}  \Bigr)$ is a structure tensor of a left-symmetric algebra.
\end{Theorem}

Let us explain how we expect it to use, in the simplest  3-dimensional case, though as mentioned above the 3-dimensional version of Problem \ref{prob:b1} is claimed to be solved in  \cite{Abdurasulov, AitAissaElBourkadiMansouri2026}; the general approach is expected to work in all dimensions. 
Theorem \ref{t:above} reduces  the classification problem   to a description of Nijenhuis operators on  $\mathbb C^3(x,y,z)$ or $\mathbb R^3(x,y,z)$ (up to linear change of variables).  Sometimes it will be convenient to represent such an operator in the form 
$$
L = x L_1 + y L_2 + z L_3,
$$
where $L_1, L_2, L_3$  are constant matrices.  Also it might be convenient to consider the linear span of these matrices $\mathcal L = \operatorname{Span}(L_1, L_2, L_3)$ as a vector subspace, associated with our operator $L$, in the space of $3\times 3$ matrices.  Notice that  $L(x,y,z)\in \mathcal L$.

Let us mention some important properties of  $L(x,y,z)$:

\begin{itemize}

\item At a generic point,  $L$ has a certain algebraic type (Segre characteristic).  
\item If $L$ has a constant eigenvalue, then it is equal to zero.
\item  Different non-constant eigenvalues are functionally independent  (in contrast with two previous properties, this one essentially uses the fact  that $L$ is Nijenhuis).
\item  The number of (algebraically=functionally) independent coefficients  $\sigma_1, \sigma_2, \sigma_3$ of the characteristic polynomial  $$\chi_L(t)=\det(t\,\Id - L)=t^3 - \sigma_1t^2 - \sigma_2 t - \sigma_3$$ equals the number of different non-constant eigenvalues.
\end{itemize}

These properties allow us to distinguish the following types of Nijenhuis operators in dimension 3 (depending on the properties of $L$ in a neighbourhood of a generic point).  Below, $\lambda, \mu, \nu$  denote the eigenvalues of $L$ (at a generic point),  we also recall that in our notation:
$$
\sigma_1 = \lambda + \mu + \nu, \quad \sigma_2 = -\lambda\mu - \lambda\nu - \mu\nu, \quad \sigma_3 = \lambda \mu \nu.
$$

\begin{enumerate}

\item[\bf Type 1.]  $\lambda, \mu, \nu$
$$
\mbox{ $\sigma_1,\sigma_2,\sigma_3$ are independent.}
$$

\item[\bf Type 2.]   $\lambda, \mu, 0$
$$
\mbox{$\sigma_1,\sigma_2$ are independent, $\sigma_3=0$.}
$$

\item[\bf Type 3.]  $\lambda, \mu=\lambda, \nu$   (one or two Jordan $\lambda$-blocks)
$$
\mbox{$\sigma_1,\sigma_2$ are independent, whereas $\sigma_3$ depends on them.} 
$$
\begin{itemize}
\item[\bf 3a:] one Jordan $\lambda$-block of size $2\times 2$
\item[\bf 3b:] two trivial Jordan $\lambda$-blocks of size  $1\times 1$
\end{itemize}

\item[\bf Type 4.]  $\lambda, \lambda, \lambda$  (i.e. $\lambda = \mu = \nu\ne 0$) 

In this case $\lambda$ is a linear function,  which can be taken to be the first coordinate  $x$. Then
 $$
 \sigma_1 = 3x, \quad \sigma_2 = -3 x^2, \quad \sigma_3 = x^3.
 $$
There are three subcases:
\begin{itemize}
\item[\bf 4a:] one Jordan $\lambda$-block of size $3\times 3$
\item[\bf 4b:] two  Jordan $\lambda$-blocks of size $2\times 2$ and $1\times 1$
\item[\bf 4c:] three trivial Jordan $\lambda$-blocks of size $1\times 1$
\end{itemize}

\item[\bf Type 5.]  $\lambda, \mu=\lambda, 0$  

In this case $\lambda$ is a linear function,  which can be taken to be the first coordinate  $x$. Then
$$
\sigma_1 = 2x, \quad \sigma_2 = -x^2, \quad \sigma_3 = 0
$$

There are two subcases:
\begin{itemize}
\item[\bf 5a:] one Jordan $\lambda$-block of size $2\times 2$
\item[\bf 5b:] two trivial Jordan $\lambda$-blocks of size $1\times 1$
\end{itemize}

\item[\bf Type 6.]  $\lambda, 0, 0$  (one or two Jordan $0$-blocks)
$$
\sigma_1 = x, \quad  \sigma_2 = \sigma_3 = 0.
$$
There are two subcases:
\begin{itemize}
\item[\bf 6a:] one Jordan $0$-block of size $2\times 2$
\item[\bf 6b:] two trivial Jordan $0$-blocks of size $1\times 1$
\end{itemize}

\item[\bf Type 7.]  $0, 0 , 0$   nilpotent case
$$
\sigma_1=\sigma_2=\sigma_3 = 0.
$$
There are three subcases:
\begin{itemize}
\item[\bf 7a:] one Jordan $0$-block of size $3\times 3$
\item[\bf 7b:] two  Jordan $0$-blocks of size $2\times 2$ and $1\times 1$
\item[\bf 7c:] three trivial Jordan $0$-blocks of size $1\times 1$
\end{itemize}
\end{enumerate}

The general scheme is as follows:
\begin{itemize}
\item[Step 1.] We start with an arbitrary operator $L$ whose components are linear functions in  $x,y,z$. Thus, we have $27 = 9 \times 3$ unknown parameters.  They satisfy (quadratic) Nijenhuis relations which we need to resolve.  We are also allowed to make arbitrary linear coordinate transformations. In other words, we solve the problem modulo the action of  $GL(3,\mathbb C)$.

\item[Step 2.] One of the main tools for description of $L$ will be the fundamental identity:
$$
J L = \begin{pmatrix} \sigma_1 & 1 & 0 \\ \sigma_2 & 0 & 1 \\ \sigma_3 & 0 & 0\end{pmatrix} J, \quad \mbox{where } J = 
\begin{pmatrix}  
\frac{\partial \sigma_1}{\partial x} & \frac{\partial \sigma_1}{\partial y} & \frac{\partial \sigma_1}{\partial z} \\    
\frac{\partial \sigma_2}{\partial x} & \frac{\partial \sigma_2}{\partial y} & \frac{\partial \sigma_2}{\partial z} \\
\frac{\partial \sigma_3}{\partial x} & \frac{\partial \sigma_3}{\partial y} & \frac{\partial \sigma_3}{\partial z} \\
\end{pmatrix}
$$ 

If $\sigma_1, \sigma_2, \sigma_3$ are given,  then this allows us (almost always) to reconstruct some of the components of $L$ and to obtain relations for other components.   To get these relations, one needs to use the fact that  $\sigma_1, \sigma_2,\sigma_3$ are the coefficients of the characteristic polynomial of $L$.  This allows us to substantially reduce the number of parameters.

\item[Step 3.] After this the operator $L$ can be written with a smaller number of parameters, and we next obtain a system of (non-linear) equations on them by writing the Nijenhuis relation.
Resolving them we get an answer with even smaller number of parameters. If this resolution is correct, then for any choice of parameters,
we get one of the left-symmetric algebras we are interested in.

\item[Step 4.]  Some of the remaining parameters can be killed by using an additional change of variables, like scaling.

\end{itemize}

Steps 3 and 4 can sometimes be interchanged, or do first 4, 3 and then 4 again.

The first type was completely treated by Serena Scapucci \cite{Scapucci} and Sofia Degtyareva  et al \cite{Degtyareva, Degtyareva2}.

\subsection{Some related problems from  the `Open problems 2018' paper \cite{BMMT}}

A general question (which would immediately clarify the local structure of singular points for a wide class of Nijenhuis operators)  is as follows.

\begin{Problem}
\label{prob:b2}
Consider a smooth map $\Phi=(\sigma_1,\dots,\sigma_n):  U(0) \to \R^n$, where $U(0)$ is a neighbourhood of the origin $0\in \R^n$.  We assume that $0$ is a singular point of $\Phi$, that is, $\operatorname{rank} \dd \Phi(0) < n$, but almost all points $x\in U(0)$ are regular. The problem is to describe those maps $\Phi$ for which all the components of the Nijenhuis operator $L(x)$ defined by \eqref{eq10.1} are smooth  functions  (notice that $L$ is well defined and smooth at regular points of $\Phi$, and we are interested in necessary and sufficient conditions for $L$ to be smoothly extendable onto the whole neighbourhood $U(0)$).
\end{Problem}

For a given left-symmetric algebra with the structure constants $b^i_{jk}$, consider the coefficients $\sigma_1, \dots, \sigma_n$ of the characteristic polynomial of  $$L(x)= \Bigl(L^i_j = \sum_k b^i_{jk} x_k\Bigr).$$  It is easy to see that $\sigma_k$ is a homogeneous polynomial of degree $k$.   Consider the class of left-symmetric algebras for which these polynomials are algebraically independent.  From \eqref{eq10.1} we conclude that the classification of such left-symmetric algebras is equivalent to that of collections of polynomials $\sigma_1, \dots, \sigma_n$, $\deg \sigma_k=k$, satisfying the following rather exceptional property.

Consider the Jacobian matrix $$J = \left(   \frac{\partial (\sigma_1,\dots, \sigma_n)}{\partial(x_1,\dots, x_n)}\right).$$  It is easy to see that $D=\det J$ is a homogeneous polynomial of degree $\tfrac{n(n-1)}{2}$.   Consider the matrix
$$
L=J^{-1}  \begin{pmatrix}  \sigma_1 & 1 & & \\
\vdots &   & \ddots & \\
\sigma_{n-1} & & & 1\\
\sigma_n & 0 & \dots & 0
\end{pmatrix} J, \quad\mbox{where } J = \left(   \frac{\partial (\sigma_1,\dots, \sigma_n)}{\partial(x_1,\dots, x_n)}\right)
$$
It is easy to check that the entries of this matrix are rational functions of the form $L_{ij}=\frac{Q_{ij}}{D}$ where $Q_{ij}$ are homogeneous polynomials of degree $\deg D + 1$. In other words, the entries of $L$ are homogeneous rational functions of degree one.  For some special polynomials $\sigma_1, \dots, \sigma_n$, a miracle happens:  all $Q_{ij}$ turn out to be divisible by $D$ and the entries of $L$ become linear functions in $x_1,\dots, x_n$.  Such $L$ are automatically Nijenhuis operators (equivalently, LSA's).

\begin{Problem}
\label{prob:b3}
Describe/classify the collections of algebraically independent homogeneous polynomials $\sigma_1, \dots, \sigma_n$, $\deg \sigma_k = k$, in $n$ variables $x_1,\dots, x_n$  with the required property.
\end{Problem}

As mentioned above, in dimension 3 this problem was solved by S. Scapucci \cite{Scapucci} and S. Degtyareva \cite{Degtyareva}.

\begin{Problem} 
\label{prob:b4}
Is it true that the LSAs possessing the above property are non-degenerate?
\end{Problem}

We know only one example of a non-degenerate LSA in arbitrary dimension, namely, the diagonal LSA.  The thing is that we have no  tools to construct them. In this regard, the answer to the following `naive' question could be quite helpful.

\begin{Problem}
\label{prob:b5} 
Is it true that the direct sum of two non-degenerate LSAs is non-degenerate?
\end{Problem}

\section{Characteristic map and its singularities}

Before starting this discussion, let us recall some basic concepts in symplectic topology of integrable systems. We just want to emphasise some analogy between the momentum mapping in the theory of integrable Hamiltonian systems and the characteristic map in Nijenhuis geometry. 

A finite-dimensional integrable system is given by $n$ Poisson commuting functions $f_1,\dots, f_n$ on a symplectic manifold $(\mathsf M^{2n}, \omega)$.   We assume that these functions are functionally independent almost everywhere.

These functions define a map $\mathcal F=(f_1,\dots,f_n) : \mathsf M^{2n} \to \R^n$ called the {\it momentum mapping} of the given integrable system.  The goal of symplectic topology of IHS is basically to study the properties of this map and, more specifically,  the {\it singular Lagrangian fibration} defined by $\mathcal F$.  The key word here is {\it singular}, which means that we are interested in those points $x\in \mathsf  M^{2n}$ where the rank of $\dd\mathcal F$ drops, i.e., becomes smaller than $n$.

The word {\it Lagrangian} is also important and basically amounts to the fact that $f_1,\dots,f_n$ pairwise Poisson commute.  This property makes the singularities of $\mathcal F$ very special.   Recall that at regular points, we have the following fundamental fact.

\begin{Proposition}[Sophus Lie, see e.g.~\cite{LibermannMarle1987}]
\label{prop:b3}
If $\mathsf p\in\mathsf M^{2n}$ is not singular (= not critical), then the functions $f_1,\dots, f_n$ can be reduced to a very simple canonical form. Namely, there exists a coordinate system $q_1,\dots,q_n, \, p_1,\dots, p_n$ such that 
$$
f_i = p_i\quad\mbox{and}\quad \omega = \sum \dd q_i \wedge \dd p_i.
$$ 
\end{Proposition}

The counterpart of the momentum mapping in the context of Nijenhuis geometry is the characteristic map which is defined as follows.

Let $L$ be a Nijenhuis operator on a manifold $\mathsf M^n$ and 
$$
\chi_L(t) = \det (t\cdot\Id - L) = t^n - \sigma_1(x)t^{n-1} - \sigma_2(x)t^{n-2}  - \dots - \sigma_{n-1}(x)t - \sigma_n(x)
$$ 
be its characteristic polynomial.  The characteristic map associated with $L$ is the map defined by the coefficients of the characteristic polynomial, that is,
$$
\Phi = (\sigma_1,\sigma_2,\dots,\sigma_n) : \mathsf M^n \to \R^n.
$$

A point $\mathsf p\in \mathsf M^n$ is called differentially non-degenerate  if it is not critical for the characteristic map $\Phi$, i.e.,
$$
\operatorname{rank}\dd \Phi(\mathsf p) = n.
$$

One of the fundamental problems in Nijenhuis geometry is to study singularities (critical points) of the characteristic map and, in particular, their normal forms.  

The importance of the characteristic map is explained by the following fact.

\begin{Proposition}[\cite{nij1}]
Let the characteristic map $\mathcal F$ be explicitly given in a certain coordinate system $x_1,\dots, x_n$.  Then in this coordinate system the following identity holds:
\begin{equation}
\label{eq:fundidentity}
J L = S J,
\end{equation}
where $L=L(x)$ is the matrix of the Nijenhuis operator in this coordinate system, $$J = \Bigl( \frac{\partial \sigma_i}{\partial x_j}  \Bigr)$$ is the Jacobi  matrix of the characteristic map $\Phi$ and $S$ is the (first) companion matrix of the characteristic polynomial of $L$, i.e.,
$$
S = \begin{pmatrix}
\sigma_1 & 1 & & & \\
\sigma_2 & 0&1 & & \\
\vdots & \vdots & \ddots &\ddots &  \\
\sigma_{n-1} &0 & \dots & 0 & 1 \\
\sigma_n & 0 &\dots & 0 & 0
\end{pmatrix}
$$
\end{Proposition}

This relation means that $L$ can be uniquely reconstructed from the characteristic map by the formula $L = J^{-1} S J$ if almost all points  are non-critical. 

An analog of Proposition \ref{prop:b3} can be stated as follows  (normal form of $\Phi$ at a regular (=non-critical) point).

\begin{Proposition}[\cite{nij1}]
Let $\mathsf p\in\mathsf M^n$ be a non-critical point of the characteristic map $\Phi=(\sigma_1,\dots,\sigma_n)$. Then there exists a coordinate system $
x_1,\dots, x_n$ such that $\sigma_i = x_i$ and
\begin{equation}
\label{eq:Lfirst}
L = \Lfirst =
\begin{pmatrix}
x_1 & 1 & & & \\
x_2 & 0&1 & & \\
\vdots & \vdots & \ddots &\ddots &  \\
x_{n-1} &0 & \dots & 0 & 1 \\
x_n & 0 &\dots & 0 & 0
\end{pmatrix}
\end{equation}
\end{Proposition}

Formula \eqref{eq:fundidentity} shows that,  in the case when the coefficients of the characteristic polynomial of $L$ are functionally independent,   the Nijenhuis operator $L$ can be explicitly reconstructed from $\sigma_1, \dots, \sigma_n$.   Notice that, by continuity, this is still true if $\sigma_1, \dots, \sigma_n$ are independent almost everywhere on $\mathsf M$.  In particular, the singularities of $L$  will be completely determined by the singularities of the smooth map $\Phi=(\sigma_1, \dots, \sigma_n) : M \to \R^n$.  However,  the singularities of $\Phi$ must be very special.

We can illustrate this with a simple example.
Let $\dim M=2$, so that we have only two coefficients
$\sigma_1$ and $\sigma_2$. Consider a point where
$d\sigma_1$ and $d\sigma_2$ are linearly dependent
and $d\Phi$ has rank $1$. Assume, in addition,
that $d\sigma_1\neq 0$, so that there exists a
local coordinate system $(x,y)$ such that
$\sigma_1=x$ and $\sigma_2=f(x,y)$.

Then the identity $JL=SJ$ gives
\[
J=
\begin{pmatrix}
1&0\\
f_x&f_y
\end{pmatrix},
\qquad
S=
\begin{pmatrix}
x&1\\
f&0
\end{pmatrix}.
\]
Hence, wherever $f_y\neq 0$,
\[
{
L=J^{-1}SJ=
\begin{pmatrix}
x+f_x&f_y\\[2mm]
\dfrac{f-f_x(x+f_x)}{f_y}&-f_x
\end{pmatrix}.
}
\]

Thus, $\sigma_2=f(x,y)$ is admissible if and only if
\[
\frac{f-f_x(x+f_x)}{f_y}
\]
admits a smooth extension across the set $f_y=0$.
In other words, to describe singular points of $L$
in this rather simple particular case, we need
to solve the following problem.

\begin{Problem} \label{problem:b6}
Describe all smooth functions $f(x,y)$ of two
variables for which
\[
{
\frac{f-f_x(x+f_x)}{f_y}
}
\]
admits a smooth extension. This question is,
of course, local, and we assume that $f_y$
vanishes at the point under consideration,
but is nonzero on a dense open subset.
\end{Problem}

A more general version of this question is Problem \ref{prob:b2}, which is basically equivalent to 
\begin{Problem}
\label{prob:b8}
Classify singular points  (find local normal forms) for Nijenhuis operators $L$ such that $L$ is differentially non-degenerate almost everywhere.
\end{Problem}

The assumption of {\it differential non-degeneracy almost everywhere} is expected to significantly simplify the local description problem since this problem, as just explained, reduces to the analysis of singularities of a characteristic map $\Phi$.   It may happen that there exists a magic recipe that allows one to effectively construct singularities of $\Phi$ with the required property, but now everything seems to be very obscure. 
However, we could try  (as one does in the theory of integrable systems) to distinguish some special classes of singularities, which are in some sense non-degenerate, and study them first.  For some progress in this direction,  see \cite{Akpan1, Akpan2, Akpan3}.

As follows from the above discussion, we have two types of special points on a Nijenhuis manifold $\mathsf M$.  Namely,  there are singular points at which the algebraic type of $L$ changes and also critical points for $\Phi$ at which $\operatorname{rank}\dd \Phi < n$. These sets are closely related, but do not necessarily coincide.

\begin{Problem}
\label{prob:b9}  
Let $\mathsf{Sing}\subset \mathsf M$ be the set of singular points of a Nijenhuis operator $L$ and 
$\mathsf{Crit}\subset \mathsf M$ the set of critical points of the characteristic map $\Phi = (\sigma_1,\sigma_2,\dots,\sigma_n)$ (i.e. such that $\operatorname{rank} \dd  \Phi(x)<n$).  What is the relationship between these sets?   

Let $L$ be $\gl$-regular and suppose $L(\mathsf p)$ is  a nilpotent Jordan block. In many examples of this kind (all?),  the critical set  $\mathsf{Crit}$ is a smooth hypersurface defined by the equation $\det L=0$ (as mentioned below, without loss of generality, one can assume that in a suitable coordinate system $\det L = z^n$).  Is it a general fact? 
\end{Problem}

There are some situations where the momentum mapping $\mathcal F$ and the characteristic mapping $\Phi$ essentially coincide (after some natural amendments in the definition of $\Phi$).  Recall that a Hamiltonian system is called {\it bi-Hamiltonian} ~\cite{Magri1978}
if it is Hamiltonian w.r.t. two different compatible symplectic (Poisson) structures $\omega_0$ and $\omega_1$.  Consider the so-called {\it recursion operator} $P = \omega_0^{-1}\omega_1$.  The compatibility condition implies that $P$ is a Nijenhuis operator.    Due to the fact that $\omega_0$ and $\omega_1$ are skew-symmetric,  the characteristic polynomial of $P$ is a full square:
$$
\chi_P(t) = \det (t\cdot\Id - P) = (t^n  - \sigma_1 t^{n-1} - \ldots - \sigma_n)^2. 
$$  
The characteristic map in the original sense does not make much sense in this situation, but it still makes sense to consider
the map $$\Phi : (\sigma_1,\dots,\sigma_n): \mathsf M^{2n} \to \R^n.$$   It is a well-known fact that  $\sigma_1,\dots,\sigma_n$ are Poisson commuting first integrals of every bi-Hamiltonian system associated with $\omega_0$ and $\omega_1$.   Thus, the above map is exactly the  momentum mapping.

\begin{Problem}
\label{prob:b10}
It would be interesting to study singularities of $\Phi$ in this case.  For instance, can these singularities be non-degenerate in the Eliasson sense? 
\end{Problem}

Recall that even if $\mathsf p \in \mathsf M^{2n}$ is not critical for the momentum/characteristic map $\Phi$,  this point can still be singular in the sense of Nijenhuis geometry, that is, the algebraic type of $P$ may change at $\mathsf p$.  This case was studied in \cite{nijapp1}.

There is a very natural construction due to Turiel~\cite{Turiel1994} which
allows one to construct an integrable bi-Hamiltonian system from an arbitrary Nijenhuis operator $L$.

Take an arbitrary Nijenhuis operator $L$ on a manifold $\mathsf M$ and construct a new (non-canonical)  2-form  on $T^*\mathsf M$ by setting $\omega_1 = \dd ( L^i_j p_i \dd q_j)$.
The form $\omega_1$  is always closed. It is symplectic wherever $L$ is
invertible, and  is  degenerate otherwise.
 
\begin{Problem}
\label{prob:b11}
How general is  this construction?  Which bi-Hamiltonian (symplectic) structures can be reduced to the Turiel form?  How are 
the singularities (normal forms) of the momentum mapping are related to the singular points of the original operator $L$?   
\end{Problem}

In this context, there is one quite interesting phenomenon which is not quite understood yet.  

\begin{Problem}
\label{prob:b12}
Let $L$ be a nilpotent Nijenhuis operator that `consists' of several Jordan $0$-blocks.  Then there are examples of such $L$ which cannot be reduced to a constant form.   Consider the natural lift of this operator from $\mathsf M$ to the cotangent bundle $T^*\mathsf M$.  By definition, this lift $\widehat L$ is the recursion operator for the symplectic structures $$\omega_0 = \sum_i \dd p_i \wedge \dd q_i$$ and $$\omega_1 = \dd ( L^i_j p_i \dd q_j),$$ that is, $\widehat L=\omega_0^{-1}\omega_1$.  According to Turiel~\cite{Turiel1994}, $\widehat L$ can always be reduced to a constant form in a neighborhood of almost every point. In other words, as a result of lifting, the operator becomes much simpler. Why?  
\end{Problem}

However,  in the context of applications, the most interesting problems are related to studying singularities (more generally, singular Lagrangian fibrations) for integrable Hamiltonian systems constructed from Nijenhuis operators as in Section \ref{ihs}.

\subsubsection*{Acknowledgements. } A.\,B. was  supported by the Ministry of Science and Higher Education of the Republic of Kazakhstan (grant No. AP23483476). V.\,M.\ was supported by the DFG (project 529233771), and  by the ARC  (Discovery Programme DP210100951). The work of A.K. was carried out at Lomonosov Moscow State University  with the support of the Russian Science Foundation under grant No. 24-71-10100 (\url{https://rscf.ru/en/project/24-71-10100/}). A part of the work was done at CIMAT Guanajuato during the conference {\it FDIS VIII} and workshop  on {\it Integrable Systems, Liouville Foliations, and Nijenhuis Geometry} (both August 2025) and the workshops {\it Integrable Systems and Differential Geometry}  at La Trobe University  and   {\it Nijenhuis geometry, Haantjes geometry and separation of variables} at MATRIX (both  February 2026). We thank the organisers for the support and the venues for perfect working conditions, and the participants for useful discussions.

\subsubsection*{Declaration on the use of generative AI.}
As explained in the abstract, this manuscript grew out of preparations
for workshops held in 2025 and early 2026. Most of the text was written
before publicly available generative AI tools could provide substantial
assistance with this work. During the final preparation of the manuscript,
the authors used ChatGPT (6.2 Astra) to help check mathematical arguments,
search for references, and improve the grammar. The authors carefully
reviewed the resulting suggestions and take full responsibility for the
final text. 
\printbibliography

@article {alonso,
    AUTHOR = {Mart\'{\i}nez Alonso, L.},
     TITLE = {Schr\"{o}dinger spectral problems with energy-dependent
              potentials as sources of nonlinear {H}amiltonian evolution
              equations},
   JOURNAL = {J. Math. Phys.},
  FJOURNAL = {Journal of Mathematical Physics},
    VOLUME = {21},
      YEAR = {1980},
    NUMBER = {9},
     PAGES = {2342--2349},
      ISSN = {0022-2488,1089-7658},
   MRCLASS = {58F07 (35Q20 81C05 82A45)},
  MRNUMBER = {585585},
MRREVIEWER = {J.\ E.\ Marsden},
     
       URL = {https://doi.org/10.1063/1.524690},
}

@article {Fomenko_80,
    AUTHOR = {Bolsinov, Alexey  V.},
     TITLE = {Geodesically equivalent metrics and {N}ijenhuis geometry},
   JOURNAL = {Mat. Sb.},
  FJOURNAL = {Matematicheski\u i\ Sbornik},
    VOLUME = {216},
      YEAR = {2025},
    NUMBER = {5},
     PAGES = {5--32},
      ISSN = {0368-8666,2305-2783},
   MRCLASS = {53B20 (53C15 53D17 70H06)},
  MRNUMBER = {4933098},
       DOI = {10.4213/sm10220},
       URL = {https://doi.org/10.4213/sm10220},
}

@article {fordy,
    AUTHOR = {Antonowicz, Marek and Fordy, Allan P.},
     TITLE = {Coupled {K}d{V} equations with multi-{H}amiltonian structures},
   JOURNAL = {Phys. D},
  FJOURNAL = {Physica D. Nonlinear Phenomena},
    VOLUME = {28},
      YEAR = {1987},
    NUMBER = {3},
     PAGES = {345--357},
      ISSN = {0167-2789,1872-8022},
   MRCLASS = {58F07 (35Q20)},
  MRNUMBER = {914454},
MRREVIEWER = {Malcolm\ R.\ Adams},
    
       URL = {https://doi.org/10.1016/0167-2789(87)90023-6},
}

@article {af,
    AUTHOR = {Antonowicz, Marek and Fordy, Allan P.},
     TITLE = {A family of completely integrable multi-{H}amiltonian systems},
   JOURNAL = {Phys. Lett. A},
  FJOURNAL = {Physics Letters. A},
    VOLUME = {122},
      YEAR = {1987},
    NUMBER = {2},
     PAGES = {95--99},
      ISSN = {0375-9601,1873-2429},
   MRCLASS = {58F07},
  MRNUMBER = {897492},
MRREVIEWER = {Alexander\ A.\ Pankov},
   
       URL = {https://doi.org/10.1016/0375-9601(87)90783-3},
}

@article {BM2003,
    AUTHOR = {Bolsinov, Alexey V. and Matveev, Vladimir S.},
     TITLE = {Geometrical interpretation of {B}enenti systems},
   JOURNAL = {J. Geom. Phys.},
  FJOURNAL = {Journal of Geometry and Physics},
    VOLUME = {44},
      YEAR = {2003},
    NUMBER = {4},
     PAGES = {489--506},
      ISSN = {0393-0440,1879-1662},
   MRCLASS = {37J35 (53D25 70G45 81R12)},
  MRNUMBER = {1943174},
MRREVIEWER = {Leo\ T.\ Butler},

       URL = {https://doi.org/10.1016/S0393-0440(02)00054-2},
}

@article {BM2006,
    AUTHOR = {Blaszak, Maciej and Marciniak, Krzysztof},
     TITLE = {From {S}t\"{a}ckel systems to integrable hierarchies of
              {PDE}'s: {B}enenti class of separation relations},
   JOURNAL = {J. Math. Phys.},
  FJOURNAL = {Journal of Mathematical Physics},
    VOLUME = {47},
      YEAR = {2006},
    NUMBER = {3},
     PAGES = {032904, 26},
      ISSN = {0022-2488,1089-7658},
   MRCLASS = {37J35 (35Q51 37K10)},
  MRNUMBER = {2219792},
MRREVIEWER = {Giovanni\ Rastelli},

       URL = {https://doi.org/10.1063/1.2176908},
}

@article {MB2010,
    AUTHOR = {Marciniak, Krzysztof and Blaszak, Maciej},
     TITLE = {Construction of coupled {H}arry {D}ym hierarchy and its
              solutions from {S}t\"ackel systems},
   JOURNAL = {Nonlinear Anal.},
  FJOURNAL = {Nonlinear Analysis. Theory, Methods \& Applications. An
              International Multidisciplinary Journal},
    VOLUME = {73},
      YEAR = {2010},
    NUMBER = {9},
     PAGES = {3004--3017},
      ISSN = {0362-546X,1873-5215},
   MRCLASS = {37K10 (35Q53 37J35)},
  MRNUMBER = {2678661},
MRREVIEWER = {Robert\ Brouzet},

       URL = {https://doi.org/10.1016/j.na.2010.06.067},
}

@article {BM2008,
    AUTHOR = {Blaszak, Maciej and Marciniak, Krzysztof},
     TITLE = {St\"ackel systems generating coupled {K}d{V} hierarchies and
              their finite-gap and rational solutions},
   JOURNAL = {J. Phys. A},
  FJOURNAL = {Journal of Physics. A. Mathematical and Theoretical},
    VOLUME = {41},
      YEAR = {2008},
    NUMBER = {48},
     PAGES = {485202, 17},
      ISSN = {1751-8113,1751-8121},
   MRCLASS = {37K10 (35F50 35Q53 37J35 37K20)},
  MRNUMBER = {2515879},
MRREVIEWER = {Artur\ Sergyeyev},

       URL = {https://doi.org/10.1088/1751-8113/41/48/485202},
}

@article {Matveev2020,
    AUTHOR = {Matveev, Vladimir S.},
     TITLE = {Quantum integrability for the {B}eltrami-{L}aplace operators
              of projectively equivalent metrics of arbitrary signatures},
   JOURNAL = {Chebyshevski\u i\ Sb.},
  FJOURNAL = {Chebyshevski\u i\ Sbornik},
    VOLUME = {21},
      YEAR = {2020},
    NUMBER = {2},
     PAGES = {275--289},
      ISSN = {2226-8383,2587-7119},
   MRCLASS = {58J60 (32Q15 53C22)},
  MRNUMBER = {4188508},
MRREVIEWER = {Michael\ Kunzinger},
       DOI = {10.22405/2226-8383-2020-21-2-275-289},
       URL = {https://doi.org/10.22405/2226-8383-2020-21-2-275-289},
}

@book {Bbook,
    AUTHOR = {B\l aszak, Maciej},
     TITLE = {Quantum versus classical mechanics and integrability
              problems---towards a unification of approaches and tools},
 PUBLISHER = {Springer, Cham},
      YEAR = {2019},
     PAGES = {xiii+460},
   MRCLASS = {70S05 (53D50 53D55 70-2 70H06 81S10)},
  MRNUMBER = {3932148},
MRREVIEWER = {Jeffrey\ K.\ Lawson},
}

@article {BS2023,
    AUTHOR = {Szablikowski, Blazej M. and Blaszak, Maciej and
              Marciniak, Krzysztof},
     TITLE = {Stationary coupled {K}d{V} systems and their {S}t\"ackel
              representations},
   JOURNAL = {Stud. Appl. Math.},
  FJOURNAL = {Studies in Applied Mathematics},
    VOLUME = {153},
      YEAR = {2024},
    NUMBER = {1},
     PAGES = {Paper No. e12698, 54},
      ISSN = {0022-2526,1467-9590},
   MRCLASS = {58J70 (35Q53 58A20)},
  MRNUMBER = {4773432},
}

@article {Dubrovin1975,
    AUTHOR = {Dubrovin, Boris  A.},
     TITLE = {A periodic problem for the {K}orteweg-de {V}ries equation in a
              class of finite band potentials},
   JOURNAL = {Funkcional. Anal. i Prilo\v zen.},
  FJOURNAL = {Akademija Nauk SSSR. Funkcional\cprime nyi Analiz i ego
              Prilo\v zenija},
    VOLUME = {9},
      YEAR = {1975},
    NUMBER = {3},
     PAGES = {41--51},
      ISSN = {0374-1990},
   MRCLASS = {34B30 (35Q99 58F05)},
  MRNUMBER = {486780},
MRREVIEWER = {K.\ A.\ Pyragas},
}

@article {GG81,
    AUTHOR = {Gromoll, Detlef and Grove, Karsten},
     TITLE = {On metrics on {$S\sp{2}$}\ all of whose geodesics are closed},
   JOURNAL = {Invent. Math.},
  FJOURNAL = {Inventiones Mathematicae},
    VOLUME = {65},
      YEAR = {1981},
    NUMBER = {1},
     PAGES = {175--177},
      ISSN = {0020-9910,1432-1297},
   MRCLASS = {58E10 (53C22)},
  MRNUMBER = {636885},
MRREVIEWER = {Gudlaugur\ Thorbergsson},
       DOI = {10.1007/BF01389300},
       URL = {https://doi.org/10.1007/BF01389300},
}

@book {Besse,
    AUTHOR = {Besse, Arthur L.},
     TITLE = {Manifolds all of whose geodesics are closed},
    SERIES = {Ergebnisse der Mathematik und ihrer Grenzgebiete [Results in
              Mathematics and Related Areas]},
    VOLUME = {93},
      NOTE = {With appendices by D. B. A. Epstein, J.-P. Bourguignon, L.
              B\'erard-Bergery, M. Berger and J. L. Kazdan},
 PUBLISHER = {Springer-Verlag, Berlin-New York},
      YEAR = {1978},
     PAGES = {ix+262},
      ISBN = {3-540-08158-5},
   MRCLASS = {53C20 (53C22 58G99)},
  MRNUMBER = {496885},
MRREVIEWER = {R.\ L.\ Bishop},
}

@article {nijenhuis,
    AUTHOR = {Nijenhuis, Albert},
     TITLE = {{$X_{n-1}$}-forming sets of eigenvectors},
      NOTE = {Nederl. Akad. Wetensch. Proc. Ser. A {\bf 54}},
   JOURNAL = {Indag. Math.},
  FJOURNAL = {},
    VOLUME = {13},
      YEAR = {1951},
     PAGES = {200--212},
   MRCLASS = {53.0X},
  MRNUMBER = {43540},
MRREVIEWER = {J.\ A.\ Schouten},
}

@article {nij2,
    AUTHOR = {Konyaev, Andrey Yu.},
     TITLE = {Nijenhuis geometry {II}: {L}eft-symmetric algebras and
              linearization problem for {N}ijenhuis operators},
   JOURNAL = {Differential Geom. Appl.},
  FJOURNAL = {Differential Geometry and its Applications},
    VOLUME = {74},
      YEAR = {2021},
     PAGES = {Paper No. 101706, 32},
      ISSN = {0926-2245,1872-6984},
   MRCLASS = {53B99 (17D25)},
  MRNUMBER = {4178992},
MRREVIEWER = {Liviu\ Octavian\ Popescu},
       URL = {https://doi.org/10.1016/j.difgeo.2020.101706}
}

@article {nij1,
    AUTHOR = {Bolsinov, Alexey V. and Konyaev, Andrey Yu. and Matveev,
              Vladimir S.},
     TITLE = {Nijenhuis geometry},
   JOURNAL = {Adv. Math.},
  FJOURNAL = {Advances in Mathematics},
    VOLUME = {394},
      YEAR = {2022},
     PAGES = {Paper No. 108001, 52},
      ISSN = {0001-8708,1090-2082},
   MRCLASS = {53D17},
  MRNUMBER = {4355721},
MRREVIEWER = {Mohamed\ Boucetta},
  
URL = {https://doi.org/10.1016/j.aim.2021.108001}
}

@article {nijapp1,
    AUTHOR = {Bolsinov, Alexey V. and Konyaev, Andrey Yu. and Matveev,
              Vladimir S.},
     TITLE = {Applications of {N}ijenhuis geometry: non-degenerate singular
              points of {P}oisson-{N}ijenhuis structures},
   JOURNAL = {Eur. J. Math.},
  FJOURNAL = {European Journal of Mathematics},
    VOLUME = {8},
      YEAR = {2022},
    NUMBER = {4},
     PAGES = {1355--1376},
      ISSN = {2199-675X,2199-6768},
   MRCLASS = {53D20 (37J06 53D22 70H06)},
  MRNUMBER = {4510609},
MRREVIEWER = {Liviu\ Octavian\ Popescu},
 
       URL = {https://doi.org/10.1007/s40879-020-00429-6},
}

@article {nijapp2,
    AUTHOR = {Bolsinov, Alexey V. and Konyaev, Andrey Yu. and Matveev,
              Vladimir S.},
     TITLE = {Applications of {N}ijenhuis geometry {II}: maximal pencils of
              multi-{H}amiltonian structures of hydrodynamic type},
   JOURNAL = {Nonlinearity},
  FJOURNAL = {Nonlinearity},
    VOLUME = {34},
      YEAR = {2021},
    NUMBER = {8},
     PAGES = {5136--5162},
      ISSN = {0951-7715,1361-6544},
   MRCLASS = {37K06 (37Kxx 53A20 53B 53B20 53D17)},
  MRNUMBER = {4281470},
MRREVIEWER = {I.\ S.\ Krasil\cprime shchik},

       URL = {https://doi.org/10.1088/1361-6544/abed39},
}

@article {nijapp3,
    AUTHOR = {Bolsinov, Alexey V. and Konyaev, Andrey Yu. and Matveev,
              Vladimir S.},
     TITLE = {Applications of {N}ijenhuis geometry {III}: {F}robenius
              pencils and compatible non-homogeneous {P}oisson structures},
   JOURNAL = {J. Geom. Anal.},
  FJOURNAL = {Journal of Geometric Analysis},
    VOLUME = {33},
      YEAR = {2023},
    NUMBER = {6},
     PAGES = {Paper No. 193, 52},
      ISSN = {1050-6926,1559-002X},
   MRCLASS = {53D17 (37K06 37K10 37K25 37K50 53A20 53B10)},
  MRNUMBER = {4576387},
       URL = {https://doi.org/10.1007/s12220-023-01237-6},
}

@article {CI10,
    AUTHOR = {Holm, Darryl D. and Ivanov, Rossen I.},
     TITLE = {Smooth and peaked solitons of the {CH} equation},
   JOURNAL = {J. Phys. A},
  FJOURNAL = {Journal of Physics. A. Mathematical and Theoretical},
    VOLUME = {43},
      YEAR = {2010},
    NUMBER = {43},
     PAGES = {434003, 18},
      ISSN = {1751-8113,1751-8121},
   MRCLASS = {35Q53 (35C08 37K20 37K40)},
  MRNUMBER = {2727777},
       DOI = {10.1088/1751-8113/43/43/434003},
       URL = {https://doi.org/10.1088/1751-8113/43/43/434003},
}

@article {CK,
    AUTHOR = {Cie\'slak, Tomasz and Kry\'nski, Wojciech},
     TITLE = {Geometric aspects of two- and threepeakons},
   JOURNAL = {Nonlinearity},
  FJOURNAL = {Nonlinearity},
    VOLUME = {34},
      YEAR = {2021},
    NUMBER = {9},
     PAGES = {6685--6704},
      ISSN = {0951-7715,1361-6544},
   MRCLASS = {70H06 (37J39 53C22)},
  MRNUMBER = {4304494},
MRREVIEWER = {Igor\ Leite\ Freire},
       DOI = {10.1088/1361-6544/ac149e},
       URL = {https://doi.org/10.1088/1361-6544/ac149e},
}

@article {HI10,
    AUTHOR = {Holm, Darryl  D. and Ivanov, Rossen I.},
     TITLE = {Multi-component generalizations of the {CH} equation:
              geometrical aspects, peakons and numerical examples},
   JOURNAL = {J. Phys. A},
  FJOURNAL = {Journal of Physics. A. Mathematical and Theoretical},
    VOLUME = {43},
      YEAR = {2010},
    NUMBER = {49},
     PAGES = {492001, 20},
      ISSN = {1751-8113,1751-8121},
   MRCLASS = {37K10 (35Q53 37K40 76B15)},
  MRNUMBER = {2740351},
MRREVIEWER = {Dmitry\ E.\ Pelinovsky},
       DOI = {10.1088/1751-8113/43/49/492001},
       URL = {https://doi.org/10.1088/1751-8113/43/49/492001},
}

@article {MN25,
    AUTHOR = {Matveev, Vladimir S. and Nikolayevsky, Yuri},
     TITLE = {Killing tensors on reducible spaces},
   JOURNAL = {Manuscripta Math.},
  FJOURNAL = {Manuscripta Mathematica},
    VOLUME = {176},
      YEAR = {2025},
    NUMBER = {1},
     PAGES = {Paper No. 3, 9},
      ISSN = {0025-2611,1432-1785},
   MRCLASS = {53C21 (37J35 53C29 70H06)},
  MRNUMBER = {4837949},
MRREVIEWER = {Masami\ Sekizawa},
       DOI = {10.1007/s00229-024-01599-3},
       URL = {https://doi.org/10.1007/s00229-024-01599-3},
}

@article{DullinRichterVeselovWaalkens2001,
  author  = {Dullin, Holger R. and Richter, Peter H. and
             Veselov, Alexander P. and Waalkens, Holger},
  title   = {Actions of the {Neumann} systems via {Picard--Fuchs} equations},
  journal = {Physica D: Nonlinear Phenomena},
  volume  = {155},
  year    = {2001},
  pages   = {159--183}
}

@article {CH93,
    AUTHOR = {Camassa, Roberto and Holm, Darryl D.},
     TITLE = {An integrable shallow water equation with peaked solitons},
   JOURNAL = {Phys. Rev. Lett.},
  FJOURNAL = {Physical Review Letters},
    VOLUME = {71},
      YEAR = {1993},
    NUMBER = {11},
     PAGES = {1661--1664},
      ISSN = {0031-9007,1079-7114},
   MRCLASS = {35Q51 (58F07 76B15 76B25)},
  MRNUMBER = {1234453},
       DOI = {10.1103/PhysRevLett.71.1661},
       URL = {https://doi.org/10.1103/PhysRevLett.71.1661},
}

@article{KM25,
  author  = {Konyaev, Andrey Yu. and Matveev, Vladimir S.},
  title   = {Lax pairs for {BKM} hierarchy},
  journal = {Stud. Appl. Math.},
  volume  = {157},
  number  = {3},
  year    = {2026},
  pages   = {e70314},
  doi     = {10.1111/sapm.70314}
}

@article{Degtyareva,
  title = {Classification of Three-Dimensional Linear Nijenhuis Operators with
Functionally Independent Invariants},
  author = {Degtyareva, Sofia D.},
  journal = {Moscow Univ. Math. Bull.},
  volume = {79},
  number = {},
  pages = {192-197},
  numpages = {0},
  year = {2024},
  %month = {},
  publisher = {},
  url = {https://doi.org/10.3103/S0027132224700256}
}

@misc{Abdurasulov,
  author        = {Abdurasulov, Kobiljon and
                   Daukeyeva, Nigora and
                   Saydaliyev, Azamat},
  title         = {The geometric classification of
                   left-symmetric algebras and superalgebras},
  year          = {2026},
  eprint        = {2607.28254},
  archivePrefix = {arXiv},
  primaryClass  = {math.RA},
  url           = {https://arxiv.org/abs/2607.28254}
}

@article{Degtyareva2,
  author  = {Astashov, E. A. and Degtiareva, S. D.},
  title   = {Singularities of three-dimensional linear
             {Nijenhuis} operators with functionally
             independent invariants},
  journal = {Chebyshevskii Sbornik},
  volume  = {26},
  number  = {5},
  pages   = {73--83},
  year    = {2025},
  doi     = {10.22405/2226-8383-2025-26-5-73-83},
  note    = {In Russian}
}

@article {nijapp4,
    AUTHOR = {Bolsinov, Alexey  V. and Konyaev, Andrey Yu. and Matveev,
              Vladimir  S.},
     TITLE = {Applications of {N}ijenhuis geometry {IV}: multicomponent
              {K}d{V} and {C}amassa-{H}olm equations},
   JOURNAL = {Dyn. Partial Differ. Equ.},
  FJOURNAL = {Dynamics of Partial Differential Equations},
    VOLUME = {20},
      YEAR = {2023},
    NUMBER = {1},
     PAGES = {73--98},
      ISSN = {1548-159X,2163-7873},
   MRCLASS = {35Q53},
  MRNUMBER = {4530061},

       URL = {https://doi.org/10.4310/dpde.2023.v20.n1.a4},
}

@article {VBM,
    AUTHOR = {Matveev, Vladimir B.},
     TITLE = {30 years of finite-gap integration theory},
   JOURNAL = {Philos. Trans. R. Soc. Lond. Ser. A Math. Phys. Eng. Sci.},
  FJOURNAL = {Philosophical Transactions of the Royal Society of London.
              Series A. Mathematical, Physical and Engineering Sciences},
    VOLUME = {366},
      YEAR = {2008},
    NUMBER = {1867},
     PAGES = {837--875},
      ISSN = {1364-503X,1471-2962},
   MRCLASS = {37K10 (34C08 35Q51 37K15)},
  MRNUMBER = {2377678},
   
       URL = {https://doi.org/10.1098/rsta.2007.2055},
}

@article {nij3,
    AUTHOR = {Bolsinov, Alexey V. and Konyaev, Andrey Yu. and Matveev,
              Vladimir S.},
     TITLE = {Nijenhuis geometry {III}: gl-regular {N}ijenhuis operators},
   JOURNAL = {Rev. Mat. Iberoam.},
  FJOURNAL = {Revista Matem\'{a}tica Iberoamericana},
    VOLUME = {40},
      YEAR = {2024},
    NUMBER = {1},
     PAGES = {155--188},
      ISSN = {0213-2230,2235-0616},
   MRCLASS = {53C15 (35A30 35B06 35F20 35N10)},
  MRNUMBER = {4705683},
       URL = {https://doi.org/10.4171/rmi/1416}
}

@article {m,
    AUTHOR = {Magri, Franco},
     TITLE = {Lenard chains for classical integrable systems},
   JOURNAL = {Teoret. Mat. Fiz.},
  FJOURNAL = {Teoreticheskaya i Matematicheskaya Fizika},
    VOLUME = {137},
      YEAR = {2003},
    NUMBER = {3},
     PAGES = {424--432},
      ISSN = {0564-6162,2305-3135},
   MRCLASS = {37J35 (70H06)},
  MRNUMBER = {2084151},
 
       URL = {https://doi.org/10.1023/B:TAMP.0000007919.80743.1e},
}

@article {separation,
    AUTHOR = {Bolsinov, Alexey V. and Konyaev, Andrey Yu. and Matveev,
              Vladimir S.},
     TITLE = {Orthogonal separation of variables for spaces of constant
              curvature},
   JOURNAL = {Forum Math.},
  FJOURNAL = {Forum Mathematicum},
    VOLUME = {37},
      YEAR = {2025},
    NUMBER = {1},
     PAGES = {13--41},
      ISSN = {0933-7741,1435-5337},
   MRCLASS = {53D17 (37Jxx 53B50 70G45 70H15 70H20 70S10)},
  MRNUMBER = {4846649},
       URL = {https://doi.org/10.1515/forum-2023-0300},
}

@article{gdh,
  title = {An Integrable Shallow Water Equation with Linear and Nonlinear Dispersion},
  author = {Dullin, Holger R. and Gottwald, Georg A. and Holm, Darryl D.},
  journal = {Phys. Rev. Lett.},
  volume = {87},
  issue = {19},
  pages = {194501},
  numpages = {4},
  year = {2001},
  publisher = {American Physical Society},

  url = {https://link.aps.org/doi/10.1103/PhysRevLett.87.194501}
}

@article {Krichever1976,
    AUTHOR = {Krichever, Igor Moiseevich},
     TITLE = {An algebraic-geometric construction of the {Z}aharov-\v Sabat
              equations and their periodic solutions},
   JOURNAL = {Dokl. Akad. Nauk SSSR},
  FJOURNAL = {Doklady Akademii Nauk SSSR},
    VOLUME = {227},
      YEAR = {1976},
    NUMBER = {2},
     PAGES = {291--294},
      ISSN = {0002-3264},
   MRCLASS = {58F05 (34B25 35Q99)},
  MRNUMBER = {413178},
MRREVIEWER = {Shunichi\ Tanaka},
}

@article {Krichever1977,
    AUTHOR = {Krichever, Igor Moiseevich},
     TITLE = {Methods of algebraic geometry in the theory of nonlinear
              equations},
   JOURNAL = {Uspehi Mat. Nauk},
  FJOURNAL = {Akademija Nauk SSSR i Moskovskoe Matemati\v ceskoe Ob\v s\v
              cestvo. Uspehi Matemati\v ceskih Nauk},
    VOLUME = {32},
      YEAR = {1977},
    NUMBER = {6(198)},
     PAGES = {183--208, 287},
      ISSN = {0042-1316},
   MRCLASS = {58F05 (14H99 35Q99)},
  MRNUMBER = {516323},
MRREVIEWER = {Mark\ Pankov},
}

@article {BMMT,
    AUTHOR = {Bolsinov, Alexey and Matveev, Vladimir S. and Miranda, Eva and
              Tabachnikov, Serge},
     TITLE = {Open problems, questions and challenges in finite-dimensional
              integrable systems},
   JOURNAL = {Philos. Trans. Roy. Soc. A},
  FJOURNAL = {Philosophical Transactions of the Royal Society A.
              Mathematical, Physical and Engineering Sciences},
    VOLUME = {376},
      YEAR = {2018},
    NUMBER = {2131},
     PAGES = {20170430, 40},
      ISSN = {1364-503X,1471-2962},
   MRCLASS = {37J99 (53D17)},
  MRNUMBER = {3868426},

       URL = {https://doi.org/10.1098/rsta.2017.0430},
}

@article {l,
    AUTHOR = {Arsie, Alessandro and Buryak, Alexandr and Lorenzoni, Paolo
              and Rossi, Paolo},
     TITLE = {Riemannian {F}-manifolds, bi-flat {F}-manifolds, and flat
              pencils of metrics},
   JOURNAL = {Int. Math. Res. Not. IMRN},
  FJOURNAL = {International Mathematics Research Notices. IMRN},
      YEAR = {2022},
       VOLUME = {2022},
    NUMBER = {21},
     PAGES = {16730--16778},
      ISSN = {1073-7928,1687-0247},
   MRCLASS = {53D45},
  MRNUMBER = {4504906},
MRREVIEWER = {Robert\ Brouzet},
       URL = {https://doi.org/10.1093/imrn/rnab203},
}

@incollection {Deift2019,
    AUTHOR = {Deift, Percy A.},
     TITLE = {Three lectures on "{F}ifty years of {K}d{V}: an integrable
              system"},
 BOOKTITLE = {Nonlinear dispersive partial differential equations and
              inverse scattering},
    SERIES = {Fields Inst. Commun.},
    VOLUME = {83},
     PAGES = {3--38},
 PUBLISHER = {Springer, New York},
      YEAR = {2019},
   
   MRCLASS = {37J99 (35Q53 35Q55)},
  MRNUMBER = {3931833},
MRREVIEWER = {Luiz\ Gustavo\ Farah},
}

@article {MT01,
    AUTHOR = {Matveev, Vladimir S. and Topalov, Peter J.},
     TITLE = {Quantum integrability of {B}eltrami-{L}aplace operator as
              geodesic equivalence},
   JOURNAL = {Math. Z.},
  FJOURNAL = {Mathematische Zeitschrift},
    VOLUME = {238},
      YEAR = {2001},
    NUMBER = {4},
     PAGES = {833--866},
      ISSN = {0025-5874,1432-1823},
   MRCLASS = {58J60 (37J35 53C22)},
  MRNUMBER = {1872577},
MRREVIEWER = {Iskander\ A.\ Taimanov},
   
       URL = {https://doi.org/10.1007/s002090100280},
}

@article{BKM-JordanChevalley,
  author  = {Bolsinov, Alexey V. and Konyaev, Andrey Yu. and
             Matveev, Vladimir S.},
  title   = {On the {Jordan--Chevalley} decomposition problem
             for operator fields in small dimensions and
             {Tempesta--Tondo} conjecture},
  journal = {Journal of Geometry and Physics},
  volume  = {218},
  pages   = {105656},
  year    = {2025},
  doi     = {10.1016/j.geomphys.2025.105656}
}

@misc{DMPS2026,
  author        = {Dullin, Holger R. and Matveev, Vladimir S.
                   and Palshin, Gleb P. and Scapucci, Serena},
  title         = {Rotationally Symmetric {Zoll} Metrics
                   with a Cubic Integral},
  year          = {2026},
  eprint        = {2609.31481},
  archivePrefix = {arXiv},
  primaryClass  = {math.DG},
  doi           = {10.48550/arXiv.2609.31481},
  note          = {Submitted to Journal of Geometry and Physics}
}

@article {HrMa,
    AUTHOR = {Hryniv, Rostyslav O. and Manko, Stepan S.},
     TITLE = {Inverse scattering on the half-line for energy-dependent
              {S}chr\"{o}dinger equations},
   JOURNAL = {Inverse Problems},
  FJOURNAL = {Inverse Problems. An International Journal on the Theory and
              Practice of Inverse Problems, Inverse Methods and Computerized
              Inversion of Data},
    VOLUME = {36},
      YEAR = {2020},
    NUMBER = {9},
     PAGES = {095002, 15},
      ISSN = {0266-5611,1361-6420},
   MRCLASS = {34L25 (34L40 65L09)},
  MRNUMBER = {4149872},
MRREVIEWER = {Namig\ J.\ Guliyev},

       URL = {https://doi.org/10.1088/1361-6420/aba416},
}

@article{Staeckel1893,
  author  = {St\"ackel, Paul},
  title   = {Ueber die {B}ewegung eines {P}unktes in einer $n$-fachen {M}annigfaltigkeit},
  journal = {Math. Ann.},
  volume  = {42},
  year    = {1893},
  pages   = {537--563},
  doi     = {10.1007/BF01447379}
}

@article {IMM,
    AUTHOR = {Ibort, A. and Magri, F. and Marmo, G.},
     TITLE = {Bihamiltonian structures and {S}t\"{a}ckel separability},
   JOURNAL = {J. Geom. Phys.},
  FJOURNAL = {Journal of Geometry and Physics},
    VOLUME = {33},
      YEAR = {2000},
    NUMBER = {3-4},
     PAGES = {210--228},
      ISSN = {0393-0440,1879-1662},
   MRCLASS = {37J35 (37J05 70H06)},
  MRNUMBER = {1747040},
MRREVIEWER = {Andrey\ V.\ Tsiganov},
    
       URL = {https://doi.org/10.1016/S0393-0440(99)00051-0},
}

@article{ItsMatveev1975,
  author  = {Its, A. R. and Matveev, V. B.},
  title   = {Schr\"odinger operators with finite-gap spectrum
             and {$N$}-soliton solutions of the
             {Korteweg--de Vries} equation},
  journal = {Theoret. and Math. Phys.},
  volume  = {23},
  number  = {1},
  year    = {1975},
  pages   = {343--355},
  doi     = {10.1007/BF01038218},
  note    = {Russian original: Teoret. Mat. Fiz.
             23 (1975), no. 1, 51--68}
}

@incollection {moser80,
    AUTHOR = {Moser, Jürgen},
     TITLE = {Various aspects of integrable {H}amiltonian systems},
 BOOKTITLE = {Dynamical systems ({C}.{I}.{M}.{E}. {S}ummer {S}chool,
              {B}ressanone, 1978)},
    SERIES = {Progr. Math.},
    VOLUME = {8},
     PAGES = {233--289},
 PUBLISHER = {Birkh\"{a}user, Boston, MA},
      YEAR = {1980},
      ISBN = {3-7643-3024-4},
   MRCLASS = {58F07 (34C35 34C40 70Hxx)},
  MRNUMBER = {589592},
MRREVIEWER = {Helmut\ R\"{u}ssmann},
}

@article {moser81,
    AUTHOR = {Moser, Jürgen},
     TITLE = {Various aspects of integrable {H}amiltonian systems},
      NOTE = {Translated from the English by A. P. Veselov},
   JOURNAL = {Uspekhi Mat. Nauk},
  FJOURNAL = {Akademiya Nauk SSSR i Moskovskoe Matematicheskoe Obshchestvo.
              Uspekhi Matematicheskikh Nauk},
    VOLUME = {36},
      YEAR = {1981},
    NUMBER = {5(221)},
     PAGES = {109--151, 248},
      ISSN = {0042-1316},
   MRCLASS = {58F05 (34C35 34C40 70Hxx)},
  MRNUMBER = {637435},
}

@article {nij4,
    AUTHOR = {Bolsinov, Alexey V. and Konyaev, Andrey Yu. and Matveev,
              Vladimir S.},
     TITLE = {Nijenhuis geometry {IV}: conservation laws, symmetries and
              integration of certain non-diagonalisable systems of
              hydrodynamic type in quadratures},
   JOURNAL = {Nonlinearity},
  FJOURNAL = {Nonlinearity},
    VOLUME = {37},
      YEAR = {2024},
    NUMBER = {10},
     PAGES = {Paper No. 105003, 27},
      ISSN = {0951-7715,1361-6544},
   MRCLASS = {53B30 (37K06 47N20 53A20 53B10 58J70)},
  MRNUMBER = {4792045},
}

@misc{BKMHaantjes2026,
  author        = {Bolsinov, Alexey V. and
                   Konyaev, Andrey Yu. and
                   Matveev, Vladimir S.},
  title         = {Local description of
                   $\mathfrak{gl}$-regular
                   {Haantjes} operators},
  year          = {2026},
  eprint        = {2606.04985},
  archivePrefix = {arXiv},
  primaryClass  = {math.DG},
  doi           = {10.48550/arXiv.2606.04985}
}

@article{MarchesielloSnobl2018,
  author  = {Marchesiello, Antonella and
             {\v S}nobl, Libor},
  title   = {An Infinite Family of Maximally
             Superintegrable Systems in a Magnetic
             Field with Higher Order Integrals},
  journal = {SIGMA},
  volume  = {14},
  pages   = {092},
  year    = {2018},
  doi     = {10.3842/SIGMA.2018.092}
}

@article{MarchesielloSnobl2020,
  author  = {Marchesiello, Antonella and
             {\v S}nobl, Libor},
  title   = {Classical Superintegrable Systems
             in a Magnetic Field that Separate
             in Cartesian Coordinates},
  journal = {SIGMA},
  volume  = {16},
  pages   = {015},
  year    = {2020},
  doi     = {10.3842/SIGMA.2020.015}
}

@article{Valent2017,
  author  = {Valent, Galliano},
  title   = {Superintegrable Models on Riemannian
             Surfaces of Revolution with Integrals
             of any Integer Degree (I)},
  journal = {Regular and Chaotic Dynamics},
  volume  = {22},
  number  = {4},
  pages   = {319--352},
  year    = {2017},
  doi     = {10.1134/S1560354717040013}
}

@misc{KTF_BKMDuality2026,
  author = {Bolsinov, Alexey V. and Konyaev, Andrey Yu. and Matveev, Vladimir S.},
  title = {Duality of Operator Frobenius Algebras and Solution of the Eisenhart--St{"a}ckel Problem in the Non-Diagonal Case},
  year = {2026}, eprint = {2604.03195}, archivePrefix = {arXiv}, primaryClass = {math.DG}
}

@article{KTF_Vollmer2026,
  author = {Vollmer, Andreas},
  title = {Manifolds with a Commutative and Associative Product Structure that Encodes Superintegrable Hamiltonian Systems},
  journal = {Annales Henri Poincar{\'e}}, volume = {27}, pages = {2589--2604},
  year = {2026}, doi = {10.1007/s00023-025-01608-5}
}

@misc{KTF_JiangEtAlI2025,
  author = {Jiang, Lingrui and Liu, Si-Qi and Tian, Yingchao and Zhang, Youjin},
  title = {Generalized Frobenius Manifold Structures on the Orbit Spaces of Affine Weyl Groups {I}},
  year = {2025}, eprint = {2506.13656}, archivePrefix = {arXiv}, primaryClass = {math.DG}
}

@article{KTF_JiangEtAlII2026,
  author = {Jiang, Lingrui and Liu, Si-Qi and Tian, Yingchao and Zhang, Youjin},
  title = {Generalized Frobenius Manifold Structures on the Orbit Spaces of Affine Weyl Groups {II}},
  journal = {Journal of Geometry and Physics}, year = {2026},
  pages = {105877}, doi = {10.1016/j.geomphys.2026.105877},
  eprint = {2602.21532}, archivePrefix = {arXiv}
}

@article{KTF_BialyMironov2011,
  author = {Bialy, Misha and Mironov, Andrey E.},
  title = {Rich Quasi-linear System for Integrable Geodesic Flows on {2}-Torus},
  journal = {Discrete and Continuous Dynamical Systems}, volume = {29}, number = {1},
  pages = {81--90}, year = {2011}, doi = {10.3934/dcds.2011.29.81}
}

@article{KTF_FerapontovPavlovVitolo2016,
  author = {Ferapontov, E. V. and Pavlov, M. V. and Vitolo, R. F.},
  title = {Towards the Classification of Homogeneous Third-Order Hamiltonian Operators},
  journal = {International Mathematics Research Notices}, volume = {2016}, number = {22},
  pages = {6829--6855}, year = {2016}, doi = {10.1093/imrn/rnv369}
}

@article{Matveev1996,
  author  = {Matveev, Vladimir S.},
  title   = {Integrable {Hamiltonian} systems with two degrees
             of freedom. The topological structure of saturated
             neighbourhoods of points of focus-focus and
             saddle-saddle type},
  journal = {Sb. Math.},
  volume  = {187},
  number  = {4},
  year    = {1996},
  pages   = {495--524},
  doi     = {10.1070/SM1996v187n04ABEH000122},
  note    = {Russian original: Mat. Sb. 187 (1996),
             no. 4, 29--58}
}

@article{BolsinovMatveev1996,
  author  = {Bolsinov, Alexey  V. and Matveev, Vladimir  S.},
  title   = {Singularities of momentum maps of integrable
             {Hamiltonian} systems with two degrees of freedom},
  journal = {Zap. Nauchn. Sem. POMI},
  volume  = {235},
  year    = {1996},
  pages   = {54--86},
  note    = {English translation: J. Math. Sci. (New York)
             94 (1999), no. 4, 1477--1500}
}

@article{JaulentMiodek1976,
  author  = {Jaulent, M. and Miodek, I.},
  title   = {Nonlinear evolution equations associated with
             ``energy-dependent {Schr\"odinger} potentials''},
  journal = {Lett. Math. Phys.},
  volume  = {1},
  number  = {3},
  year    = {1976},
  pages   = {243--250},
  doi     = {10.1007/BF00417611}
}

@article{Zung1997,
  author  = {Zung, Nguyen Tien},
  title   = {A note on focus-focus singularities},
  journal = {Differential Geom. Appl.},
  volume  = {7},
  number  = {2},
  year    = {1997},
  pages   = {123--130},
  doi     = {10.1016/S0926-2245(96)00042-3}
}

@article{VuNgoc2003,
  author  = {V{\~u} Ng{\d{o}}c, San},
  title   = {On semi-global invariants for focus-focus singularities},
  journal = {Topology},
  volume  = {42},
  number  = {2},
  year    = {2003},
  pages   = {365--380},
  doi     = {10.1016/S0040-9383(01)00026-X}
}

@misc{KTF_Eastwood2023,
  author = {Eastwood, Michael}, title = {Killing Tensors on Complex Projective Space},
  year = {2023}, eprint = {2309.00589}, archivePrefix = {arXiv}, primaryClass = {math.DG}
}

@article{KTF_MatveevNikolayevsky2024,
  author = {Matveev, Vladimir S. and Nikolayevsky, Yuri},
  title = {Quadratic Killing Tensors on Symmetric Spaces which are not Generated by Killing Vector Fields},
  journal = {Comptes Rendus Math{\'e}matique}, volume = {362}, pages = {1043--1049},
  year = {2024}, doi = {10.5802/crmath.624}
}

@misc{KTF_MatveevNikolayevsky2026,
  author = {Matveev, Vladimir S. and Nikolayevsky, Yuri},
  title = {On Killing Tensors on Riemannian Symmetric Spaces},
  year = {2026}, eprint = {2604.27950}, archivePrefix = {arXiv}, primaryClass = {math.DG}
}

@misc{KTF_EastwoodLeistner2026,
  author = {Eastwood, Michael and Leistner, Thomas},
  title = {Prolongation and Killing Two-Tensors},
  year = {2026}, eprint = {2604.17855}, archivePrefix = {arXiv}, primaryClass = {math.DG}
}

@misc{KTF_DearricottNikolayevsky2026,
  author = {Dearricott, Owen and Nikolayevsky, Yuri},
  title = {Killing Tensors on Projective Spaces},
  year = {2026}, eprint = {2607.05819}, archivePrefix = {arXiv}, primaryClass = {math.DG}
}

@misc{KTF_MatveevNguyenNikolayevsky2026,
  author = {Matveev, Vladimir S. and Nguyen, An Ky and Nikolayevsky, Yuri},
  title = {Quadratic Killing Tensors on Classical Lie Groups are Decomposable},
  year = {2026}, eprint = {2607.04208}, archivePrefix = {arXiv}, primaryClass = {math.DG}
}

@article{KTF_BertoCuevasDelBarcoMoroianu2026,
  author = {Berto-Cuevas, Renan and del Barco, Viviana and Moroianu, Andrei},
  title = {Symmetric Killing Tensors on Almost Abelian Lie Groups},
  journal = {Differential Geometry and its Applications}, volume = {104},
  pages = {102416}, year = {2026}, doi = {10.1016/j.difgeo.2026.102416}
}

@article {KTF_NguyenNikolayevsky2026,
    AUTHOR = {Nguyen, An Ky and Nikolayevsky, Yuri},
     TITLE = {Quadratic {K}illing tensors on some symmetric spaces of higher
              rank},
   JOURNAL = {J. Geom. Phys.},
  FJOURNAL = {Journal of Geometry and Physics},
    VOLUME = {231},
      YEAR = {2027},
     PAGES = {Paper No. 106011},
      ISSN = {0393-0440,1879-1662},
   MRCLASS = {53C35 (53B20 53D25)},
  MRNUMBER = {5140018},
       DOI = {10.1016/j.geomphys.2026.106011},
       URL = {https://doi.org/10.1016/j.geomphys.2026.106011},
}

@misc{AitAissaElBourkadiMansouri2026,
  author        = {A{\"i}t Aissa, T. and
                   El Bourkadi, S. and
                   Mansouri, M. W.},
  title         = {Three-Dimensional Real Affine Lie Groups},
  year          = {2026},
  eprint        = {2606.29317},
  archivePrefix = {arXiv},
  primaryClass  = {math.SG},
  url           = {https://arxiv.org/abs/2606.29317}
}

@article{KTF_LorenzoniPerlettiVanGemst2025,
  author  = {Lorenzoni, Paolo and Perletti, Sara and van Gemst, Karoline},
  title   = {Integrable Hierarchies and {$F$}-Manifolds with Compatible Connection},
  journal = {Communications in Mathematical Physics},
  volume  = {406},
  pages   = {91},
  year    = {2025},
  doi     = {10.1007/s00220-025-05262-0}
}

@article{KTF_LorenzoniPerlettiVanGemstSemiHam2026,
  author  = {Lorenzoni, Paolo and Perletti, Sara and van Gemst, Karoline},
  title   = {Semi-{Hamiltonian} Properties of a Class of Non-diagonalisable Systems of Hydrodynamic Type},
  journal = {Communications in Mathematical Physics},
  volume  = {407},
  pages   = {46},
  year    = {2026},
  doi     = {10.1007/s00220-025-05544-7}
}

@article{KTF_LorenzoniPerlettiVanGemstHodograph2026,
  author  = {Lorenzoni, Paolo and Perletti, Sara and van Gemst, Karoline},
  title   = {The Generalised Hodograph Method for Non-diagonalisable Integrable Systems of Hydrodynamic Type},
  journal = {Nonlinearity},
  volume  = {39},
  number  = {5},
  pages   = {055005},
  year    = {2026},
  doi     = {10.1088/1361-6544/ae65fd}
}

@misc{KTF_ArsieLorenzoniCyclic2026,
  author        = {Arsie, Alessandro and Lorenzoni, Paolo},
  title         = {Cyclic {$F$}-Manifolds, Distinguished Connections and Integrability},
  year          = {2026},
  eprint        = {2605.25277},
  archivePrefix = {arXiv},
  primaryClass  = {math-ph},
  url           = {https://arxiv.org/abs/2605.25277}
}

@misc{KTF_ArsieLorenzoniHaantjes2026,
  author        = {Arsie, Alessandro and Lorenzoni, Paolo},
  title         = {Haantjes Torsion and Integrability: A Proof of {Bolsinov--Konyaev--Matveev}'s Conjecture},
  year          = {2026},
  eprint        = {2607.29373},
  archivePrefix = {arXiv},
  primaryClass  = {math-ph},
  url           = {https://arxiv.org/abs/2607.29373}
}

@article{Qiao1997GeneralizedNeumann,
  author  = {Qiao, Zhijun},
  title   = {Completely Integrable Generalized {C. Neumann}
             Systems on Several Symplectic Submanifolds},
  journal = {Journal of Nonlinear Mathematical Physics},
  volume  = {4},
  number  = {1-2},
  pages   = {221--230},
  year    = {1997},
  doi     = {10.2991/jnmp.1997.4.1-2.36}
}

@article{Qiao2003CamassaHolm,
  author  = {Qiao, Zhijun},
  title   = {The {Camassa--Holm} Hierarchy, {$N$}-Dimensional
             Integrable Systems, and Algebro-Geometric Solution
             on a Symplectic Submanifold},
  journal = {Communications in Mathematical Physics},
  volume  = {239},
  pages   = {309--341},
  year    = {2003},
  doi     = {10.1007/s00220-003-0880-y}
}

@article{KTF_BKMRiemannInvariants2026,
  author  = {Bolsinov, Alexey V. and Konyaev, Andrey Yu. and Matveev, Vladimir S.},
  title   = {Existence of Riemannian Invariants for Integrable Systems of Hydrodynamic Type},
  journal = {Physica D: Nonlinear Phenomena},
  volume  = {495},
  pages   = {135304},
  year    = {2026},
  doi     = {10.1016/j.physd.2026.135304}
}

@misc{AkpanManenti2026,
  author        = {Akpan, Dinmukhammed and Manenti, Matteo},
  title         = {On Three-Dimensional
                   {$\operatorname{gl}$}-Regular
                   Nijenhuis Operators},
  year          = {2026},
  eprint        = {2607.04913},
  archivePrefix = {arXiv},
  primaryClass  = {math.DG},
  url           = {https://arxiv.org/abs/2607.04913}
}

@article{Valent2021,
  author  = {Valent, Galliano},
  title   = {Superintegrable Geodesic Flows versus
             {Zoll} Metrics},
  journal = {Journal of Geometry and Physics},
  volume  = {159},
  pages   = {103873},
  year    = {2021},
  doi     = {10.1016/j.geomphys.2020.103873}
}

@article{QuaschnerSteneker2026,
  author        = {Quaschner, Manuel and Steneker, Wijnand},
  title         = {Jet-density of finite-gap solutions
                   for classes of {BKM} systems},
  journal       = {Nonlinearity},
  volume        = {39},
  number        = {8},
  pages         = {085027},
  year          = {2026},
  eprint        = {2604.16233},
  archivePrefix = {arXiv},
  url           = {https://arxiv.org/abs/2604.16233}
}

@article {BKM_magnetic,
    AUTHOR = {Bolsinov, Alexey V. and Konyaev, Andrey Yu. and Matveev,
              Vladimir S.},
     TITLE = {Integrability of the magnetic geodesic flow on the sphere with
              a constant 2-form},
   JOURNAL = {Int. J. Geom. Methods Mod. Phys.},
  FJOURNAL = {International Journal of Geometric Methods in Modern Physics},
    VOLUME = {23},
      YEAR = {2026},
    NUMBER = {15},
     PAGES = {Paper No. 2650041, 13},
      ISSN = {0219-8878,1793-6977},
   MRCLASS = {37J35 (70H06)},
  MRNUMBER = {5105102},
       DOI = {10.1142/S0219887826500416},
       URL = {https://doi.org/10.1142/S0219887826500416},
}

@article{KubuMarchesielloSnobl2023,
  author  = {Kub{\r u}, Ond{\v r}ej and
             Marchesiello, Antonella and
             {\v S}nobl, Libor},
  title   = {New classes of quadratically
             integrable systems in magnetic fields:
             The generalized cylindrical and
             spherical cases},
  journal = {Annals of Physics},
  volume  = {451},
  pages   = {169264},
  year    = {2023},
  doi     = {10.1016/j.aop.2023.169264}
}

@article{MarchesielloReyesSnobl2024,
  author  = {Marchesiello, Antonella and
             Reyes, Daniel and {\v S}nobl, Libor},
  title   = {Superintegrable families of magnetic
             monopoles with non-radial potential
             in curved background},
  journal = {Journal of Geometry and Physics},
  volume  = {203},
  pages   = {105261},
  year    = {2024},
  doi     = {10.1016/j.geomphys.2024.105261}
}

@article{AgapovDemina2027,
  author  = {Agapov, Sergei and Demina, Maria},
  title   = {On superintegrable magnetic geodesic
             flows on 2-surfaces},
  journal = {Journal of Geometry and Physics},
  volume  = {232},
  pages   = {106036},
  year    = {2027},
  doi     = {10.1016/j.geomphys.2026.106036}
}

@article {nab,
    AUTHOR = {Nabiev, A. A.},
     TITLE = {Inverse scattering problem for the {S}chr\"{o}dinger-type
              equation with a polynomial energy-dependent potential},
   JOURNAL = {Inverse Problems},
  FJOURNAL = {Inverse Problems. An International Journal on the Theory and
              Practice of Inverse Problems, Inverse Methods and Computerized
              Inversion of Data},
    VOLUME = {22},
      YEAR = {2006},
    NUMBER = {6},
     PAGES = {2055--2068},
      ISSN = {0266-5611,1361-6420},
   MRCLASS = {34A55 (34L25 34L40)},
  MRNUMBER = {2277529},
MRREVIEWER = {Mehmet\ Can},
    
       URL = {https://doi.org/10.1088/0266-5611/22/6/009},
}

@article {MT97,
    AUTHOR = {Matveev, Vladimir  S. and Topalov, Peter  I.},
     TITLE = {Trajectory equivalence and corresponding integrals},
   JOURNAL = {Regul. Chaotic Dyn.},
  FJOURNAL = {Regular \& Chaotic Dynamics. Regulyarnaya \& Khaoticheskaya
              Dinamika},
    VOLUME = {3},
      YEAR = {1998},
    NUMBER = {2},
     PAGES = {30--45},
      ISSN = {1560-3547,1468-4845},
   MRCLASS = {37J35 (53D25 70H06)},
  MRNUMBER = {1693470},
MRREVIEWER = {Iskander\ A.\ Taimanov},
  
       URL = {https://doi.org/10.1070/rd1998v003n02ABEH000069},
}

@article {novikov74,
    AUTHOR = {Novikov, Sergei P.},
     TITLE = {A periodic problem for the {K}orteweg-de {V}ries equation.
              {I}},
   JOURNAL = {Funkcional. Anal. i Prilo\v{z}en.},
  FJOURNAL = {Akademija Nauk SSSR. Funkcional\cprime nyi Analiz i ego
              Prilo\v{z}enija},
    VOLUME = {8},
      YEAR = {1974},
    NUMBER = {3},
     PAGES = {54--66},
      ISSN = {0374-1990},
   MRCLASS = {35Q99},
  MRNUMBER = {382878},
MRREVIEWER = {H.\ Flaschka},
}

@incollection {Kosmann,
    AUTHOR = {Kosmann-Schwarzbach, Yvette},
     TITLE = {Beyond recursion operators},
 BOOKTITLE = {Geometric methods in physics {XXXVI}},
    SERIES = {Trends Math.},
     PAGES = {167--180},
 PUBLISHER = {Birkh\"auser/Springer, Cham},
      YEAR = {2019},
   MRCLASS = {53A55 (01A60 17B70 39B12 53D45 70G45)},
  MRNUMBER = {3991134},
       DOI = {10.1007/978-3-030-01156-7\_18},
       URL = {https://doi.org/10.1007/978-3-030-01156-7_18},
}

@article{DubrovinMatveevNovikov1976,
  author  = {Dubrovin, Boris  A. and Matveev, Vladimir  B. and Novikov, Sergei  P.},
  title   = {Non-linear equations of {Korteweg--de Vries} type,
             finite-zone linear operators, and {Abelian} varieties},
  journal = {Russian Math. Surveys},
  volume  = {31},
  number  = {1},
  year    = {1976},
  pages   = {59--146},
  doi     = {10.1070/RM1976v031n01ABEH001446},
  note    = {Russian original: Uspekhi Mat. Nauk
             31 (1976), no. 1(187), 55--136}
}

@article {Bogo2004,
    AUTHOR = {Bogoyavlenski\u i, Oleg I.},
     TITLE = {General algebraic identities for the {N}ijenhuis and
              {H}aantjes tensors},
   JOURNAL = {Izv. Ross. Akad. Nauk Ser. Mat.},
  FJOURNAL = {Izvestiya Rossiiskoi Akademii Nauk. Seriya Matematicheskaya},
    VOLUME = {68},
      YEAR = {2004},
    NUMBER = {6},
     PAGES = {71--84},
      ISSN = {1607-0046,2587-5906},
   MRCLASS = {53C15},
  MRNUMBER = {2108522},
MRREVIEWER = {Olga\ V.\ Radko},
       DOI = {10.1070/IM2004v068n06ABEH000511},
       URL = {https://doi.org/10.1070/IM2004v068n06ABEH000511},
}

@article{magnetic,
  author    = {Bolsinov, Alexey V. and Konyaev, Andrey Yu. and Matveev, Vladimir S.},
  title     = {Integrability of the magnetic geodesic flow on the sphere with a constant 2-form},
  journal   = {International Journal of Geometric Methods in Modern Physics},
  year      = {2026},
  pages     = {2650041},
  note      = {Article ID 2650041 (13 pages)},
  doi       = {10.1142/S0219887826500416},
  issn      = {0219-8878},
  publisher = {World Scientific Publishing Company}
}

@article {BKW1983,
    AUTHOR = {Boyer, C. P. and Kalnins, E. G. and Winternitz, P.},
     TITLE = {Completely integrable relativistic {H}amiltonian systems and
              separation of variables in {H}ermitian hyperbolic spaces},
   JOURNAL = {J. Math. Phys.},
  FJOURNAL = {Journal of Mathematical Physics},
    VOLUME = {24},
      YEAR = {1983},
    NUMBER = {8},
     PAGES = {2022--2034},
      ISSN = {0022-2488,1089-7658},
   MRCLASS = {70H20 (20G45 53B30 53C80 58F07)},
  MRNUMBER = {713533},
MRREVIEWER = {Mauro\ Francaviglia},
       DOI = {10.1063/1.525943},
       URL = {https://doi.org/10.1063/1.525943},
}

@article {BKW1985,
    AUTHOR = {Boyer, C. P. and Kalnins, E. G. and Winternitz, P.},
     TITLE = {Separation of variables for the {H}amilton-{J}acobi equation
              on complex projective spaces},
   JOURNAL = {SIAM J. Math. Anal.},
  FJOURNAL = {SIAM Journal on Mathematical Analysis},
    VOLUME = {16},
      YEAR = {1985},
    NUMBER = {1},
     PAGES = {93--109},
   MRCLASS = {58G30 (35A30 53B21 53C57)},
  MRNUMBER = {772870},
MRREVIEWER = {Yvette\ Kosmann-Schwarzbach},
       DOI = {10.1137/0516006},
       URL = {https://doi.org/10.1137/0516006},
}

@article{WOR1994,
    AUTHOR = {Winternitz, P. and del Olmo, M. A. and Rodr\'iguez, M. A.},
     TITLE = {Maximal abelian subgroups of {${\rm {SU}}(p,q)$} and integrable
              {H}amiltonian systems},
    JOURNAL = { Noncompact {L}ie groups and some of their applications ({S}an
              {A}ntonio, {TX}, 1993), NATO Adv. Sci. Inst. Ser. C: Math. Phys. Sci.},
    VOLUME = {429},
     PAGES = {181--198},
 PUBLISHER = {Kluwer Acad. Publ., Dordrecht},
      YEAR = {1994},
   MRCLASS = {58F07 (22E70)},
  MRNUMBER = {1306524},
MRREVIEWER = {Malcolm\ R.\ Adams},
}

@article{BDMN,
  author  = {Bolsinov, Alexey and Dullin, Holger R. and Matveev, Vladimir S. and Nikolayevsky, Yuri},
  title   = {On separation of variables for symmetric spaces of rank~1},
  journal = {Comptes Rendus. Math{\'e}matique},
  volume  = {364},
  pages   = {267--278},
  year    = {2026},
  doi     = {10.5802/crmath.823}
}

@article{ORW1993,
    AUTHOR = {del Olmo, M. A. and Rodriguez, M. A. and Winternitz, P.},
     TITLE = {Integrable systems based on {${\rm {SU}}(p,q)$} homogeneous
              manifolds},
   JOURNAL = {J. Math. Phys.},
  FJOURNAL = {Journal of Mathematical Physics},
    VOLUME = {34},
      YEAR = {1993},
    NUMBER = {11},
     PAGES = {5118--5139},
   MRCLASS = {58F07 (58F05 70H05 81Q99)},
  MRNUMBER = {1243130},
MRREVIEWER = {Akira\ Yoshioka},
    DOI = {10.1063/1.530346},
       URL = {https://doi.org/10.1063/1.530346},
}

@article {ORWZ1990,
    AUTHOR = {del Olmo, M. A. and Rodr\'iguez, M. A. and Winternitz, P. and
              Zassenhaus, H.},
     TITLE = {Maximal abelian subalgebras of pseudounitary {L}ie algebras},
   JOURNAL = {Linear Algebra Appl.},
  FJOURNAL = {Linear Algebra and its Applications},
    VOLUME = {135},
      YEAR = {1990},
     PAGES = {79--151},
      ISSN = {0024-3795,1873-1856},
   MRCLASS = {17B05 (17B30 81R05)},
  MRNUMBER = {1061531},
MRREVIEWER = {Joris\ Van der Jeugt},
       DOI = {10.1016/0024-3795(90)90117-U},
       URL = {https://doi.org/10.1016/0024-3795(90)90117-U},
}

@article {sewa,
    AUTHOR = {Matveev, Vladimir S. and Shevchishin, Vsevolod V.},
     TITLE = {Two-dimensional superintegrable metrics with one linear and
              one cubic integral},
   JOURNAL = {J. Geom. Phys.},
  FJOURNAL = {Journal of Geometry and Physics},
    VOLUME = {61},
      YEAR = {2011},
    NUMBER = {8},
     PAGES = {1353--1377},
      ISSN = {0393-0440,1879-1662},
   MRCLASS = {53D25 (37J15 53B20 70H06 70H33)},
  MRNUMBER = {2802477},
MRREVIEWER = {Martin\ L\'egar\'e},
       DOI = {10.1016/j.geomphys.2011.02.012},
       URL = {https://doi.org/10.1016/j.geomphys.2011.02.012},
}

@article {TT2021,
    AUTHOR = {Tempesta, Piergiulio and Tondo, Giorgio},
     TITLE = {Haantjes algebras and diagonalization},
   JOURNAL = {J. Geom. Phys.},
  FJOURNAL = {Journal of Geometry and Physics},
    VOLUME = {160},
      YEAR = {2021},
     PAGES = {Paper No. 103968, 21},
      ISSN = {0393-0440,1879-1662},
   MRCLASS = {53A45 (58C40)},
  MRNUMBER = {4170275},
       DOI = {10.1016/j.geomphys.2020.103968},
       URL = {https://doi.org/10.1016/j.geomphys.2020.103968},
}

@article {MZZ,
    AUTHOR = {Marquette, Ian and Zhang, Junze and Zhang, Yao-Zhong},
     TITLE = {Algebraic approach and exact solutions of superintegrable
              systems in 2{D} {D}arboux spaces},
   JOURNAL = {J. Phys. A},
  FJOURNAL = {Journal of Physics. A. Mathematical and Theoretical},
    VOLUME = {56},
      YEAR = {2023},
    NUMBER = {35},
     PAGES = {Paper No. 355201, 50},
      ISSN = {1751-8113,1751-8121},
   MRCLASS = {81R12 (70H06 81Q35)},
  MRNUMBER = {4634076},
MRREVIEWER = {V.\ S.\ Otchik},
}

@article {KKMW,
    AUTHOR = {Kalnins, Ernie G. and Kress, Jonathan M. and Miller, Jr., Willard and
              Winternitz, Pavel},
     TITLE = {Superintegrable systems in {D}arboux spaces},
   JOURNAL = {J. Math. Phys.},
  FJOURNAL = {Journal of Mathematical Physics},
    VOLUME = {44},
      YEAR = {2003},
    NUMBER = {12},
     PAGES = {5811--5848},
      ISSN = {0022-2488,1089-7658},
   MRCLASS = {37J35 (70G45 70H06 70H20)},
  MRNUMBER = {2023556},
MRREVIEWER = {Anatoly\ G.\ Meshkov},
       DOI = {10.1063/1.1619580},
       URL = {https://doi.org/10.1063/1.1619580},
}

@article {Matveev2025,
    AUTHOR = {Matveev, Vladimir S.},
     TITLE = {Real analyticity of 2-dimensional superintegrable metrics and
              solution of two {B}olsinov-{K}ozlov-{F}omenko conjectures},
   JOURNAL = {Regul. Chaotic Dyn.},
  FJOURNAL = {Regular and Chaotic Dynamics. International Scientific
              Journal},
    VOLUME = {30},
      YEAR = {2025},
    NUMBER = {4},
     PAGES = {677--687},
      ISSN = {1560-3547,1468-4845},
   MRCLASS = {37J35 (53D25 70H06)},
  MRNUMBER = {4944200},
       DOI = {10.1134/S1560354725040148},
       URL = {https://doi.org/10.1134/S1560354725040148},
}

@unpublished{BKM26a,
  author        = {Bolsinov, Alexey V. and Konyaev, Andrey Yu. and Matveev, Vladimir S.},
  title         = {{St\"ackel} problem for non-diagonal {Killing} tensors:
                   {Yano--Patterson} lifts, algebra of strong symmetries
                   and quadratic in momenta integrals},
  year          = {2025},
  eprint        = {2512.17609},
  archivePrefix = {arXiv},
  primaryClass  = {nlin.SI},
  doi           = {10.48550/arXiv.2512.17609},
  note          = {Accepted for publication in Journal of Nonlinear Science}
}

@article {TT2022,
    AUTHOR = {Tempesta, Piergiulio and Tondo, Giorgio},
     TITLE = {Haantjes algebras of classical integrable systems},
   JOURNAL = {Ann. Mat. Pura Appl. (4)},
  FJOURNAL = {Annali di Matematica Pura ed Applicata. Series IV},
    VOLUME = {201},
      YEAR = {2022},
    NUMBER = {1},
     PAGES = {57--90},
      ISSN = {0373-3114,1618-1891},
   MRCLASS = {37J37 (70H06)},
  MRNUMBER = {4375001},
MRREVIEWER = {Michel\ Nguiffo B. Boyom},
       DOI = {10.1007/s10231-021-01107-4},
       URL = {https://doi.org/10.1007/s10231-021-01107-4},
}

@article {Pavlov2014,
    AUTHOR = {Pavlov, Maxim V.},
     TITLE = {Integrable dispersive chains and energy dependent
              {S}chr\"odinger operator},
   JOURNAL = {J. Phys. A},
  FJOURNAL = {Journal of Physics. A. Mathematical and Theoretical},
    VOLUME = {47},
      YEAR = {2014},
    NUMBER = {29},
     PAGES = {295204, 22},
      ISSN = {1751-8113,1751-8121},
   MRCLASS = {35Q53 (35A30)},
  MRNUMBER = {3230263},

       URL = {https://doi.org/10.1088/1751-8113/47/29/295204},
}

@article {MY79,
    AUTHOR = {Matveev, V. B. and Yavor, M. I.},
     TITLE = {Solutions presque p\'eriodiques et \`a{} {$N$}-solitons de
              l'\'equation hydrodynamique non lin\'eaire de {K}aup},
   JOURNAL = {Ann. Inst. H. Poincar\'e{} Sect. A (N.S.)},
  FJOURNAL = {Annales de l'Institut Henri Poincar\'e. Section A. Physique
              Th\'eorique. Nouvelle S\'erie},
    VOLUME = {31},
      YEAR = {1979},
    NUMBER = {1},
     PAGES = {25--41},
      ISSN = {0246-0211},
   MRCLASS = {35Q20 (76B25)},
  MRNUMBER = {557050},
MRREVIEWER = {W.\ G\"uttinger},
}

@article {Nijenhuis1955,
    AUTHOR = {Nijenhuis, A.},
     TITLE = {Jacobi-type identities for bilinear differential concomitants
              of certain tensor fields.  {II}},
      NOTE = {Nederl. Akad. Wetensch. Proc. Ser. A {{\bf{5}}8}},
   JOURNAL = {Indag. Math.},
    VOLUME = {17},
      YEAR = {1955},
     PAGES = { 398--403},
   MRCLASS = {53.1X},
  MRNUMBER = {0074879},
MRREVIEWER = {K. Yano},
}

@article {Haantjes,
    AUTHOR = {Haantjes, J.},
     TITLE = {On {$X_m$}-forming sets of eigenvectors},
      NOTE = {Nederl. Akad. Wetensch. Proc. Ser. A {\bf 58}},
   JOURNAL = {Indag. Math.},
  FJOURNAL = {},
    VOLUME = {17},
      YEAR = {1955},
     PAGES = {158--162},
   MRCLASS = {53.0X},
  MRNUMBER = {70232},
MRREVIEWER = {K.\ Yano},
}

@article {veselov,
    AUTHOR = {Veselov, Alexander  P.},
     TITLE = {Finite-zone potentials and integrable systems on a sphere with
              quadratic potential},
   JOURNAL = {Funktsional. Anal. i Prilozhen.},
  FJOURNAL = {Akademiya Nauk SSSR. Funktsional\cprime ny\u{\i} Analiz i ego
              Prilozheniya},
    VOLUME = {14},
      YEAR = {1980},
    NUMBER = {1},
     PAGES = {48--50},
      ISSN = {0374-1990},
   MRCLASS = {58F07 (34B25)},
  MRNUMBER = {565097},
MRREVIEWER = {Mark\ Pankov},
}

@article {BKMreduction,
    AUTHOR = {Bolsinov, Alexey  V. and Konyaev, Adnrey  Yu. and Matveev,
              Vladimir S.},
     TITLE = {Finite-dimensional reductions and finite-gap-type solutions of
              multicomponent integrable {PDE}s},
   JOURNAL = {Stud. Appl. Math.},
  FJOURNAL = {Studies in Applied Mathematics},
    VOLUME = {155},
      YEAR = {2025},
    NUMBER = {2},
     PAGES = {Paper No. e70100, 30},
      ISSN = {0022-2526,1467-9590},
   MRCLASS = {37Kxx},
  MRNUMBER = {4949662},
       URL = {https://doi.org/10.1111/sapm.70100},
}

@article {KM2025,
    AUTHOR = {    Kress, Jonathan M.  and   Matveev, Vladimir S.},
     TITLE = {Self-adjoint quantization of Stäckel integrable systems},
       JOURNAL =   {J. Phys. A: Math. Theor.},   VOLUME= {58},   YEAR = {2025}, PAGES= {345202}, url ={  https://arxiv.org/abs/2501.18082}
}

@incollection{Scapucci,
  author    = {Scapucci, Serena},
  title     = {Classification of differentially non-degenerate left-symmetric algebras in dimension 3},
  booktitle = {2024 MATRIX Annals, Part I},
  editor    = {Wood, David R. and Etheridge, Alison M. and de Gier, Jan and Joshi, Nalini},
  series    = {MATRIX Book Series},
  volume    = {7},
  publisher = {Springer},
  year      = {2026},
  pages     = {289--315},
  url       = {https://www.matrix-inst.org.au/wp_Matrix2016/wp-content/uploads/2024/2024_Dullin/Scapucci.pdf}
}

@article {NTT2021,
    AUTHOR = {Nozaleda, Daniel Reyes and Tempesta, Piergiulio and Tondo,
              Giorgio},
     TITLE = {Classical multiseparable {H}amiltonian systems,
              superintegrability and {H}aantjes geometry},
   JOURNAL = {Commun. Nonlinear Sci. Numer. Simul.},
  FJOURNAL = {Communications in Nonlinear Science and Numerical Simulation},
    VOLUME = {104},
      YEAR = {2022},
     PAGES = {Paper No. 106021, 25},
      ISSN = {1007-5704,1878-7274},
   MRCLASS = {37J06 (34A26 37J35)},
  MRNUMBER = {4318302},
MRREVIEWER = {Patricia\ Yanguas},
       DOI = {10.1016/j.cnsns.2021.106021},
       URL = {https://doi.org/10.1016/j.cnsns.2021.106021},
}

@article {Magri2018,
    AUTHOR = {Magri, Franco},
     TITLE = {Haantjes manifolds with symmetry},
   JOURNAL = {Teoret. Mat. Fiz.},
  FJOURNAL = {Teoreticheskaya i Matematicheskaya Fizika},
    VOLUME = {196},
      YEAR = {2018},
    NUMBER = {2},
     PAGES = {313--327},
      ISSN = {0564-6162,2305-3135},
   MRCLASS = {53Z05},
  MRNUMBER = {3833559},
       DOI = {10.4213/tmf9499},
       URL = {https://doi.org/10.4213/tmf9499},
}

@article {Akpan3,
    AUTHOR = {Akpan, Dinmukhammed},
     TITLE = {Singularities of two-dimensional {N}ijenhuis operators},
   JOURNAL = {Eur. J. Math.},
  FJOURNAL = {European Journal of Mathematics},
    VOLUME = {8},
      YEAR = {2022},
    NUMBER = {4},
     PAGES = {1328--1340},
      ISSN = {2199-675X,2199-6768},
   MRCLASS = {58K20 (53B20 53B30)},
  MRNUMBER = {4510607},
MRREVIEWER = {Vladimir\ S.\ Matveev},
    
       URL = {https://doi.org/10.1007/s40879-022-00585-x},
}

@article {Akpan2,
    AUTHOR = {Akpan, Dinmukhammed Zh.},
     TITLE = {Almost differentially nondegenerate {N}ijenhuis operators},
   JOURNAL = {Russ. J. Math. Phys.},
  FJOURNAL = {Russian Journal of Mathematical Physics},
    VOLUME = {29},
      YEAR = {2022},
    NUMBER = {4},
     PAGES = {413--416},
      ISSN = {1061-9208,1555-6638},
   MRCLASS = {53C15 (35Q83)},
  MRNUMBER = {4539410},
MRREVIEWER = {Vladimir\ S.\ Matveev},
       DOI = {10.1134/s106192082204001x},
       URL = {https://doi.org/10.1134/s106192082204001x},
}

@article {Akpan1,
    AUTHOR = {Akpan, Dinmukhammed  and Oshemkov, Andrey},
     TITLE = {Elementary differential singularities of three-dimensional
              {N}ijenhuis operators},
   JOURNAL = {Russ. J. Math. Phys.},
  FJOURNAL = {Russian Journal of Mathematical Physics},
    VOLUME = {30},
      YEAR = {2023},
    NUMBER = {4},
     PAGES = {425--431},
      ISSN = {1061-9208,1555-6638},
   MRCLASS = {53A55},
  MRNUMBER = {4683119},
MRREVIEWER = {Jan\ Kurek},
       DOI = {10.1134/S1061920823040015},
       URL = {https://doi.org/10.1134/S1061920823040015},
}

@article {nijapp5,
    AUTHOR = {Bolsinov, Alexey V. and Konyaev, Andrey Yu. and Matveev,
              Vladimir S.},
     TITLE = {Applications of {N}ijenhuis geometry {V}: geodesic equivalence
              and finite-dimensional reductions of integrable quasilinear
              systems},
   JOURNAL = {J. Nonlinear Sci.},
  FJOURNAL = {Journal of Nonlinear Science},
    VOLUME = {34},
      YEAR = {2024},
    NUMBER = {2},
     PAGES = {Paper No. 33, 18},
      ISSN = {0938-8974,1432-1467},
   MRCLASS = {37K06 (37K10 37K25 37K50)},
  MRNUMBER = {4696318},
       DOI = {10.1007/s00332-023-10008-0},
       URL = {https://doi.org/10.1007/s00332-023-10008-0},
}

@incollection{BKMmatrix,
  author    = {Bolsinov, Alexey V. and Konyaev, Andrey Yu. and Matveev, Vladimir S.},
  title     = {Research problems on relations between {Nijenhuis} geometry and integrable systems},
  booktitle = {2024 {MATRIX} Annals, Part I},
  editor    = {Wood, David R. and Etheridge, Alison M. and de Gier, Jan and Joshi, Nalini},
  series    = {MATRIX Book Series},
  volume    = {7},
  pages     = {229--268},
  publisher = {Springer},
  address   = {Cham},
  year      = {2026},
  doi       = {10.1007/978-3-032-16202-1_8}
}

@incollection{BKMmoser,
  author    = {Konyaev, Andrey Yu. and Matveev, Vladimir S.},
  title     = {Thoughts about potentials with finite-band spectrum and finite-dimensional reductions of integrable systems},
  booktitle = {2024 {MATRIX} Annals, Part I},
  editor    = {Wood, David R. and Etheridge, Alison M. and de Gier, Jan and Joshi, Nalini},
  series    = {MATRIX Book Series},
  volume    = {7},
  pages     = {201--228},
  publisher = {Springer},
  address   = {Cham},
  year      = {2026},
  doi       = {10.1007/978-3-032-16202-1_7}
}

@article{Lax1975,
  author  = {Lax, Peter D.},
  title   = {Periodic solutions of the {KdV} equation},
  journal = {Comm. Pure Appl. Math.},
  volume  = {28},
  number  = {1},
  year    = {1975},
  pages   = {141--188},
  doi     = {10.1002/cpa.3160280105}
}

@article{McKeanVanMoerbeke1975,
  author  = {McKean, H. P. and van Moerbeke, P.},
  title   = {The spectrum of {Hill}'s equation},
  journal = {Invent. Math.},
  volume  = {30},
  number  = {3},
  year    = {1975},
  pages   = {217--274},
  doi     = {10.1007/BF01425567}
}

@article{Neumann1859,
  author  = {Neumann, C.},
  title   = {De problemate quodam mechanico, quod ad primam integralium ultraellipticorum classem revocatur},
  journal = {J. Reine Angew. Math.},
  volume  = {56},
  year    = {1859},
  pages   = {46--63}
}

@incollection{Uhlenbeck1982,
  author    = {Uhlenbeck, Karen K.},
  title     = {Equivariant harmonic maps into spheres},
  booktitle = {Harmonic Maps (New Orleans, La., 1980)},
  series    = {Lecture Notes in Math.},
  volume    = {949},
  publisher = {Springer},
  address   = {Berlin--New York},
  year      = {1982},
  pages     = {146--158}
}

@incollection{Moser1980Chern,
  author    = {Moser, Jürgen},
  title     = {Geometry of quadrics and spectral theory},
  booktitle = {The Chern Symposium 1979 (Berkeley, Calif., 1979)},
  publisher = {Springer},
  address   = {New York--Berlin},
  year      = {1980},
  pages     = {147--188},
  doi       = {10.1007/978-1-4613-8109-9_7}
}

@article{Carter1977,
  author  = {Carter, Brandon},
  title   = {Killing tensor quantum numbers and conserved currents in curved space},
  journal = {Phys. Rev. D},
  volume  = {16},
  number  = {12},
  year    = {1977},
  pages   = {3395--3414},
  doi     = {10.1103/PhysRevD.16.3395}
}

@article{DuvalValent2005,
  author  = {Duval, C. and Valent, G.},
  title   = {Quantum integrability of quadratic {Killing} tensors},
  journal = {J. Math. Phys.},
  volume  = {46},
  number  = {5},
  year    = {2005},
  pages   = {053516},
  doi     = {10.1063/1.1899986}
}

@article{BLZ1996,
  author  = {Bazhanov, Vladimir V. and Lukyanov, Sergei L. and Zamolodchikov, Alexander B.},
  title   = {Integrable structure of conformal field theory, quantum {KdV} theory and thermodynamic {Bethe} ansatz},
  journal = {Comm. Math. Phys.},
  volume  = {177},
  number  = {2},
  year    = {1996},
  pages   = {381--398},
  doi     = {10.1007/BF02101898}
}

@article{TempestaTondo2022higher,
  author  = {Tempesta, Piergiulio and Tondo, Giorgio},
  title   = {Higher {Haantjes} brackets and integrability},
  journal = {Comm. Math. Phys.},
  volume  = {389},
  number  = {3},
  year    = {2022},
  pages   = {1647--1671},
  doi     = {10.1007/s00220-021-04233-5}
}

@article{Magri2014,
  author  = {Magri, Franco},
  title   = {{Haantjes} manifolds},
  journal = {J. Phys.: Conf. Ser.},
  volume  = {482},
  number  = {1},
  year    = {2014},
  pages   = {012028},
  doi     = {10.1088/1742-6596/482/1/012028}
}

@article{Zoll1903,
  author  = {Zoll, Otto},
  title   = {Ueber {Fl\"achen} mit {Scharen} geschlossener geod\"atischer {Linien}},
  journal = {Math. Ann.},
  volume  = {57},
  number  = {1},
  year    = {1903},
  pages   = {108--133},
  doi     = {10.1007/BF01449019}
}

@article{KressSchoebelVollmer2023,
  author  = {Kress, Jonathan and Sch\"obel, Konrad and Vollmer, Andreas},
  title   = {An algebraic geometric foundation for a classification of second-order superintegrable systems in arbitrary dimension},
  journal = {J. Geom. Anal.},
  volume  = {33},
  year    = {2023},
  pages   = {Paper No. 360},
  doi     = {10.1007/s12220-023-01413-8}
}

@article{KressSchoebelVollmer2024,
  author  = {Kress, Jonathan and Sch\"obel, Konrad and Vollmer, Andreas},
  title   = {Algebraic conditions for conformal superintegrability in arbitrary dimension},
  journal = {Comm. Math. Phys.},
  volume  = {405},
  year    = {2024},
  pages   = {Paper No. 92},
  doi     = {10.1007/s00220-023-04872-w}
}

@article{Takeuchi1983,
  author  = {Takeuchi, Masaru},
  title   = {Killing tensor fields on spaces of constant curvature},
  journal = {Tsukuba J. Math.},
  volume  = {7},
  number  = {2},
  year    = {1983},
  pages   = {233--255}
}

@article{Thompson1986,
  author  = {Thompson, Gerard},
  title   = {Killing tensors in spaces of constant curvature},
  journal = {J. Math. Phys.},
  volume  = {27},
  number  = {11},
  year    = {1986},
  pages   = {2693--2699},
  doi     = {10.1063/1.527288}
}

@article{HertlingManin1999,
  author  = {Hertling, Claus and Manin, Yuri},
  title   = {Weak {Frobenius} manifolds},
  journal = {Internat. Math. Res. Notices},
  year    = {1999},
  number  = {6},
  pages   = {277--286},
  doi     = {10.1155/S1073792899000148}
}

@article{Tsarev1991,
  author  = {Tsarev, Sergei  P.},
  title   = {The geometry of {Hamiltonian} systems of hydrodynamic type. {The} generalized hodograph method},
  journal = {Math. USSR-Izv.},
  volume  = {37},
  number  = {2},
  year    = {1991},
  pages   = {397--419},
  doi     = {10.1070/IM1991v037n02ABEH002069},
  note    = {Russian original: Izv. Akad. Nauk SSSR Ser. Mat. 54 (1990), no. 5, 1048--1068}
}

@article{Vinberg1963,
  author  = {Vinberg, Ernst  B.},
  title   = {The theory of homogeneous convex cones},
  journal = {Trans. Moscow Math. Soc.},
  volume  = {1963},
  year    = {1963},
  pages   = {340--403},
  note    = {English translation published in 1965. Russian original: Trudy Moskov. Mat. Obshch. 12 (1963), 303--358}
}

@article{Burde2006,
  author  = {Burde, Dietrich},
  title   = {Left-symmetric algebras, or pre-{Lie} algebras in geometry and physics},
  journal = {Cent. Eur. J. Math.},
  volume  = {4},
  number  = {3},
  year    = {2006},
  pages   = {323--357},
  doi     = {10.2478/s11533-006-0014-9}
}

@article{Bai2009,
  author  = {Bai, Chengming},
  title   = {Bijective 1-cocycles and classification of 3-dimensional left-symmetric algebras},
  journal = {Comm. Algebra},
  volume  = {37},
  number  = {3},
  year    = {2009},
  pages   = {1016--1057},
  doi     = {10.1080/00927870802279030}
}

@misc{Winterhalder1997,
  author        = {Winterhalder, Axel},
  title         = {Linear {Nijenhuis}-tensors and the construction of integrable systems},
  year          = {1997},
  eprint        = {physics/9709008},
  archivePrefix = {arXiv},
  note          = {Freiburg Preprint THEP 97/16}
}

@book{LibermannMarle1987,
  author    = {Libermann, Paulette and Marle, Charles-Michel},
  title     = {Symplectic Geometry and Analytical Mechanics},
  series    = {Mathematics and its Applications},
  volume    = {35},
  publisher = {D. Reidel},
  address   = {Dordrecht},
  year      = {1987},
  doi       = {10.1007/978-94-009-3807-6}
}

@article{Magri1978,
  author  = {Magri, Franco},
  title   = {A simple model of the integrable {Hamiltonian} equation},
  journal = {J. Math. Phys.},
  volume  = {19},
  number  = {5},
  year    = {1978},
  pages   = {1156--1162},
  doi     = {10.1063/1.523777}
}

@article{Turiel1994,
  author  = {Turiel, Francisco-Javier},
  title   = {Classification locale simultan\'ee de deux formes symplectiques compatibles},
  journal = {Manuscripta Math.},
  volume  = {82},
  number  = {3--4},
  year    = {1994},
  pages   = {349--362},
  doi     = {10.1007/BF02567706}
}

@article{Eliasson1990,
  author  = {Eliasson, L. H.},
  title   = {Normal forms for {Hamiltonian} systems with {Poisson} commuting integrals---elliptic case},
  journal = {Comment. Math. Helv.},
  volume  = {65},
  number  = {1},
  year    = {1990},
  pages   = {4--35},
  doi     = {10.1007/BF02566590}
}

@book{BolsinovFomenko2004,
  author    = {Bolsinov, Alexey V. and Fomenko, Anatoly T.},
  title     = {Integrable {Hamiltonian} Systems: Geometry, Topology, Classification},
  publisher = {Chapman \& Hall/CRC},
  address   = {Boca Raton, FL},
  year      = {2004},
  isbn      = {978-0-415-29805-6},
  doi       = {10.1201/9780203643426}
}

@article{Duistermaat1980,
  author  = {Duistermaat, J. J.},
  title   = {On global action-angle coordinates},
  journal = {Comm. Pure Appl. Math.},
  volume  = {33},
  number  = {6},
  year    = {1980},
  pages   = {687--706},
  doi     = {10.1002/cpa.3160330602}
}

@incollection{Kruskal1975,
  author    = {Kruskal, Martin D.},
  title     = {Nonlinear wave equations},
  booktitle = {Dynamical Systems, Theory and Applications (Battelle Seattle 1974 Rencontres)},
  editor    = {Moser, J.},
  series    = {Lecture Notes in Phys.},
  volume    = {38},
  publisher = {Springer},
  address   = {Berlin--Heidelberg},
  year      = {1975},
  pages     = {310--354},
  doi       = {10.1007/3-540-07171-7_9}
}

@article{Kaup1975,
  author  = {Kaup, D. J.},
  title   = {A higher-order water-wave equation and the method for solving it},
  journal = {Progr. Theoret. Phys.},
  volume  = {54},
  number  = {2},
  year    = {1975},
  pages   = {396--408},
  doi     = {10.1143/PTP.54.396}
}

@article{FuchssteinerFokas1981,
  author  = {Fuchssteiner, B. and Fokas, A. S.},
  title   = {Symplectic structures, their {B\"acklund} transformations and hereditary symmetries},
  journal = {Phys. D},
  volume  = {4},
  number  = {1},
  year    = {1981},
  pages   = {47--66},
  doi     = {10.1016/0167-2789(81)90004-X}
}

@article{ChenLiuZhang2006,
  author  = {Chen, Ming and Liu, Si-Qi and Zhang, Youjin},
  title   = {A two-component generalization of the {Camassa--Holm} equation and its solutions},
  journal = {Lett. Math. Phys.},
  volume  = {75},
  number  = {1},
  year    = {2006},
  pages   = {1--15},
  doi     = {10.1007/s11005-005-0041-7}
}

@article{Benenti1992,
  author  = {Benenti, Sergio},
  title   = {Inertia tensors and {St\"ackel} systems in the {Euclidean} spaces},
  journal = {Rend. Sem. Mat. Univ. Politec. Torino},
  volume  = {50},
  number  = {4},
  year    = {1992},
  pages   = {315--341}
}

@article{Benenti2005,
  author  = {Benenti, Sergio},
  title   = {Special symmetric two-tensors, equivalent dynamical systems, cofactor and bi-cofactor systems},
  journal = {Acta Appl. Math.},
  volume  = {87},
  number  = {1--3},
  year    = {2005},
  pages   = {33--91},
  doi     = {10.1007/s10440-005-1138-9}
}

@article{GGKM1967,
  author  = {Gardner, Clifford S. and Greene, John M. and Kruskal, Martin D. and Miura, Robert M.},
  title   = {Method for solving the {Korteweg-deVries} equation},
  journal = {Phys. Rev. Lett.},
  volume  = {19},
  number  = {19},
  year    = {1967},
  pages   = {1095--1097},
  doi     = {10.1103/PhysRevLett.19.1095}
}

@article{Hirota1971,
  author  = {Hirota, Ryogo},
  title   = {Exact solution of the {Korteweg--de Vries} equation for multiple collisions of solitons},
  journal = {Phys. Rev. Lett.},
  volume  = {27},
  number  = {18},
  year    = {1971},
  pages   = {1192--1194},
  doi     = {10.1103/PhysRevLett.27.1192}
}
\end{document}